\documentclass[12pt,letterpaper,openany]{book} 
\usepackage[utf8]{inputenc}
\usepackage{serina}
\usepackage{bbm}
\usepackage[backend=biber,style=alphabetic]{biblatex}
\usepackage{hyperref}
\allowdisplaybreaks

\newcommand{\nc}{\newcommand}
\nc{\on}{\operatorname}
\nc{\mbb}{\mathbb}
\nc{\mfk}{\mathfrak}
\nc{\mcal}{\mathcal}
\nc{\BR}{\mbb{R}}
\nc{\BN}{\mbb{N}}
\nc{\BC}{\mbb{C}}
\nc{\BZ}{\mbb{Z}}
\DeclareMathOperator{\Hom}{Hom}
\DeclareMathOperator{\coker}{coker}
\DeclareMathOperator{\End}{End}
\DeclareMathOperator{\id}{id}
\DeclareMathOperator{\Ob}{Ob}
\DeclareMathOperator{\ev}{ev}
\DeclareMathOperator{\coev}{coev}
\DeclareMathOperator{\Rep}{Rep}
\DeclareMathOperator{\Fun}{Fun}
\DeclareMathOperator{\Tr}{Tr}
\DeclareMathOperator{\Ind}{Ind}
\DeclareMathOperator{\colim}{colim}
\DeclareMathOperator{\ch}{ch}
\newcommand{\Los}{\L{}o\'{s}}
\usetikzlibrary{decorations.pathreplacing}

\title{An Introduction to Deligne Categories}
\author{Serina Hu; Advisors: Professor Michael Hopkins, Professor Pavel Etingof (MIT)}
\date{March 2021}
\newlength{\drop}

\newcommand*{\titleTMB}{\begingroup
  \drop=0.1\textheight
  \centering
  \settowidth{\unitlength}{\LARGE An Introduction to Deligne Categories}
  \vspace*{\baselineskip}
  {\large\scshape Serina Hu}\\[\baselineskip]
  \rule{\unitlength}{1.6pt}\vspace*{-\baselineskip}\vspace*{2pt}
  \rule{\unitlength}{0.4pt}\\[\baselineskip]
  {\LARGE An Introduction to Deligne Categories}
  \rule{\unitlength}{0.4pt}\vspace*{-\baselineskip}\vspace{3.2pt}
  \rule{\unitlength}{1.6pt}\\[\baselineskip]
  {\large\scshape Advised by Professor Michael Hopkins}\\
  {\large\scshape \&}\\
    {\large \scshape Professor Pavel Etingof (MIT)}\\
  \vfill
  \includegraphics[width=0.3\textwidth]{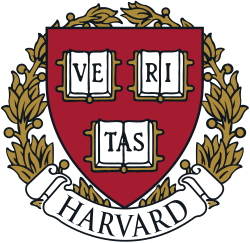}
  \vfill
  {\large An undergraduate thesis submitted in partial fulfillment of the honors requirement for the degree of Bachelor of Arts in Mathematics}\\[\baselineskip]
  {\small\scshape March 2021}\par
  \vspace*{\drop}
\endgroup}

\begin{document}
\begin{titlepage}
  \titleTMB{}
\end{titlepage}

\tableofcontents
\chapter*{Introduction}
\markboth{Introduction}{Introduction}
One way to study the representations of a group (or some other algebraic object) is to consider the category of its representations, $\Rep(G)$. This category is endowed with the structure of a symmetric tensor category: roughly, it is equipped with direct sums, a symmetric tensor product, and duals. In such a category, we can define a categorical notion of ``dimension,'' which agrees with our typical notion of dimension on vector spaces. This leads us to two (related) questions: do all symmetric tensor categories arise as (super)group representation categories, and do all symmetric tensor categories have only objects of integer dimension? A theorem of Deligne (\cite{deligne_tensor}) says that in characteristic zero, a symmetric category is equivalent to a representation category if and only if it satisfies some constraints (known as sub-exponential growth) on the dimensions of tensor powers of objects. Moreover, Deligne and Milne constructed examples of categories that don't satisfy these constraints and have objects of non-integer dimension, first in \cite{deligne_milne}. These are now known as the Deligne categories and notated $\Rep(S_t)$, $\Rep(GL_t)$, $\Rep(O_t)$, and $\Rep(Sp_t)$, generalizations of the classical representation categories of $S_n$, $GL_n$, $O_n$, and $Sp_{2n}$ to complex values of $n$. Although they are called $\Rep(G_t)$, none of these categories arise as the representation category of a (super)group, so they are examples of ``exotic'' symmetric tensor categories.

What are some properties of these categories, and how do they relate to the classical representation categories? For one, Deligne categories have nice universal properties: roughly, $\Rep(S_t)$ is universal in categories generated by a symmetric (i.e. commutative) Frobenius algebra with nondegenerate trace; $\Rep(GL_t)$ is universal in categories generated by an object and its dual; and $\Rep(O_t)$ (resp. $\Rep(Sp_t)$) is universal in categories generated by an object with a symmetric (resp. skew-symmetric) isomorphism with its dual (all the objects mentioned have dimension $t$). Deligne categories also provide a way to study the (walled) Brauer algebras $B_{r, s}(t)$ and partition algebras $P_n(t)$ for non-integer parameters $t$. As we will see, these algebras and their subalgebras form the endomorphism algebras of Deligne categories, so studying properties of said categories amount to simultaneously studying these algebras for all values of $(r, s)$ or $n$. Finally, Deligne categories form the basis of ``representation theory in complex rank,'' laid out in \cite{etingof_cr1} and \cite{etingof_cr2}, which seeks to transfer problems of classical representation theory to Deligne categories and can lead to new realizations of classical constructions, e.g. in \cite{etingof_ddca}.

This work hopes to be an introduction to Deligne categories for someone familiar with classical representation theory and some category theory. In the first chapter, we motivate and define (symmetric) tensor categories, construct the Deligne categories, and prove their universal properties. In the second chapter, we define ultraproducts and use them to construct the Deligne categories in a different way; in particular, this construction works in positive characteristic. In the third chapter, we connect Deligne categories to the classical categories by discussing simple objects and explaining what happens in the $t = n$ case. In the last chapter, along with recommendations for further and related reading, we give an example of an investigation of Harish-Chandra bimodules in representation theory in complex rank, which is a joint work with Alexandra Utiralova (\cite{utiralova2022harishchandra}).

\chapter{Constructing the Deligne categories}
\label{chap:interp_constructions}

In this chapter, our goal is to construct the Deligne categories as interpolations of the standard categories and prove their universal properties. We assume familiarity with the classical representation theory of the symmetric group and Lie algebras of types $A, B, C, D$, and refer the reader to \cite{FH} otherwise. In the first section, we provide an overview of symmetric tensor categories. Most of the category theory is a formalization of concepts from ring theory and defines terminology we use in the rest of the work, so we encourage skimming it to establish familiarity with the notation. However, we direct the reader to pay special attention to the notion of the Karoubian envelope, which plays a key role in the construction. For proofs of the results, see the referenced literature. In the second and third sections, we construct the Deligne categories and prove their universal properties.

\section{Category theory preliminaries}

In this section, we define tensor categories and describe some of the properties we will be discussing later on. For proof we redirect readers to \cite{EGNO}. We assume that the reader is familiar with the basic theory of categories and functors, see for example \cite{MacLane}. Deligne categories are examples of symmetric Karoubian tensor categories; roughly, a tensor category's objects are abstractions of finite-dimensional vector spaces and can be added, multiplied, and dualized; a Karoubian category is closed under taking images of idempotents; and a symmetric tensor category has commutative tensor product. Formally,
\begin{Definition}{}{}
  A \textbf{tensor category} is a category $\mathcal{C}$ that is locally finite, $F$-linear, abelian, rigid, and monoidal with bilinear tensor product and $\End_{\mathcal{C}} \mathbf{1} = F$. All the tensor categories we are considering will also be symmetric.
\end{Definition}
We now define each of the components of the definition above.

Abelian categories are ubiquitous throughout algebra as the most general setting in which all kernels and cokernels of morphisms exist, and thus we can work with exact sequences and perform homological algebra.
\begin{Definition}{}{}
  An \textbf{abelian category} is a category $\mathcal{C}$ that satisfies the following properties:
  \begin{enumerate}
  \item Every Hom-set has the structure of an abelian group, and composition of morphisms is bilinear.
  \item All finite products exist. In fact, coproducts and products coincide, so we call such products direct sums and we have a bifunctor $\oplus: \mathcal{C} \times \mathcal{C} \to \mathcal{C}$.
  \item There exists a zero object $0$ with $\Hom(0, 0) = 0$.
  \item All kernels and cokernels exist.
  \item Every monomorphism is a kernel and every epimorphism is a cokernel.
  \end{enumerate}
  The last condition has a number of equivalent formulations; in particular, that every morphism $f: X \to Y$ has a canonical decomposition
  \begin{equation*}
    K \xrightarrow{k} X \xrightarrow{i} I \xrightarrow{j} Y \xrightarrow{c} C
  \end{equation*}
  with $K = \ker f$, $C = \coker f$, $I = \coker k = \ker c$. We call $I$ the image of $f$.
\end{Definition}
\begin{Definition}{}{}
  A category satisfying only the first three axioms is an \textbf{additive category}.
\end{Definition}

Now it also makes sense to define quotients and subobjects: $X$ is a subobject of $Y$ if there is a monomorphism $\iota : X \to Y$, and $Z$ is a quotient of $Y$ if there is an epimorphism $\pi: Y \to Z$.

These notions are familiar, and for good reason:
\begin{Theorem}{Freyd-Mitchell}{}
  Every (small) abelian category is equivalent as an additive category to a full subcategory to the category of left modules over a ring.
\end{Theorem}
Since we can think of abelian categories as categories of modules, we can define an analog of composition length for objects in abelian categories.
\begin{Definition}{}{}
  A \textbf{Jordan-Holder series} of an object $X$ is a filtration $0 = X_0 \subset X_1 \subset \cdots \subset X_n = X$ so that each $X_i / X_{i - 1}$ is simple. We say $X$ has \textbf{length} $n$, which is a well-defined notion, since it turns out that all Jordan-Holder series for $X$ have the same length. The length of an object is analogous to the composition length of a module.
\end{Definition}

\begin{Definition}{}{}
  A category $\mathcal{C}$ is an \textbf{$F$-linear category} for a field $F$ if the Hom-sets are all vector spaces over $F$ and composition of morphisms is $F$-bilinear.
\end{Definition}

When we talk about functors between ($F$-linear) abelian categories, we want them to respect the Hom-set structure, so:
\begin{Definition}{}{}
  A functor $F: \mathcal{C} \to \mathcal{D}$ between two additive categories is an \textbf{additive ($F$-linear) functor} if $\Hom_{\mathcal{C}}(X, Y) \to \Hom_{\mathcal{D}}(F(X), F(Y))$ is a homomorphism of abelian groups ($F$-vector spaces). In both cases, $F$ respects direct sums, so we have a natural isomorphism $F(X) \oplus F(Y) \to F(X \oplus Y)$.
\end{Definition}

\begin{Definition}{}{}
  A $F$-linear abelian category is \textbf{locally finite} or \textbf{Artinian} if
  \begin{enumerate}
  \item All Hom-sets are finite-dimensional
  \item All objects have finite length.
  \end{enumerate}
  In locally finite categories, many statements about arbitrary objects can be reduced to statements about simple objects.
\end{Definition}

Monoidal categories, as the name suggests, give us a unit object and the ability to ``multiply'' objects, and their axioms ``categorify'' associativity and properties of the unit in a monoid.

\begin{Definition}{}{}
  A \textbf{monoidal category} is a category $\mathcal{C}$ along with a bifunctor $\otimes : \mathcal{C} \times \mathcal{C} \to \mathcal{C}$, an associativity constraint (natural isomorphism) $a: (- \otimes -) \otimes - \xrightarrow{\sim} - \otimes (- \otimes -)$, an object $\mathbf{1}$ of $\mathcal{C}$, and an isomorphism $\iota: \mathbf{1} \otimes \mathbf{1} \to \mathbf{1}$ that satisfies the following axioms:
  \begin{enumerate}
  \item The pentagon axiom, ensuring that different parenthesizations of products cohere:
    \begin{center}
    \begin{tikzcd}
      & ((W \otimes X) \otimes Y) \otimes Z \arrow[dl, "a_{W, X, Y} \otimes \id_Z"] \arrow[dr, "a_{W \otimes X, Y, Z}"] & \\
      (W \otimes (X \otimes Y)) \otimes Z \arrow[d, "a_{W, X \otimes Y, Z}"] & & (W \otimes X) \otimes (Y \otimes Z) \arrow[d, "a_{W, X, Y \otimes Z}"] \\
      W \otimes ((X \otimes Y) \otimes Z) \arrow[rr, "\id_W \otimes a_{X, Y, Z}"] & & W \otimes (X \otimes (Y \otimes Z))
    \end{tikzcd}
  \end{center}
\item the unit axiom, that left and right multiplication $L_1 : X \mapsto \mathbf{1} \otimes X$ and $R_1 : X \mapsto X \otimes \mathbf{1}$, are equivalences $\mathcal{C} \to \mathcal{C}$.
  \end{enumerate}
  It turns out that any two choices of $\iota$ are equivalent, so we will denote a monoidal category by the quadruple $(\mathcal{C}, \otimes, \mathbf{1}, a)$, and in most of our examples we'll omit the latter three items and just refer to the category as $\mathcal{C}$.
\end{Definition}

We can also define natural isomorphisms $l : \mathbf{1} \otimes - \to -$ and $r : - \otimes \mathbf{1} \to -$ by requiring
\begin{align*}
  L_1(l_X) &: \mathbf{1} \otimes (\mathbf{1} \otimes X) \xrightarrow{a^{-1}_{\mathbf{1}, \mathbf{1}, X}} (\mathbf{1} \otimes \mathbf{1}) \otimes X \xrightarrow{\iota \otimes \id_X} \mathbf{1} \otimes X \\
  R_1(r_X) &: (X \otimes \mathbf{1}) \otimes \mathbf{1} \xrightarrow{a_{X, \mathbf{1}, \mathbf{1}}} X \otimes (\mathbf{1} \otimes \mathbf{1}) \xrightarrow{\id_X \otimes \iota} X \otimes \mathbf{1}.
\end{align*}
The above two axioms are enough to guarantee that the multiple ways of simplifying each of $X \otimes Y \otimes \mathbf{1}$, $X \otimes \mathbf{1} \otimes Y$, and $\mathbf{1} \otimes X \otimes Y$ to $X \otimes Y$ are equivalent (the triangle axioms).

They are also enough to guarantee that parenthesizations of arbitrary numbers of objects are canonically isomorphic, so we will not have to worry about parenthesizing products in the future.
\begin{Theorem}{Coherence Theorem}{}
  Let $X_1, \dots, X_n \in \Ob(\mathcal{C})$. Let $P_1$ and $P_2$ be any (ordered) parenthesized products of $X_1, \dots, X_n$ with arbitrary insertions of $\mathbf{1}$, and suppose we have two isomorphisms $\varphi, \psi: P_1 \to P_2$ consisting of asssociativity, unit, and identity morphisms and their inverses. Then $\varphi = \psi$.
\end{Theorem}
Thus ``multiplication'' works as we would expect.

Monoidal functors, analogous to maps between monoids, respect the tensor product of objects.
\begin{Definition}{}{}
  Given two monoidal categories $(\mathcal{C}, \otimes, \mathbf{1}, a, \iota)$  and $(\mathcal{C}', \otimes', \mathbf{1}', a', \iota,)$, a \textbf{monoidal functor} is a pair $(F, J)$ where $F: \mathcal{C} \to \mathcal{C}'$ is a functor and $J$ is a natural isomorphism $F(-) \otimes' F(-) \xrightarrow{sim} F(- \otimes -)$, so that
  \begin{enumerate}
  \item $F(\mathbf{1}) \cong \mathbf{1}'$,
  \item $F$ respects associativity, so the following diagram commutes for all $X, Y, Z \in \Ob(\mathcal{C})$:
    \begin{center}
    \begin{tikzcd}[column sep=huge]
      (F(X) \otimes' F(Y)) \otimes' F(Z) \arrow[r, "a'_{F(X), F(Y), F(Z)}"] \arrow[d, "J_{X, Y} \otimes' \id_{F(Z)}"] & F(X) \otimes' (F(Y) \otimes' F(Z)) \arrow[d, "\id_{F(X) \otimes' J_{Y, Z}}"] \\
      F(X \otimes Y) \otimes' F(Z) \arrow[d, "J_{X \otimes Y, Z}"] & F(X) \otimes' F(Y \otimes Z) \arrow[d, "J_{X, Y \otimes Z}"] \\
      F((X \otimes Y) \otimes Z) \arrow[r, "F(a_{X, Y, Z})"] & F(X \otimes (Y \otimes Z))
    \end{tikzcd}
  \end{center}
\end{enumerate}
  It's possible that there are many different choices of $J$ that can turn an ordinary functor $F$ into a monoidal functor. In addition, the axioms imply that $F(\mathbf{1})$ is actually canonically isomorphic to $\mathbf{1}$.
\end{Definition}

The final property in the definition of a tensor category is that the category is rigid: its objects have both left and right duals, which have properties similar to that of the dual of a finite-dimensonial vector space.
\begin{Definition}{Left and right duals}{}
  Let $X \in \Ob(\mathcal{C})$. $X^* \in \Ob(\mathcal{C})$ is a \textbf{left dual} of $X$ if there are evaluation and coevaluation maps, $\ev_X: X^* \otimes X \to \mathbf{1}$ and $\coev_X: \mathbf{1} \to X \otimes X^*$ so that the compositions (omitting associativity isomorphisms)
  \begin{align*}
    X \xrightarrow{\coev_X \otimes \id_X} &X \otimes X^* \otimes X \xrightarrow{\id_X \otimes \ev_X} X \\
    X^* \xrightarrow{\id_{X^*} \otimes \coev_X} &X^* \otimes X \otimes X^* \xrightarrow{\ev_X \otimes \id_{X^*}} X^*
  \end{align*}
  are isomorphisms. Likewise, $\prescript{*}{}{X}$ is a \textbf{right dual} of $X$ if there are evaluation and coevaluation maps $\ev'_X : X \otimes \prescript{*}{}{X} \to \mathbf{1}$ and $\coev'_X: \mathbf{1} \to \prescript{*}{}{X} \otimes X$ so that we have isomorphisms
  \begin{align*}
    X \xrightarrow{\id_X \otimes \coev'_X} &X \otimes \prescript{*}{}{X} \otimes X \xrightarrow{\ev'_X \otimes \id_X} X \\
    \prescript{*}{}{X} \xrightarrow{\coev'_X \otimes \id_{\prescript{*}{}{X}}} &\prescript{*}{}{X} \otimes X \otimes \prescript{*}{}{X} \xrightarrow{\id_{\prescript{*}{}{X}} \otimes \ev'_X} \prescript{*}{}{X}.
  \end{align*}

\end{Definition}
\begin{Definition}{}{}
  A category is \textbf{rigid} if every object admits both left and right duals.
\end{Definition}

As you might expect, duals are unique up to unique isomorphism, and a map $f: X \to Y$ induces maps $f^*: Y^* \to X^*$ and $\prescript{*}{}{f} : \prescript{*}{}{Y} \to \prescript{*}{}{X}$. $f^*$ satisfies the property that
\begin{equation}
  \label{eq:dual}
  \ev_X \circ (f^* \otimes \id_X) = \ev_Y \circ (\id_{Y^*} \otimes f), \qquad (\id_X \otimes f^*) \circ \coev_X = (f \otimes \id_{Y^*}) \circ \coev_Y
\end{equation}
and $\prescript{*}{}{f}$ satisfies an analogous property. We also have an isomorphism $\Hom(X, Y) \cong \Hom(\mathbf{1}, Y \otimes X^*)$ given by
\begin{equation}
  \label{eq:dual_hom}
  X \xrightarrow{f} Y \mapsto \mathbf{1} \xrightarrow{\coev_X} X \otimes X^* \xrightarrow{f \otimes \id_{X^*}} Y \otimes X^*.
\end{equation}
We will denote this composition by $\Xi(f)$.

Moreover, for a monoidal functor $F$, $F(X)^* \cong F(X^*)$ and $\prescript{*}{}{F(X)} \cong F(\prescript{*}{}{X})$. Moving up a level of abstraction, for two functors $F, F': \mathcal{C} \to \mathcal{D}$, where $\mathcal{C}$ is rigid monoidal and $\mathcal{D}$ is monoidal, every natural transformation $\eta: F \to F'$ is a natural isomorphism and $\eta_{X^*} = (\eta_X^{_1})^* = ((\eta_X)^*)^{-1}$.

Putting this all together, tensor categories can be thought of as categorifications of rings with some finiteness conditions. These appear all the time, as the notation of direct sum and tensor product indicate: finite-dimensional vector spaces, finite-dimensional representations of groups, and finite-dimensional representations of Lie algebras. We note that all of these have a symmetric tensor product: there is a natural commutativity isomorphism given by swapping components. In general, we have a family of isomorphisms $c_{X, Y} : X \otimes Y \to Y \otimes X$ with $c_{Y, X} \circ c_{X, Y} = \id_{X \otimes Y}$. This is known as the $c$ and must satisfy the following coherence conditions (omitting associativity isomorphisms):
\begin{equation*}
  \begin{tikzcd}
    X \otimes Y \otimes Z \arrow[dr, "c_{X, Y} \otimes \id_Z"] \arrow[rr, "c_{X, Y \otimes Z}"] & & Y \otimes Z \otimes X \\
    & Y \otimes X \otimes Z \arrow[ur, "\id_Y \otimes c_{X, Z}"] & 
  \end{tikzcd}
\end{equation*}
and
\begin{equation*}
  \begin{tikzcd}
    X \otimes Y \otimes Z \arrow[dr, "\id_X \otimes c_{Y, Z}"] \arrow[rr, "c_{X \otimes Y, Z}"] & & Z \otimes X \otimes Y \\
    & X \otimes Z \otimes Y \arrow[ur, "c_{X, Z} \otimes \id_Y"] &
  \end{tikzcd}.
\end{equation*}

Functors between symmetric tensor categories must also respect the braiding, so the following diagram must commute, where $c'$ is the braiding in the target category:
\begin{center}
\begin{tikzcd}[column sep=large]
  F(X) \otimes F(Y) \arrow[r, "c'_{F(X), F(Y)}"] \arrow[d, "J_{X, Y}"] & F(Y) \otimes F(X) \arrow[d, "J_{Y, X}"] \\
  F(X \otimes Y) \arrow[r, "F(c_{X, Y})"] & F(Y \otimes X)
\end{tikzcd}
\end{center}
Moreover, symmetric tensor categories have a notion of trace and dimension, analogous to the vector space notions.
\begin{Definition}{}{}
  Let $\mathcal{C}$ be a symmetric monoidal category and suppose $X$ has left and right duals. Then the \textbf{trace} of $f$, $\Tr(f)$, is
  \begin{equation*}
    \mathbf{1} \xrightarrow{\coev_X} X \otimes X^* \xrightarrow{f \otimes \id_{X^*}} X \otimes X^* \xrightarrow{c_{X, X^*}} X^* \otimes X \xrightarrow{\ev_X} \mathbf{1}.
  \end{equation*}
\end{Definition}
Since $\End(\mathbf{1}) = F$, the trace of a morphism is just a scalar.
\begin{Definition}{}{}
  The \textbf{dimension} of $X$, $\dim X$, is $\Tr(\id_X)$.
\end{Definition}
It is easy to see that if $\mathcal{C}$ is the category of $F$-vector spaces, this agrees with the usual vector space trace and dimension. Monoidal functors preserve dimension, a fact we will use when we think about the universal properties of the Deligne categories.

Finally, we briefly explain Karoubian categories and Karoubian envelopes. Karoubian categories are closed under images of idempotents (i.e. direct summands) and direct sums.
\begin{Definition}{}{karoubi}
  The \textbf{Karoubian envelope} (or \textbf{pseudo-abelian completion}) of an additive category $\mathcal{C}$ has the following:
  \begin{itemize}
  \item objects: pairs $(X, e)$ with $X \in \Ob(\mathcal{C})$ and $e \in \End(X)$ an idempotent, which we can interpret as ``the image of $e$''
  \item morphisms: $\Hom((X, e), (Y, f)) = f\Hom_{\mathcal{C}}(X, Y)e$
  \item composition: same as $\mathcal{C}$
  \end{itemize}
  In other words, every idempotent now has a kernel and cokernel. Moreover, for every idempotent $e: X \to X$, we can write $(X, \id) = \ker e \oplus \ker (1 - e) = (X, 1 - e) \oplus (X, e)$. Categories already satisfying this condition are \textbf{Karoubian categories}.
\end{Definition}
\begin{Remark}{}{}
  Recall that an idempotent $e$ is primitive if it cannot be written as a sum $e_1 + e_2$ where $e_1, e_2$ are nonzero idempotents and $e_1 e_2 = e_2 e_1 = 0$. Every object in a Karoubian envelope can written uniquely (up to permutation) as a finite direct sum of indecomposables, where each indecomposable corresponds to $(X, e)$ with $e$ a primitive idempotent of $X$.
\end{Remark}
\begin{Example}{}{}
  A example of a Karoubian envelope is $\Rep(S_n)$, which is the Karoubian envelope of the category whose objects are direct sums and tensor products of the standard $n$-dimensional representation and morphisms are maps between representations. The irreducible representations are then obtained as images of the Young symmetrizers, and all the (finite-dimensional) representations of $S_n$ can then be obtained as images of direct sums of Young symmetrizers.
\end{Example}

\section{Construction of $\Rep(S_t)$}
In this section, we will express some properties of $\Rep(S_n)$ in a way that is independent of the integrality of $n$; then, we can replace $n$ with an arbitrary complex number $t$ and construct the tensor category $\Rep(S_t)$. All our categories here will be $\mathbb{C}$-linear, and we will come back to the positive characteristic case in the next chapter. We combine the exposition in \cite{EGNO} and \cite{comes_ostrik}.

\subsection{Reframing $\Rep(S_n)$}

Let $V$ be the $n$-dimensional permutation representation of $S_n$ where $S_n$ acts by permuting the basis vectors.
Consider the category $\Rep(S_n)^{(0)}$ whose objects are $V^{\otimes m}$ for $m \ge 0$ and Hom-spaces are $\Hom(V^{\otimes l}, V^{\otimes m}) = \Hom_{S_n}(V^{\otimes l}, V^{\otimes m})$.
We know that all the irreducible representations of $S_n$ appear as summands in tensor powers of $V$.
So we can construct $\Rep(S_n)$ as the Karoubian envelope (\cref{def:karoubi}) of the additive completion of $\Rep(S_n)^{(0)}$.
This provides an idea for how to interpolate $\Rep(S_t)$: start with a category $\Rep(S_t)^{(0)}$ whose only objects are tensor powers of a single generating object (the analogue of $V$), define the Hom spaces, and take a completion.

Defining the Hom spaces requires us to write $\Hom_{S_n}(V^{\otimes l}, V^{\otimes m})$ without mentioning $S_n$.
In fact, this Hom space has a basis indexed by partitions of $l + m$, and we will write composition as a polynomial function of $n$. Then, replacing $n$ with $t$, we can define $\Hom_{\Rep(S_t)^{(0)}}(V^{\otimes l}, V^{\otimes m})$.
Let $X = \{1, \dots, n\}$. Then $V$ corresponds to $\Fun(X, \mathbb{C})$, the functions $X \to \mathbb{C}$. Moreover, $S_n$ acts by $\sigma(f(i)) = f(\sigma^{-1}(i))$; the standard basis of $V$ corresponds to the characteristic functions on $X$.
Notice that $V^* \cong V$, so $\Hom(V^{\otimes l}, V^{\otimes m}) \cong V^{\otimes l} \otimes_{S_n} V^{\otimes m} \cong \Fun(X^{l + m}, \mathbb{C}^{S_n})$, the $S_n$-invariant functions on $X^{l + m}$. Here, $S_n$ acts on a function by $\sigma(f(x_1, \dots, x_{l + m})) = f(\sigma^{-1}(x_1), \dots, \sigma^{-1}(x_{l + m}))$.

Notice from the above that if $f$ is $S_n$-invariant, then for every orbit $O$, $f(x)$ is the same for every $x = (x_1, \dots, x_{l + m}) \in O$. In other words, there is a basis for $\Fun(X^{l + m}, \mathbb{C})^{S_n}$ consisting of the characteristic functions $\delta_O$ corresponding to each orbit $O$ of $S_n$.
\begin{Definition}{}{}
  Each orbit $O$ is labeled by a partition $P$ of $l + m$, corresponding to which elements in the tuple $(x_1, \dots, x_{l + m})$ are equal to one another. We call these \textbf{equality patterns} and denote the orbit by $O(P)$. To emphasize that we are partitioning the sum $l + m$, in our examples we will partition the set $\{1, \dots, l, 1', \dots, m'\}$. For example, if $l = 2$ and $m = 3$, the tuple $(1, 1, 2, 5, 2)$ corresponds to the partition $([1, 2], [1', 3'], [2'])$.
\end{Definition}
When $n \ge l + m$, all partitions of $l + m$ have a corresponding orbit, and so $\Hom_{S_n} (V^{\otimes l}, V^{\otimes m})$ has a basis independent of $n$. We will later refer to this Hom-space as $\mathbb{C}P_{l, m}(n)$, the $\mathbb{C}$-vector space with basis $P_{l, m}$, the partitions of $l + m$. If $l = m$, we will abbreviate the Hom-space as $\mathbb{C}P_m(n)$. In Deligne categories, the integer parameter $n$ will be replaced by an arbitrary complex number $t$. In addition, we denote the characteristic functions by $\delta_P$ with $P$ a partition of $l + m$.\footnote{In the notation of \cite{comes_ostrik}, $\delta$ is the characteristic function for a ``perfect'' pattern; $e$, defined below, is notated by $f$ and a ``good'' pattern.}
\begin{Example}{}{}
  Here are some examples of associating $\delta_P$ with a map $V^{\otimes l} \to V^{\otimes m}$; we essentially ``roll back'' the identifications used above.
  \begin{enumerate}
  \item Let $l = 0$, $m = 5$, and $P = ([1', 2'], [3', 5'], [4'])$. Then the corresponding map $\mathbb{C} \to V^{\otimes 5}$ takes $1 \mapsto \sum_{i \ne j \ne k \in X} v_i \otimes v_i \otimes v_j \otimes v_k \otimes v_j$.
  \item Let $l = 4$, $m = 3$, and $P = ([1, 1'], [2, 4], [3], [2', 3'])$. Then the corresponding map $V^{\otimes 4} \to V^{\otimes 3}$ takes
    \begin{equation*}
      v_{i_1} \otimes v_{i_2} \otimes v_{i_3} \otimes v_{i_4} \mapsto \begin{cases} \sum_{j \in X} v_{i_1} \otimes v_j \otimes v_j & \text{if } i_2 = i_4, i_1 \ne i_2 \ne i_3 \ne j \\
        0 & \text{otherwise} \end{cases}.
    \end{equation*}
    Here $v_i$ is the $i$th basis vector of $V$.
  \end{enumerate}
\end{Example}
In order to define composition of morphisms, we will use a different basis of $\Hom_{S_n}(V^{\otimes l}, V^{\otimes m})$. We say that a partition $P$ \textbf{refines} $P'$ if $P'$ can be obtained by taking unions of some subsets in $P$ (for example, $([1, 2], [1', 3'], [2'])$ refines $([1, 2, 1', 3'], [2'])$). Our new basis consists of functions
\begin{equation*}
  e_P = \sum_{P' \ge P} \delta_{O(P)}.
\end{equation*}
In other words, $e_P$ is the characteristic function on tuples $(x_1, \dots, x_{l + m})$ so that $x_i = x_j$ if $i, j$ are in the same subset of $P$, but it can also happen that $x_i = x_j$ for $i, j$ not in the same subset. Since $e_P$ differs from $\delta_P$ by a triangular matrix with $1s$ on the diagonal, it is also a basis for $\Hom_{S_n}(V^{\otimes l}, V^{\otimes m})$.
\begin{Example}{}{correspondence}
  Adapting our examples above, essentially the element of $\Hom_{S_n}(V^{\otimes l}, V^{\otimes m})$ corresponding to $e_P$ has a less restrictive equality pattern than the one corresponding to $\delta_P$.
  \begin{enumerate}
  \item Let $l = 0$, $m = 5$, $P = ([1', 2'], [3', 5'], [4'])$. Then the corresponding map $\mathbb{C} \to V^{\otimes 5}$ takes $1 \mapsto \sum_{i, j, k \in X} v_i \otimes v_i \otimes v_j \otimes v_k \otimes v_j$.
  \item Let $l = 4$, $m = 3$, $P = ([1, 1'], [2, 4], [3], [2', 3'])$. Then the corresponding map $V^{\otimes 4} \to V^{\otimes 3}$ takes
    \begin{equation*}
      v_{i_1} \otimes v_{i_2} \otimes v_{i_3} \otimes v_{i_4} \mapsto \begin{cases} \sum_{j \in X} v_{i_1} \otimes v_j \otimes v_j & \text{if } i_2 = i_4 \\
        0 & \text{otherwise} \end{cases}.
    \end{equation*}
  \end{enumerate}
\end{Example}

Now we can actually define composition. Let $e_P: V^{\otimes l} \to V^{\otimes m}$ and $e_Q: V^{\otimes k} \to V^{\otimes l}$, and we want to say something about $e_P \circ e_Q$. $P$ corresponds to a partition of $\{1', \dots, l', 1'', \dots, m''\}$, enforcing equalities between variables $y_1, \dots, y_l, z_1, \dots, z_m \in X$. Likewise, $Q$ corresponds to a partition of $\{1, \dots, k, 1', \dots, l'\}$ and equalities between $x_1, \dots, x_k, y_1, \dots, y_l$. Now we count the possible equality patterns we could have between the $x_i$ and the $z_i$. Let $P * Q$ be the least restrictive equality pattern possible between $x_1, \dots, x_k, z_1, \dots, z_m$ that still satisfies $P$ and $Q$, corresponding to a partition of $\{1, \dots, k, 1'', \dots, m''\}$. This determines some of the $y_i$, but it leaves $N(P, Q)$ groups of $y_i$ to be equal to any element of $X$, so there are $n^{N(P, Q)}$ ways to set these other $y_i$ while satisfying equality patterns. Therefore,
\begin{Lemma}{}{}
  $e_P \circ e_Q = n^{N(P, Q)}e_{P * Q}$.
\end{Lemma}
This expresses composition in a way that we can replace $n$ with any complex number $t$. So we can finally define $\Rep(S_t)$.

\begin{Definition}{}{}
  First, we define $\Rep(S_t)^{(0)}$, whose only objects are analogues of tensor powers of $V$. The objects $[m]$ are labeled by nonnegative integers $m$; the unit object is $[0]$ and $[l] \otimes [m] = [l + m]$. $\Hom([l], [m])$ is an $\mathbb{C}$-vector space with basis $e_P$ labeled by partitions of $l + m$, and $e_P \circ e_Q = t^{N(P, Q)} e_{P * Q}$. The tensor product $e_P \otimes e_Q = e_{P \cup Q}$.
\end{Definition}

\begin{Definition}{}{}
  $\Rep(S_t)$ is the Karoubian envelope of the additive completion of $\Rep(S_t)^{(0)}$.
\end{Definition}

\subsection{Properties of $[m]$ and graphical notation for partitions}
The $e_P$ notation for a basis of $\Hom_{\Rep(S_t)^{(0)}}([l], [m])$ can also be represented graphically: we label $l$ points as $1, \dots, l$ and $m$ points as $1', \dots, m'$, and points are connected if and only if their labels are in the same part of the partition.
\begin{Example}{}{}
  \begin{enumerate}
  \item The identity $\id_m: [m] \to [m]$ is
    \begin{center}
      \begin{tikzpicture}
        \filldraw[black] (0,0) circle (2pt)
        (0, 1) circle (2pt)
        (3, 0) circle (2pt)
        (3, 1) circle (2pt);
        \draw (0,0) node[below=1pt] {1'} -- (0,1) node[above=1pt] {1};
        \draw (3,0) node[below=1pt] {m'} -- (3,1) node[above=1pt] {m};
        \draw[dotted] (1.5,0) node[below=1pt] {$\cdots$} -- (1.5,1) node[above=1pt] {$\cdots$};
      \end{tikzpicture}
    \end{center}
  \item The map $e_P: [3] \to [6]$ with $P = ([1, 3, 2'], [2, 4', 5'], [1'], [3', 6'])$ is
    \begin{center}
      \begin{tikzpicture}
        \foreach \x [count=\l] in {0,...,2} {
          \filldraw[black] (\x,1) circle (2pt) node[above=1pt] {\l};
        }
        \foreach \x [count=\l] in {0,...,5} {
          \filldraw[black] (\x,0) circle (2pt) node[below=1pt] {\l'};
        }
        \draw (0,1) .. controls (0.5,0.5) and (1.5,0.5) .. (2,1) -- (1,0);
        \draw (1,1) -- (3,0) -- (4,0);
        \draw (2,0) .. controls (3,0.5) and (4,0.5) .. (5,0);
      \end{tikzpicture}
    \end{center}
  \item For a $P$ a partition of $\{1, \dots, k, 1', \dots, l'\}$ and $Q$ a partition of $\{1', \dots, l', 1'', \dots, m''\}$, consider the diagram we get by putting $P$ on top of $Q$. Then $N(P, Q)$ is the number of connected components that consist only of elements of $\{1', \dots, l'\}$, and $P * Q$ is given by restricting to $\{1, \dots, k, 1'', \dots, m''\}$. For example, composing $e_P$ from above with the map $e_Q: [6] \to [2]$ given by the partition $([1', 3'], [2', 2''], [4', 1''], [5'], [6'])$ gives us
    \begin{equation*}
     \vcenter{\hbox{
          \begin{tikzpicture}
            \foreach \x [count=\l] in {0,...,2} {
              \filldraw[black] (\x,1) circle (2pt) node[above=1pt] {\l};
            }
            \foreach \x [count=\l] in {0,...,1} {
              \filldraw[black] (\x,-1) circle (2pt) node[below=1pt] {\l''};
            }
            \draw (0,1) .. controls (0.5,0.5) and (1.5,0.5) .. (2,1) -- (1,0);
            \draw (1,1) -- (3,0) -- (4,0);
            \draw (2,0) .. controls (3,0.5) and (4,0.5) .. (5,0);
            \draw (0,0) .. controls (0.5,-0.5) and (1.5,-0.5) .. (2,0);
            \draw (1,0) -- (1,-1);
            \draw (3,0) -- (0,-1);
            \foreach \x [count=\l] in {0,...,5} {
              \filldraw[black] (\x,0) circle (2pt) node[below=1pt] {\l'};
            }
          \end{tikzpicture}
       } }
      \to P * Q =
     \vcenter{\hbox{
          \begin{tikzpicture}
            \foreach \x [count=\l] in {0,...,2} {
              \filldraw[black] (\x,1) circle (2pt) node[above=1pt] {\l};
            }
            \foreach \x [count=\l] in {0,...,1} {
              \filldraw[black] (\x,0) circle (2pt) node[below=1pt] {\l''};
            }
            \draw (0,1) .. controls (0.5,0.5) and (1.5,0.5) .. (2,1) -- (1,0);
            \draw (1,1) -- (0,0);
          \end{tikzpicture}
        }     },
      N(P, Q) = 1
    \end{equation*}
  \end{enumerate}
\end{Example}
Visualizing maps with this notation allows us to easily compute compositions.
\begin{Proposition}{}{}
  $[m]$ is self-dual for all nonnegative integers $m$. 
\end{Proposition}
\begin{proof}
  The diagrams for $\ev_{[m]}$ and $\coev_{[m]}$ are
  \begin{equation*}
    \ev_{[m]} =
    \vcenter{\hbox{
        \begin{tikzpicture}
          \foreach \x in {0, 2, 3, 5} {
            \filldraw[black] (\x, 0) circle (2pt);
          }
          \foreach \x in {1, 4} {
            \draw (\x, 0) node {$\cdots$};
          }
          \draw (0, 0) .. controls (1, -1) and (2, -1) .. (3, 0);
          \draw (2, 0) .. controls (3, -1) and (4, -1) .. (5, 0);
          \draw (0, 0) node[above=1pt] {$1$};
          \draw (3, 0) node[above=1pt] {$m + 1$};
          \draw (2, 0) node[above=1pt] {$m$};
          \draw (5, 0) node[above=1pt] {$2m$};
        \end{tikzpicture}
     }    }
  \end{equation*}
  \begin{equation*}
    \coev_{[m]} =
    \vcenter{\hbox{
        \begin{tikzpicture}
          \foreach \x in {0, 2, 3, 5} {
            \filldraw[black] (\x, 0) circle (2pt);
          }
          \foreach \x in {1, 4} {
            \draw (\x, 0) node {$\cdots$};
          }
          \draw (0, 0) node[below=1pt] {1'} .. controls (1, 1) and (2, 1) .. (3, 0) node[below=1pt] {$(m + 1)'$};
          \draw (2, 0) node[below=1pt] {$m'$} .. controls (3, 1) and (4, 1) .. (5, 0) node[below=1pt] {$(2m)'$};
        \end{tikzpicture}
      }   }
  \end{equation*}
  Then, the composition $(\id_{[m]} \otimes \ev_{[m]}) \circ (\coev_{[m]} \otimes \id_{[m]})$ is
  \begin{equation*}
    \begin{tikzpicture}
      \foreach \x in {0, 2, 3, 5, 6, 8} {
        \filldraw[black] (\x, 0) circle (2pt);
      }
      \foreach \x in {1, 4, 7} {
        \draw (\x, 0) node {$\cdots$};
      }
      \draw (1, -1) node {$\cdots$};
      \draw (7, 1) node {$\cdots$};
      \foreach \x in {0, 2} {
        \filldraw[black] (\x, -1) circle (2pt);
        \draw (\x, 0) -- (\x, -1);
        \draw (\x, 0) .. controls ({\x + 1}, 1) and ({\x + 2}, 1) .. ({\x + 3}, 0);
      }
      \foreach \x in {6, 8} {
        \filldraw[black] (\x, 1) circle (2pt);
        \draw (\x, 0) -- (\x, 1);
        \draw (\x, 0) .. controls ({\x - 1}, -1) and ({\x - 2}, -1) .. ({\x - 3}, 0);
      }
    \end{tikzpicture}
  \end{equation*}
  i.e. the identity, and likewise $(\id_{[m]} \otimes \coev_{[m]}) \circ (\ev_{[m]} \otimes \id_{[m]}) = \id_{[m]}$.
\end{proof}

\begin{Definition}{}{}
  The trace of a map $f: [m] \to [m]$, denoted $\Tr f \in \Hom([0], [0]) = F$, is the composition
  \begin{equation*}
    [0] \xrightarrow{\coev_{[m]}} [m] \otimes [m] \xrightarrow{f \otimes \id_{[m]}} [m] \otimes [m] \xrightarrow{\ev_{[m]}} [0].
  \end{equation*}
  Pictorially, the trace looks like
  \begin{equation*}
    \vcenter{\hbox{\begin{tikzpicture}
      \foreach \x in {0, 2, 3, 5} {
        \filldraw[black] (\x, 0) circle (2pt);
        \filldraw[black] (\x, 1) circle (2pt);
        \filldraw[black] (\x, -1) circle (2pt);        
      }
      \foreach \x in {1, 4} {
        \draw (\x, 0) node {$\cdots$};
        \draw (\x, 1) node {$\cdots$};
        \draw (\x, -1) node {$\cdots$};
      }
      \foreach \x in {0, 2} {
        \draw (\x, -1) .. controls ({\x + 1}, -1.5) and ({\x + 2}, -1.5) .. ({\x + 3}, -1);
        \draw (\x, 1) .. controls ({\x + 1}, 1.5) and ({\x + 2}, 1.5) .. ({\x + 3}, 1);
        \draw (\x, 0) -- ({\x + 3}, -1);
        \draw (\x, -1) -- ({\x + 3}, 0);
        \draw[dashed] (\x, 1) -- (\x, 0);
      }
      \foreach \x in {3, 5} {
        \draw (\x, 1) -- (\x, 0);
      }
      \draw (1, 0.5) node {$P$};
    \end{tikzpicture}}}.
  \end{equation*}
  Using the law of composition described above, we see that $\Tr(e_P) = t^{l(P)}$, where $l(P)$ is the number of connected components in the above diagram; equivalently, the number of connected components in
  \begin{equation*}
    \vcenter{\hbox{\begin{tikzpicture}
      \foreach \x in {0, 2, 3, 5} {
        \filldraw[black] (\x, 0) circle (2pt);
        \filldraw[black] (\x, 1) circle (2pt);
      }
      \foreach \x in {1, 4} {
        \draw (\x, 0) node {$\cdots$};
        \draw (\x, 1) node {$\cdots$};
      }
      \foreach \x in {0, 2} {
        \draw[dashed] (\x, 1) -- (\x, 0);
      }
      \foreach \x in {3, 5} {
        \draw (\x, 1) -- (\x, 0);
      }
      \draw (3, 1) arc (0:180:0.5);
      \draw (2, 0) arc (180:360:0.5);
      \draw (5, 1) arc (0:180:2.5 and 1);
      \draw (0, 0) arc (180:360:2.5 and 1);
      \draw (1, 0.5) node {$P$};
    \end{tikzpicture}}}.
  \end{equation*}
\end{Definition}
\begin{Corollary}{}{}
  The dimension of $[m]$ is $\Tr \id_{[m]} = t^m$, thus interpolating the dimensions of tensor powers of the permutation representation for integer $n$.
\end{Corollary}

\subsection{Universal Property}

\label{sec:univ_prop_st}
Important for us in relating this construction to the ultraproduct construction is the universal property of $\Rep(S_t)$.
We first describe a universal property of $\Rep(S_n)$ to give some motivation. As above, let $X = \{1, \dots, n\}$ and $V$ be the permutation representation. We can identify $V$ with $\Fun(X, \mathbb{C})$, so it actually has a natural structure of an associative, commutative, unital algebra (from now on when we say algebra we assume these properties) through pointwise multiplication of functions. This allows us to define the trace of an element of $f \in V$ as $\Tr l_f$, where $l_f$ is left multiplication by $f$, and in particular the map
\begin{align*}
  V \otimes V &\to \mathbb{C} \\
  x \otimes y &\mapsto \Tr(xy)
\end{align*}
is nondegenerate. Note also that $\Tr 1_X$, the function taking $X \to \mathbf{1}$, is $n$.
\begin{Theorem}{}{repsn_univprop}
  Let $\mathcal{T}$ be a tensor category and $T \in \Ob(\mathcal{T})$ be an $n$-dimensional algebra with $\wedge^{n + 1}T = 0$ so that the pairing $x \otimes y \mapsto \Tr(xy)$ is nondegenerate. Suppose that $T$ also has the action of a finite group $G$. Then we have an essentially unique (i.e. unique up to isomorphism) functor $\mathcal{F}: \Rep(S_n) \to \mathcal{T}$ with $\mathcal{F}(V) \cong T$ as $G$-algebras. In this sense, $\Rep(S_n)$ is a universal category for $n$-dimensional representations of finite groups.
\end{Theorem}

$\Rep(S_t)$ has an analogous universal property, but with $t$-dimensional objects in arbitrary symmetric tensor categories. Suppose $T$ is an object satisfying the following three conditions:
\begin{enumerate}
\item $T$ is an algebra (here multiplication is given by a bilinear map $\mu: T \otimes T \to T$ which satisfies coherence conditions corresponding to associativity and commutativity).
\item $T$ admits a dual. Moreover, the map $T \otimes T \xrightarrow{\mu_T} T \xrightarrow{\Tr} 1$ is nondegenerate, making $T$ self-dual via the identification $\Hom(T \otimes T, \mathbf{1}) \sim \Hom(T, T^*)$. Here we again define $\Tr$ as the trace of left multiplication, namely the composition
  \begin{equation}
    \label{eq:trace}
    T \xrightarrow{\id_T \otimes \coev_T} T \otimes T \otimes T^* \xrightarrow{\mu_T \otimes \id_{T^*}} T \otimes T^* \cong T^* \otimes T \xrightarrow{\ev_T} \mathbf{1}.
  \end{equation}
\item $\dim T = \Tr(\id_T) = t$.
\end{enumerate}
In other words, $T$ is a $t$-dimensional symmetric Frobenius algebra with nondegenerate trace. We see that $[1]$ satisfies the above properties too - the unit map is given by the partition $\{[1]\}$, the multiplication map $[1] \otimes [1] \to [1]$ corresponds to the partition $\{[1], [2], [1']\}$, and it's easy to then check that the properties hold.
\begin{Theorem}{}{}
  Let $\mathcal{T}$ be a symmetric monoidal category and $T \in \Ob(\mathcal{T})$ an object satisfying the above three properties. Then there is a unique (up to isomorphism) additive symmetric monoidal functor $\Rep(S_t) \to \mathcal{T}$ sending $[1] \mapsto T$. In other words, there is an equivalence of categories between the category of isomorphism classes of additive symmetric monoidal functors $\Rep(S_t) \to \mathcal{T}$ and the category of isomorphism classes of $t$-dimensional symmetric Frobenius algebras with nondegenerate trace $T \in \Ob(\mathcal{T})$.
\end{Theorem}
\begin{proof}
  We want to define a functor $\mathcal{F}$ with the requirement that $\mathcal{F}([1]) = T$. Since $\Rep(S_t)$ is the Karoubian completion of $\Rep(S_t)^{(0)}$, whose objects are just $[1]^{\otimes n}$, it suffices to show that we can extend $\mathcal{F}$ to $\Rep(S_t)^{(0)}$. Because $\mathcal{F}$ will be a monoidal functor, it must send $[1]^{\otimes n} \mapsto T^{\otimes n}$. For a morphism $e_P: [1]^{\otimes l} \to [1]^{\otimes m}$, let $T_P := \mathcal{F}e_p: T^{\otimes l} \to T^{\otimes m}$.
  We claim that $T_P T_Q = t^{N(P, Q)}T_{P * Q}$. If we can do this, we have defined $\mathcal{F}$ on tensor powers of $T$, so we are done.
\end{proof}
\begin{Proposition}{}{t_compn}
  $T_P T_Q = t^{N(P, Q)}T_{P * Q}$.
\end{Proposition}
\begin{proof}
  First, we must define $e_P$ in a more ``categorical'' way so that its definition can be extended to objects of an arbitrary category. Let $P_l = P \cap \{1, \dots, l\}$ and $P_m = P \cap \{1', \dots, m'\}$. We see that $e_P$ can be written only in terms of multiplication and its dual, so
  \begin{equation*}
    e_P = \otimes_P ([1]^{\otimes P_l} \xrightarrow{m_{P_l}} [1] \xrightarrow{m_{P_m}^*} [1]^{\otimes P_m})
  \end{equation*}
  where $m_{P_l}$ is the multiplication map, and we are taking the tensor product over all the parts of $P$. Analogously,
  \begin{equation*}
    T_P = \otimes_P (T^{\otimes P_l} \xrightarrow{m_{P_l}} T \xrightarrow{m_{P_m}^*} T^{\otimes P_m}).
  \end{equation*}
  Thus, we must analyze what each composition $m_{P_m}^* \circ m_{P_l}$ looks like.
  Next, we identify $\Hom(T^{\otimes l}, T^{\otimes m})$ with $\Hom(T^{\otimes l + m}, \mathbf{1})$ by sending $f$ to the bilinear form given by
  \begin{equation*}
    B_f: T^{\otimes l} \otimes T^{\otimes m} \xrightarrow{f \otimes \id_{T^{\otimes m}}} T^{\otimes m} \otimes T^{\otimes m} \xrightarrow{\ev_{T^{\otimes m}}} \mathbf{1}.
  \end{equation*}
  This bijection of Hom-spaces interacts nicely with composition and duals, and will allow us to characterize $T_P$ in such a way that the desired composition law is easy to see.
  \begin{Proposition}{}{paren}
    The above identification has the following properties:
    \begin{enumerate}
    \item For a composition $gf : T^{\otimes k} \to T^{\otimes l} \to T^{\otimes m}$, $B_{gf}$ is the composition
      \begin{equation*}
        T^{\otimes k} \otimes T{\otimes m} \xrightarrow{\id_{T^{\otimes k}} \otimes \coev_{T^{\otimes l}} \otimes \id_{T^{\otimes m}}} T^{\otimes k} \otimes T^{\otimes l} \otimes T^{\otimes l} \otimes T^{\otimes m} \xrightarrow{B_f \otimes B_g} \mathbf{1}
      \end{equation*}
    \item $B_{f^*} = B_f$
    \item For $f: T^{\otimes l} \to T^{\otimes m}$ and $g: T^{\otimes k} \to T^{\otimes m}$, $B_{g^* f}$ simplifies to
      \begin{equation*}
        T^{\otimes l} \otimes T^{\otimes k} \xrightarrow{f \circ g} T^{\otimes m} \otimes T^{\otimes m} \xrightarrow{\ev_{T^{\otimes m}}} \mathbf{1}.
      \end{equation*}
    \end{enumerate}
  \end{Proposition}
  We omit the proof because it involves much diagram chasing.
  
  In order to describe $B_{T_P}$ and $B_{T_Q}$, we need to know three things: what $B_{m_p}$ looks like, where $m_p: T^{\otimes p} \to T$ sends $t_1 \otimes \cdots \otimes t_p \mapsto \prod t_1 \cdots t_p$, and how $B_{T_P}$ interacts with coevaluation, which will tell us how to compute its composition with $B_{T_Q}$.

  Recall that $T$ is self-dual via the map $T \otimes T \xrightarrow{\mu} T \xrightarrow{\Tr} \mathbf{1}$, where $\mu$ is the algebra multiplication map. This composition is $\ev_T: T \otimes T \to \mathbf{1}$, so $B_{m_p}$ is the composition
  \begin{equation*}
    \Tr m_{p + 1} : T^{\otimes p} \otimes T \xrightarrow{m_p \otimes \id_T} T \otimes T \xrightarrow{\mu} T \xrightarrow{\Tr} \mathbf{1}.
  \end{equation*}
  Moreover, by \cref{propn:paren}, $B_{m_q^* m_p}$ will be a map 
  \begin{equation*}
    \Tr m_{p + q}: T^{\otimes p} \otimes T^{\otimes q} \xrightarrow{m_p \otimes m_q} T \otimes T \xrightarrow{\mu} T \xrightarrow{\Tr} \mathbf{1}.
  \end{equation*}
  So in fact $B_{T_P}$ can be written nicely as $\otimes (T^{\otimes P} \xrightarrow{m_P} T \xrightarrow{\Tr} \mathbf{1})$. Then, by \cref{propn:paren}, $B_{T_Q \circ T_P}$ is a composition
  \begin{equation*}
    T^{\otimes k} \otimes T^{\otimes m} \to (T^{\otimes k} \otimes T^{\otimes l}) \otimes (T^{\otimes l} \otimes T^{\otimes m}) \to \mathbf{1}
  \end{equation*}
  where we end up ``contracting'' $\coev_T: \mathbf{1} \to T \otimes T$ $k$ times. In fact, we need only look at this for every part $R$ of $P * Q$.
  For how $B_{T_P}$ and $B_{T_Q}$ (which we can now just think of as multiplication maps) interact with coevaluation, we look at our maps $B_{m_p} = \Tr m_{p+ 1}$. By \cref{propn:paren},
  \begin{equation*}
    T^{\otimes p} \otimes T^{\otimes q} \xrightarrow{\id_{T^{\otimes p}} \otimes \coev_T \otimes \id_{T^{\otimes q}}} T^{\otimes p} \otimes T \otimes T \otimes T^{\otimes q} \xrightarrow{\Tr m_{p + 1} \otimes \Tr m_{q + 1}} \mathbf{1}
  \end{equation*}
  is by definition $B_{m_q^* m_p}$, which is
  \begin{equation*}
    \Tr m_{p + q}: T^{\otimes p} \otimes T^{\otimes q} \xrightarrow{m_p \otimes m_q} T \otimes T \xrightarrow{\mu} T \xrightarrow{\Tr} \mathbf{1}.
  \end{equation*}
  If $R$ intersects $\{1, \dots, k\}$ or $\{1, \dots, m\}$, then the corresponding part of the composition will be repeated multiplication $T^{\otimes R_k} \otimes T^{\otimes R_m} \to \mathbf{1}$, and the tensor product of all of these will be $B_{T_{P * Q}}$. However, if it only intersects $\{1, \dots, l\}$ (this corresponds to the $y_i$ that don't have to be equal to any $x_i$ or $z_i$), then the corresponding part of the composition will be the map
  \begin{equation*}
    f: \mathbf{1} \xrightarrow{\coev_T} T \otimes T \xrightarrow{\mu} T \xrightarrow{\Tr} \mathbf{1}.
  \end{equation*}
  The associativity of $\mu$ and our definition of trace in \cref{eq:trace} imply that $f = \Tr = B_{\mu \circ \coev_T}$. But consider also $B_u$, where $u$ is the unit of $T$:
    \begin{equation*}
    B_u: T = \mathbf{1} \otimes T \xrightarrow{u \otimes T} T \otimes T \xrightarrow{\mu} T \xrightarrow{\Tr} \mathbf{1}
  \end{equation*}
  is also $\Tr$. So $u = \mu \circ \coev_T$ and $f = \Tr u = t$.
  Thus, we have a total factor of $t^{N(P, Q)}$ and we are done.
\end{proof}
\begin{Remark}{}{}
  In fact, the universal property of $\Rep(S_n)$ we stated in \cref{theo:repsn_univprop} can be derived from the universal property of $\Rep(S_t)$ when $t = n$. As we will see in \cref{sec:interpolation}, we obtain $\Rep(S_n)$ from $\Rep(S_t)$ by quotienting out by a tensor ideal. In fact, this tensor ideal is the one generated by the identity map on $\wedge^{n + 1} [1]$, hence the requirement that $\wedge^{n + 1} T = 0$.
\end{Remark}
\begin{Definition}{}{interp_functor}
  The universal property defines the \textbf{interpolation functor} $\mathcal{F}_n: \Rep(S_n) \to \overline{\Rep}(S_n)$ given by $\mathcal{F}_n[1] = V$, where $V$ is the permutation representation of $S_n$ and $\overline{\Rep}(S_n)$ is the classical category of representations of $S_n$.
\end{Definition}
We will discuss the interpolation functor in more detail in \cref{sec:interpolation}.
\begin{Remark}{}{}
  For $t \in \mathbb{Z}_{\ge 0}$, $\Rep(S_t)$ is not an abelian category. However, we can obtain the abelian envelope $\Rep^{ab}(S_t)$, as defined in \cite{deligne_ct} 8.19, by closing under subquotients, and in \cref{sec:interpolation} we will show how to show how to recover $\Rep(S_n)$, an abelian category, by quotienting out a tensor ideal.
\end{Remark}

\section{Analogous constructions for $GL_t$, $O_t$, $Sp_t$}
\label{sec:glosp_interp}
Using $\Rep(S_t)$ as a model, we can also construct $\Rep(GL_t)$, $\Rep(O_t)$, and $\Rep(Sp_t)$. These have very similar constructions: we describe the Hom-spaces between tensor powers of the standard $n$-dimensional representation (and its dual, for $GL_t$) and the composition law, then take a Karoubian completion. We combine the exposition in \cite{EGNO} and \cite{comes_wilson}.

\subsection{Construction and universal property of $\Rep(GL_t)$}
\label{sec:gl_construction}
Again, we first analyze tensor powers of the standard representation $V$ for the classical category $\Rep(GL_n)$. Every irreducible representation of $GL_n$ lies in $V^{\otimes l} \otimes V^{*\otimes m}$ for some $l, m$. This leads us to consider the category generated by tensor powers of both $V$ and $V^*$, and we want to characterize
\begin{equation*}
  \Hom_{GL_n}(V^{\otimes l_1} \otimes V^{* \otimes m_1}, V^{\otimes l_2} \otimes V^{* \otimes m_2}) = \Hom_{GL_n}(V^{\otimes (l_1 + m_2)}, V^{\otimes (l_2 + m_1)}) = \left(V^{\otimes (l_1 + m_2)} \otimes V^{* \otimes (l_2 + m_1)}\right)^{GL_n}.
\end{equation*}
This is nonzero only when $r := l_1 + m_2 = l_2 + m_1$, in which case we are characterizing $\End_{GL_n}(V^{\otimes r})$. By Schur-Weyl duality, $\End_{GL_n}(V^{\otimes r}) \cong \mathbb{C}S_r$ when $n \ge r$. Moreover, treating $\left(V^{\otimes r} \otimes V^{* \otimes r}\right)^{GL_n}$ as $GL_n$-invariant polynomials, the fundamental theorem of invariant theory gives that the elements $\sum_{i_1, \dots, i_r} e_{i_1} \otimes \cdots \otimes e_{i_r} \otimes e^*_{i_{\sigma(1)}} \otimes \cdots e^*_{i_{\sigma(r)}}$ span $\left(V^{\otimes r} \otimes V^{* \otimes r}\right)^{GL_n}$, where $e$ is a basis for $V$, $e^*$ its dual basis, and $\sigma \in S_r$. These form a basis when $n \ge r$.

Now we define composition of morphisms. Again, we think about how to represent the above maps pictorially. Every permutation $\sigma \in S_r$ corresponds to a pairing from $\{1, \dots, l_1 + m_2\}$ to $\{1, \dots, l_2 + m_1\}$. Representing copies of $V$ as black dots and copies of $V^*$ as white dots, we get a diagram between two rows; the first has $l_1$ black dots and $m_1$ white dots while the second has $l_2$ black dots and $m_2$ white dots, and our permutation connects pairs of these dots. We will call these black-white diagrams.
\begin{Example}{}{}
  The evaluation map $V^* \otimes V \to \mathbb{C}$ is represented by
  \begin{equation*}
    \begin{tikzpicture}
      \filldraw[black] (1, 0) circle (2pt);
      \draw (1, 0) arc (0:170:0.5);
      \filldraw[white] (0, 0) circle (2pt);
      \draw (0, 0) circle (2pt);
    \end{tikzpicture}.
  \end{equation*}
  Similarly, the coevaluation map $\mathbb{C} \to V \otimes V^*$ is represented by
  \begin{equation*}
    \begin{tikzpicture}
      \filldraw[black] (0, 0) circle (2pt);
      \draw (0, 0) arc (180:360:0.5);
      \filldraw[white] (1, 0) circle (2pt);
      \draw (1, 0) circle (2pt);
    \end{tikzpicture}.
  \end{equation*}
  The identity map $V^{\otimes l} \otimes V^{* \otimes m} \to V^{\otimes l} \to V^{* \otimes m}$ is represented by
  \begin{equation*}
    \begin{tikzpicture}
      \foreach \x in {0, 2} {
        \filldraw[black] (\x, 1) circle (2pt) -- (\x, 0) circle (2pt);
      }
      \foreach \x in {1, 4} {
        \draw (\x, 0) node {$\cdots$};
        \draw (\x, 1) node {$\cdots$};        
      }
      \foreach \x in {3, 5} {
        \draw (\x, 0) -- (\x, 1);
        \filldraw[white] (\x, 1) circle (2pt);
        \filldraw[white] (\x, 0) circle (2pt);
        \draw (\x, 0) circle (2pt);
        \draw (\x, 1) circle (2pt);
      }
      \draw[decorate,decoration={brace,amplitude=10pt,mirror},xshift=0pt,yshift=-4pt] (0, 0) -- (2, 0) node[midway,yshift=-20pt] {$l$};
      \draw[decorate,decoration={brace,amplitude=10pt,mirror},xshift=0pt,yshift=-4pt] (3, 0) -- (5, 0) node[midway,yshift=-20pt] {$m$};
    \end{tikzpicture}.
  \end{equation*}
\end{Example}
In general, for any black-white diagram, an edge connecting two dots in the top row corresponds to evaluation, an edge connecting two dots in the bottom row corresponds to coevaluation, and an edge connecting two dots in different rows corresponds to the identity map. As in the symmetric group case, composition concatenates two diagrams, then adds factors of $n$ based on what happens in the ``middle'' row. Consider $\ev_V \circ c_{V, V^*} \circ \coev_V$, which we know is multiplication by $n$; diagrammatically, this composition contains a closed loop:
\begin{equation*}
  \begin{tikzpicture}
    \filldraw[black] (1, 1) circle (2pt);
    \filldraw[black] (0, 0) circle (2pt);
    \draw (0, 0) -- (1, 1) arc (0:180:0.5) -- (1, 0) arc (360:180:0.5);
    \filldraw[white] (1, 0) circle (2pt);
    \filldraw[white] (0, 1) circle (2pt);
    \draw (1, 0) circle (2pt);
    \draw (0, 1) circle (2pt);
  \end{tikzpicture}.
\end{equation*}
In general, every closed loop simplifies to a $\coev_V \circ \ev_V$ composition, thus adding a factor of $n$.

\begin{Proposition}{}{}
  Let $r = l_1 + m_2 = l_2 + m_1$ and $\sigma, \sigma' \in S_r$ correspond to two basis elements $e_\sigma, e_{\sigma'}$ of $\Hom_{GL_n}(V^{\otimes l_1} \otimes V^{* \otimes m_1}, V^{\otimes l_2} \otimes V^{* \otimes m_2})$. Then $e_\sigma \circ e_{\sigma'} = n^l e_{\sigma \circ \sigma'}$, where $l$ is the number of loops removed when we concatenate the diagrams corresponding to $\sigma, \sigma'$ and remove the middle row.
\end{Proposition}

Now we can construct $\Rep(GL_t)$.
\begin{Definition}{}{}
  We first define the tensor category $\Rep(GL_t)^{(0)}$. Its objects are of the form $[l, m]$ for nonnegative integers $l, m$, and the tensor product $[l_1, m_1] \otimes [l_2, m_2] = [l_1 + l_2, m_1 + m_2]$. Its Hom-spaces are
  \begin{equation*}
    \Hom([l_1, m_1], [l_2, m_2]) = \begin{cases} 0 & l_1 + m_2 \ne l_2 + m_1 \\
      \mathbb{C}S_{l_1 + m_2} & l_1 + m_2 = l_2 + m_1 \end{cases}
  \end{equation*}
  and composition is given as above, except every removed loop gives a factor of $t$ instead of $n$. The tensor product on basis elements of the Hom-space is the ``union'' of their diagrams (and thus their permutations). The unit object is $[0, 0]$.

  $\Rep(GL_t)$ is the Karoubian envelope of the additive completion of $\Rep(GL_t)^{(0)}$.
\end{Definition}
\begin{Definition}{}{}
  The \textbf{walled Brauer algebra} $B_{l, m}(t)$ is $\mathbb{C}S_{l, m}$ with the composition law above: i.e. $\End_{\Rep(GL_t)}([l, m])$.
\end{Definition}

By a similar analysis as the $S_t$ case, we see that the dimension of $[r, s]$ is $t^{r + s}$, and the dual of $[r, s]$ is $[s, r]$. For example, $[1, 0]$, corresponding to the standard representation in $\Rep(GL_n)$, has dimension $t$.

Like $\Rep(S_t)$, the category $\Rep(GL_t)$ has a universal property. Roughly, it is the universal category generated by an object of dimension $t$ and its dual.
\begin{Theorem}{}{}
  Let $\mathcal{T}$ be a rigid additive symmetric monoidal category. Then there is an equivalence of categories between the category of isomorphism classes of symmetric additive monoidal functors $\Rep(GL_t) \to \mathcal{T}$ and isomorphism classes of objects of dimension $t$ in $\mathcal{T}$, given on objects by $\mathcal{F} \mapsto \mathcal{F}([1, 0])$ and morphisms by $\eta \mapsto \eta_{[1, 0]}$.
\end{Theorem}
\begin{proof}
  We show that the functor $-([1, 0]): \mathcal{F} \mapsto \mathcal{F}([1, 0])$ is indeed an equivalence of categories. Because $\Rep(GL_t)$ is the Karoubian envelope of the additive completion of $\Rep(GL_t)^{(0)}$, it suffices to define $\mathcal{F}$ on this smaller category. Recall that $\Rep(GL_t)^{(0)}$ is the category generated by tensor products and powers of $[1, 0]$ and $[0, 1]$. We show that $\mathcal{F}$ is essentially surjective (\cref{propn:glt-surjective}), full (\cref{propn:glt-full}), and faithful(\cref{propn:glt-faithful}).
\end{proof}
\begin{Proposition}{}{glt-surjective}
  $-([1, 0])$ is essentially surjective.
\end{Proposition}
\begin{proof}
  Let $T$ be any object of dimension $t$ in $\mathcal{T}$; it has a dual $T^*$ and evaluation/coevaluation maps by assumption. Because monoidal functors preserve duals, we must have $\mathcal{F}([0, 1]) = T^*$, $\mathcal{F}(\ev_{[1, 0]}) = \ev_T$, $\mathcal{F}(\coev_{[1, 0]}) = \coev_T$. This is well-defined, for both $T$ and $[1, 0]$ have dimension $t$ and we know $\mathcal{F}$ must preserve  dimension. In fact, this is enough to determine $\mathcal{F}$, because any diagram is made up only of tensor products of identity, evaluation, and coevaluation maps. Thus our functor is essentially surjective.
\end{proof}

\begin{Proposition}{}{glt-full}
  $-([1, 0])$ is full.
\end{Proposition}
\begin{proof}
  Let $f: T \to U$ be a morphism between two objects $T, U$ of dimension $t$ in $\mathcal{T}$, and $\mathcal{F}_X, \mathcal{F}_Y$ be the two functors sending $[1, 0]$ to $X, Y$ respectively. We want to construct a natural transformation $\eta: \mathcal{F}_X \to \mathcal{F}_Y$ so that $\eta_{[1, 0]} = f$. We will define $\eta$ as
  \begin{equation*}
    \eta_{[1, 0]} = f, \qquad \eta_{[0, 1]} = (f^{-1})^*, \qquad \eta_{[0, 0]} = \id_{\mathbf{1}}
  \end{equation*}
  and extend to the other objects in $\Rep(GL_t)^{(0)}$ by $\eta_{[r, s] \otimes [l, m]} := \eta_{[r, s]} \otimes \eta_{[l, m]}$. Now we must verify that $\eta$ is natural.

  First, we claim that $\eta$ respects the braiding; then we can reorder our tensor products while checking naturality. Let $A, B$ be two objects in $\Rep(GL_t)^{(0)}$ and $c, d$ be the braidings in $\Rep(GL_t)^{(0)}$ and $\mathcal{T}$ respectively. The diagram
  \begin{equation*}
    \begin{tikzcd}
      \mathcal{F}_X(A \otimes B) \arrow[r, equal] \arrow[d, "\mathcal{F}_X(c_{A, B})"'] & \mathcal{F}_X A \otimes \mathcal{F}_X B \arrow[r, "\eta_A \otimes \eta_B"] \arrow[d, "d_{\mathcal{F}_X A, \mathcal{F}_X B}"'] & \mathcal{F}_Y A \otimes \mathcal{F}_Y B \arrow[r, equal] \arrow[d, "d_{\mathcal{F}_Y A, \mathcal{F}_Y B}"] & \mathcal{F}_Y (A \otimes B) \arrow[d, "\mathcal{F}_Y(c_{A, B})"] \\
      \mathcal{F}_X(B \otimes A) \arrow[r, equal] & \mathcal{F}_X B \otimes \mathcal{F}_X A \arrow[r, "\eta_B \otimes \eta_A"'] & \mathcal{F}_Y B \otimes \mathcal{F}_Y A \arrow[r, equal] & \mathcal{F}_Y(B \otimes A)
    \end{tikzcd}
  \end{equation*}
  commutes by the definitions and properties of a symmetric monoidal category. Therefore, $\eta$ respects the braiding.

  Now we only need to check naturality of $\eta$ on four maps: the identity maps for $[1, 0]$ and $[0, 1]$, and coevaluation and evaluation maps for $[1, 0]$. It's clear that $\eta$ is natural on the identity maps. On evaluation, naturality means
  \begin{equation*}
    \begin{tikzcd}
      X^* \otimes X \arrow[r, "\ev_X"] \arrow[d, "(f^{-1})^* \otimes f"'] & \mathbf{1} \arrow[d, equal] \\
      Y^* \otimes Y \arrow[r, "\ev_Y"'] & \mathbf{1}
    \end{tikzcd}
  \end{equation*}
  commutes, but that is equivalent to \cref{eq:dual}. Naturality of coevaluation is analogous. So we are done.
\end{proof}

\begin{Proposition}{}{glt-faithful}
  $-([1, 0])$ is faithful.
\end{Proposition}
\begin{proof}
  Every natural transformation $\eta: \mathcal{F}_X \to \mathcal{F}_Y$ is determined by its value on $[1, 0]$ and $[0, 1]$; moreover, because $[0, 1] = [1, 0]^*$, $\eta_{[0, 1]} = (\eta_{[1, 0]}^{-1})^*$, so $\eta$ is determined uniquely by $\eta_{[1, 0]}$.
\end{proof}

\begin{Corollary}{}{}
  The universal property of $\Rep(GL_t)$ defines for us three interpolation functors. If $t = n \in \mathbb{Z}_{\ge 0}$, we have a functor $\Rep(GL_n) \to \overline{\Rep}(GL_n)$. If $t = n \in \mathbb{Z}_{< 0}$, we have a functor $\Rep(GL_n) \to \overline{\Rep}(GL_n, \mathbb{Z}/2\mathbb{Z})$, the category of $\mathbb{Z}/2\mathbb{Z}$-graded representations of $GL_n$, graded by the action of $-I$.\footnote{In this category, every representation has an even part, where $-I$ acts as $\id$, and an odd part, where $-I$ acts as $-\id$.} These two functors will be discussed in \cref{sec:interpolation}.

  We also have a functor $\Rep(GL_{t = m - n}) \to \overline{\Rep}(GL_{m|n})$. Supergroups are outside the scope of this paper, but we refer the reader to \cite{comes_wilson} for a discussion of this functor.
\end{Corollary}

\subsection{Construction and universal properties of $\Rep(O_t)$ and $\Rep(Sp_t)$}

We proceed similarly to the construction of $\Rep(GL_t)$. Let $V$ be the standard $n$-dimensional representation of $O_n$; since every irreducible $O_n$-representation occurs in $V^{\otimes m}$ for some $m$, we will first look at $\Hom_{O_n}(V^{\otimes l}, V^{\otimes m})$. Here $V$ is self-dual, so
\begin{equation*}
  \Hom_{O_n}(V^{\otimes l}, V^{\otimes m}) = (V^{\otimes l + m})^{O_n}
\end{equation*}
and is nonzero iff $l + m$ is even, i.e. $l + m = 2k$ for some $k$, and then $(V^{\otimes 2k})^{O_n} = \End_{O_n}(V^{\otimes k})$. If $n \ge k$, this endomorphism algebra is isomorphic to the Brauer algebra $B_k(n)$. It has a basis consisting of perfect matchings of $2k$ points, which we draw as two rows of $l + m$ points connected by edges corresponding to matchings. When $n < k$, these matchings span the endomorphism algebra. As in our analysis of Hom-spaces of $GL_n$, composition of two matchings is concatenation of their diagrams and the removal of the middle row, where every removed loop multiplies by a factor of $n$.
\begin{Definition}{}{}
  We first define $\Rep(O_t)^{(0)}$. Its objects are of the form $[m]$ for nonnegative integers $m$, with tensor product $[l] \otimes [m] = [l + m]$. Its Hom-spaces are
  \begin{equation*}
    \Hom([l], [m]) = \begin{cases} B_k(t) & l + m = 2k, k \in \mathbb{Z} \\
      0 & \text{otherwise} \end{cases}
  \end{equation*}
  and composition is given as above, except every removed loop adds a factor of $t$. Again, tensor product on tensor basis elements of the Hom-space is the union of the diagrams, and the unit object is $[0]$.

  $\Rep(O_t)$ is the Karoubian envelope of the additive completion of $\Rep(O_t)^{(0)}$.
\end{Definition}

$Sp_t$ is constructed completely analogously, since its Hom-spaces are also equivalent to Brauer algebras, so we omit the costruction here. Moreover, $\Rep(O_t)$ and $\Rep(Sp_{-t})$ are equivalent tensor categories via a functor that makes the braiding detect parity. The normal braiding on $\Rep(O_t)$ is the commutative map $c_{r, s}: [r] \otimes [s] \to [s] \otimes [r]$ given by swapping the factors. Now change the braiding to $c'$, defined by
\begin{equation*}
  c'_{r, s} = \begin{cases} -c_{r, s} & r, s \equiv 1 \pmod{2} \\ c_{r, s} & \text{otherwise}.\end{cases}
\end{equation*}
This new braiding makes $\Rep(O_t)$ equivalent to $\Rep(Sp_{-t})$ as a symmetric monoidal category. So proving results for $\Rep(O_t)$ also proves them for $\Rep(Sp_t)$ as long as we are careful about the action of the new braiding.

Both $\Rep(O_t)$ and $\Rep(Sp_t)$ have universal properties. Roughly, $\Rep(O_t)$ is universal among categories generated by a single object with a symmetric isomorphism to its dual, and $\Rep(Sp_t)$ is universal among categories generated by a single object with a skew-symmetric isomorphism to its dual.
\begin{Theorem}{}{}
  Let $\mathcal{T}$ be a rigid symmetric tensor category. Then there is an equivalence of categories between the category of isomorphism classes of additive symmetric monoidal functors $\Rep(O_t) \to \mathcal{T}$ (respectively $\Rep(Sp_t) \to \mathcal{T}$) and the category of isomorphism classes of $t$-dimensional objects $t \in \Ob(\mathcal{T})$ so that $T \cong T^*$ via a symmetric (respectively skew-symmetric) isomorphism, given by $\mathcal{F} \mapsto \mathcal{F}([1])$.
\end{Theorem}
\begin{Corollary}{}{}
  The universal property defines an interpolation functor $\Rep(O_{t = p - 2q}) \to \overline{\Rep}(OSp_{p|2q})$, the representations of the orthosymplectic supergroup. When $p$ or $q = 0$, we have interpolation functors to $\overline{\Rep}(O_t)$ and $\overline{\Rep}(Sp_t)$, the classical representation categories.
\end{Corollary}


\chapter{Ultraproduct construction of the Deligne categories}

In this chapter, we will present a different construction of the Deligne categories via a concept from logic, the ultraproduct, and show that the categories we get are equivalent to the ones we constructed in the first chapter. This allows us to also think of Deligne categories as a limit of classical categories. Moreover, the ultraproduct construction defines the Deligne categories in characteristic $p > 0$, whereas the previous interpolation construction would fail. In \cref{chap:research}, we show how the ultraproduct can be used to locate and construct objects with desirable properties by passing between the classical and Deligne categories.

\section{Ultraproducts and filtered ultraproducts}
\label{sec:ultraproduct_bg}

In this section, we define ultraproducts, a notion from logic that has seen much use in algebra because of its ability to preserve first-order algebraic properties. We will use these to construct $\Rep(S_t)$, $\Rep(GL_t)$, $\Rep(SO_t)$, and $\Rep(Sp_t)$. We follow the exposition in \cite{kalinov}, and we direct the reader to \cite{hans_ultraproduct} to learn more about ultraproducts in algebra.

\begin{Definition}{Ultrafilters}{}
  An \textbf{ultrafilter} $\mathcal{F}$ on a set $F$ is a collection of subsets satisfying the following properties:
  \begin{enumerate}
  \item $F \in \mathcal{F}$.
  \item $\mathcal{F}$ is closed under taking supersets, i.e. if $A \in \mathcal{F}$ and $A \subset B$, then $B \in \mathcal{F}$.
  \item $\mathcal{F}$ is closed under finite intersections, i.e. if $A, B \in \mathcal{F}$, then $A \cap B \in \mathcal{F}$.
  \item $\mathcal{F}$ contains exactly one of a set and its complement, i.e. exactly one of $A, F \setminus A$ is in $\mathcal{F}$ for $A \subset F$.
  \end{enumerate}
\end{Definition}

\begin{Definition}{}{}
  An ultrafilter of the form $\mathcal{F}_x = \{A | x \in A\}$ for $x \in F$ is a \textbf{principal ultrafilter}.
\end{Definition}

\begin{Proposition}{}{}
  If an ultrafilter contains a finite set, then it is principal.
\end{Proposition}
\begin{proof}
  Let $A$ be the finite set of smallest cardinality in $\mathcal{F}$; by definition $A$ cannot be the empty set, and if $A$ is a singleton, we're done. So suppose $A$ has cardinality at least $2$ and let $x \in A$. Then $F \setminus \{x\} \in \mathcal{F}$ because $\{x\} \notin \mathcal{F}$ by assumption. However, since $\mathcal{F}$ is closed under finite intersection, $A \setminus \{x\} \in \mathcal{F}$, so we have found a set of smaller cardinality, contradiction.
\end{proof}

The existence of non-principal ultrafilters is a consequence of the Axiom of Choice (see \cite{rubin}). From the above, non-principal ultrafilters are composed entirely of infinite subsets of $F$ and contain all the cofinite sets of $F$.

We will work entirely with non-principal ultrafilters over $\mathbb{N}$. Let $\mathcal{F}$ be an ultrafilter over $\mathbb{N}$; we note two additional properties of $\mathcal{F}$. When we say a statement is true for almost all $n$, we mean it is true for some set $A \in \mathcal{F}$. This expands our usual notion of almost all, which is all but finitely many, to other infinite subsets of $\mathbb{N}$.
\begin{Proposition}{}{condis}
\begin{enumerate}
\item If two statements are true for almost all $n \in \mathbb{N}$, then their conjunction is true for almost all $n$ as well.
\item If the disjunction of a finite number of statements is true for almost all $n$, then at least one statement is true for almost all $n$.
\end{enumerate}
\end{Proposition}
\begin{proof}
  ~\begin{enumerate}
  \item If the two statements are true for $A, B \in \mathcal{F}$ respectively, their conjunction is true for $A \cap B \in \mathcal{F}$.
  \item Let the statements be $s_i$. For contradiction, suppose no sets on which the $s_i$ hold are in $\mathcal{F}$. So let $A_i = \{n \in \mathbb{N} | s_i(n) = 0\}$. Each $A_i \in \mathcal{F}$, so $A := \bigcap A_i \in \mathcal{F}$. This represents the natural numbers for which all the $s_i$, and therefore their disjunction, are $0$. By assumption, $\mathbb{N} \setminus A$, the numbers $n$ for which $s_i(n) = 1$, is also in $\mathcal{F}$, contradiction.
  \end{enumerate}
\end{proof}

Now we can define ultraproducts.
\begin{Definition}{}{}
  Let $S_i$ be a collection of sets labeled by the natural numbers. Consider the set of sequences defined on almost all $n$, $\prod'_{\mathcal{F}} S_i$, which consists of all sets $\{s_n\}_{n \in A}$ for $A \in \mathcal{F}$ and $s_i \in S_i$. Define an equivalence relation $\sim$ with $a \sim b$ if $a_i = b_i$ for almost all $i$, so the sequences $a$ and $b$ agree on a set in the ultrafilter. The \textbf{ultraproduct} $\prod_{\mathcal{F}} S_i$ is $\prod'_{\mathcal{F}} S_i / \sim$. To denote an equivalence class of sequences, we will use $\prod_{\mathcal{F}} s_i$.
\end{Definition}
\begin{Remark}{}{}
  Each equivalence class of sequences represents the ``germ'' of a sequence, just as we might think of germs of functions.
\end{Remark}

\begin{Example}{}{}
  For a finite set $A$, $\prod_{\mathcal{F}} A = A$. Let $s \in \prod_{\mathcal{F}} A$. By \cref{propn:condis}, there exists an $a \in A$ so that $s_n = a$ for almost all $n$, so $s_n \sim \{a\}_n$.
\end{Example}

The most important property of ultraproducts is their ability to preserve maps and properties of their constituent sets. \Los{}' theorem, formally stated below, tells us that if we have sets with algebraic structure defined by maps between them, the ultraproduct also has these maps. Moreover, as long as these maps (on sets) satisfy properties or relations for almost all $n$, the ultraproduct maps also satisfy said properties. So we can take ultraproducts of groups, rings, fields, etc. and this gives a group, ring, or field structure on the ultraproduct.

To formally state the theorem, we first explain how the ultraproduct inherits operations and relations to give some intuition.
\begin{Example}{}{ultra_operations}
Let $S_n$ be our sequence of sets and define some $k$-ary operations $\circ_n$ for almost all $n$. Then we can also define a $k$-ary operation $\circ$ on $\prod_{\mathcal{F}} S_n$ as
\begin{equation*}
  \circ(s^1, \dots, s^k) = \prod_{\mathcal{F}} \circ_n (s^1_n, \dots, s^k_n)
\end{equation*}
where the $s^i$ are equivalence classes of sequences in $\prod_{\mathcal{F}} S_n$. Likewise, if we have some $k$-ary relations $r_n$, we can define a $k$-ary relation on $\prod_{\mathcal{F}} S_n$ as
\begin{equation*}
  r(s^1, \dots, s^k) = \prod_{\mathcal{F}} r_n (s^1_n, \dots, s^k_n) \in \prod_{\mathcal{F}} \{0, 1\} = \{0, 1\}.
\end{equation*}
In particular, if $r_n$ is true for almost all $n$, then so is $r$. 
\end{Example}

\begin{Theorem}{\Los{}' Theorem}{los}
  Let $S_i^{(k)}$ be a collection of sequences of sets with $k = 1, \dots, m$ and $f_i^{(r)}$ a collection of sequences of elements with $r = 1, \dots, l$. Suppose we have a formula of a first-order language $\varphi(x_1, \dots, x_l, X_1, \dots, X_m)$ where the $x_i$ are parameters and $X_m$ are sets. Let $S^{(k)} = \prod_{\mathcal{F}} S_i^{(k)}$ and $f^{(r)} = \prod_{\mathcal{F}} f_i^{(r)}$. Then $\varphi(f_n^{(1)}, \dots, f_n^{(l)}, S_n^{(1)}, \dots, S_n^{(m)})$ is true for almost all $n$ iff $\varphi(f^{(1)}, \dots, f^{(l)}, S^{(1)}, \dots, S^{(m)})$ is true.
\end{Theorem}

We now give some applications of this theorem (and non-applications) to show what algebraic properties ultraproducts preserve. 

\begin{Example}{}{ultra_vec}
As we said above, ultraproducts of groups, rings, and fields also have the structure of groups, rings, and fields respectively. This also applies to vector spaces.  However, finiteness conditions are usually not preserved. For example, taking the ultraproduct of countably many finite-dimensional vector spaces will not necessarily yield a finite-dimensional vector space because finite-dimensionality is not expressible in first-order language. However, if instead we ask for the stricter condition that these vector spaces have bounded dimension, then almost all vector spaces have the same dimension $d$ and hence the ultraproduct is also a vector space with dimension $d$.
\end{Example}

\begin{Example}{}{}
  Let $\overline{\mathbb{Q}}$ be the algebraic closure of $\mathbb{Q}$. Then $\prod_{\mathcal{F}} \overline{\mathbb{Q}} \simeq \mathbb{C}$. The ultraproduct is an algebraically closed field by \Los{}' theorem. Moreover, it has characteristic $0$ since we can write any integer $k = \prod_{\mathcal{F}} k \ne 0$, and its cardinality is the same as the reals. There is exactly one algebraically closed field of every characteristic and cardinality up to isomorphism (Steinitz' Theorem), $\prod_{\mathcal{F}} \overline{\mathbb{Q}} \simeq \mathbb{C}$. This isomorphism is not canonical, but it does give us a way of writing every complex number as a (equivalence class of a) sequence. 
\end{Example}

\begin{Example}{}{}
  By the same argument as the above, for a sequence of distinct prime numbers $p_n$, $\prod_{\mathcal{F}} \overline{\mathbb{F}_{p_n}} \simeq \mathbb{C}$, again not canonically. The ultraproduct has characteristic $0$ because $\prod_{\mathcal{F}} k$ has almost all sequence elements $k$: $k = 0$ for a finite number of $p_n$.
\end{Example}

We will use the above two identifications when constructing the Deligne categories; before that, we explain how to take an ultraproduct of categories.

\begin{Example}{}{}
  Let $\mathcal{C}_i$ be a collection of (locally small) categories. Then the ultraproduct $\mathcal{C} := \prod_{\mathcal{F}} \mathcal{C}_i$ is a category; $\Ob(\mathcal{C}) = \prod_{\mathcal{F}} \Ob(\mathcal{C}_i)$ and $\Hom(\prod_{\mathcal{F}} X_i, \prod_{\mathcal{F}} Y_i) = \prod_{\mathcal{F}} \Hom_{\mathcal{C}_i} (X_i, Y_i)$. Composition is given by the ultraproduct of composition maps (\cref{xmpl:ultra_operations}). The ultraproduct will also satisfy the axioms of a category by \Los{}' theorem, and if the $\mathcal{C}_i$ are abelian, monoidal, etc., $\mathcal{C}$ will also be abelian, monoidal, etc. However, if the $\mathcal{C}_i$ are tensor categories, $\mathcal{C}$ will have all the properties of a tensor category except local finiteness, as in \cref{xmpl:ultra_vec}.
\end{Example}

In general, $\mathcal{C}$ will be very large, and in practice we will be considering its subcategories, e.g. bounding the objects we consider (e.g. bounded length). One way to insert a boundedness condition is by trying to give our ultraproduct a ``filtration''.

\begin{Definition}{}{}
  Let $V_n$ a sequence of filtered vector spaces $F^0 V_n \subset F^1 V_n \subset \cdots \subset F^k V_n \subset \cdots$. Then the \textbf{filtered or restricted ultraproduct} $\prod^r_{\mathcal{F}} V_n = \bigcup_{k = 0}^\infty \prod_{\mathcal{F}} F^k V_n \subset \prod_{\mathcal{F}} V_n$.
\end{Definition}
\begin{Remark}{}{}
  For our purposes, we will work with filtered ultraproducts only of objects that have universally bounded multiplicities of simple objects. In this case, filtered ultraproducts will be independent of the filtration we give, so it makes sense to not have the filtration $F$ appear in the notation. This is a technical result that is proven in \cite{utiralova2022harishchandra}.
\end{Remark}

\begin{Example}{}{}
  Filtered ultraproducts can enforce desirable stability or finiteness conditions. For example, let $V$ a vector space of countable dimension and filter $V$ by finite-dimensional vector spaces. Then while $\prod_{\mathcal{F}} V$ has uncountable dimension, $\prod^r_{\mathcal{F}} V = \bigcup_{k = 0}^\infty \prod_{\mathcal{F}} F^k V = \bigcup_{k = 0}^\infty F^k V = V$.
\end{Example}

\section{Ultraproduct construction of the Deligne categories}

As mentioned previously, we can also construct $\Rep(S_t)$ as an ultraproduct of $\Rep(S_n)$; this allows us to think of $\Rep(S_t)$ as a ``limit'' of the classical categories. Moreover, it transfers many classical constructions relatively easily to the Deligne category case: roughly, we just take an ultraproduct of the classical constructions. We follow the argument from \cite{etingof_ddca}.

Fix a non-principal ultrafilter $\mathcal{F}$ on $\mathbb{N}$. Recall from \cref{sec:ultraproduct_bg} that we have an isomorphism $\prod_{\mathcal{F}} \overline{\mathbb{Q}} \cong \mathbb{C}$, so let $t = \prod_{\mathcal{F}} i$ for integer $i$. Now $\prod_{\mathcal{F}} \Rep(S_i)$ is very large, so we will want to consider a subcategory. Analogous to the above construction, we start with the permutation representation: we take the full subcategory generated by $\prod_{\mathcal{F}} V_i$ closed under tensor powers, direct sums, and direct summands, and let this be $\Rep(S_t)$.

However, there is an additional snag: we will also need to distinguish between transcendental and algebraic $t \in \mathbb{C}$ because we are taking an ultraproduct of integers (the integer $t$ case we get by taking $\prod_{\mathcal{F}} t$).
\begin{Lemma}{}{transcendental_ultraproduct}
  Let $t$ be a transcendental number over $\mathbb{C}$. Then for any sequence of $i \in \overline{\mathbb{Q}}$ where not almost all the $i$ are equal, there exists an isomorphism $\prod_{\mathcal{F}} \overline{\mathbb{Q}} \cong \mathbb{C}$ so that $\prod_{\mathcal{F}} i = t$.
\end{Lemma}
\begin{proof}
  By assumption, there are infinitely many distinct $i$. We first claim that $\overline{t} = \prod_{\mathcal{F}} i$ must be transcendental regardless of what isomorphism $\varphi: \prod_{\mathcal{F}} \overline{\mathbb{Q}} \cong \mathbb{C}$.
  Suppose not; then $t$ is algebraic over $\mathbb{Q}$ and there exists a monic polynomial $f$ with $f(t) = 0$. Then, by \cref{theo:los} for almost all $i$, we must have $f(i) = 0$ also. In other words, $f(i) = 0$ for infinitely many distinct $i$ and $f = 0$. Now take an automorphism of $\mathbb{C}$ over $\mathbb{Q}$ that sends $\overline{t} \mapsto t$; composing this automorphism with $\varphi$, we get an isomorphism of fields where $\prod_{\mathcal{F}} i = t$.
\end{proof}
To also construct $\Rep(S_t)$ for algebraic $t$, we will take an ultraproduct of algebraic closures of finite fields instead:
\begin{Lemma}{}{}
  Let $t$ be an algebraic number over $\mathbb{Q}$, but not an integer, with minimal polynomial $q(x) \in \mathbb{Z}[x]$. Then we can choose a sequence of distinct primes $p_n$, a sequence of integers $t_n$ tending to infinity with $q(t_n) = 0 \in F_{p_n}$, and an isomorphism $\prod_{\mathcal{F}} \overline{F}_{p_n} \cong \mathbb{C}$ so that $\prod_{\mathcal{F}} t_i = t$.
\end{Lemma}
\begin{proof}
  First, given any sequence of distinct primes $p_n$, we can find an infinite sequence of integers $t_n$ tending to infinity with $q(t_n) \equiv 0 \pmod{p_n}$. We will show the stronger statement that there are infinitely many primes dividing $q(n)$ below \cref{propn:q_primes}. This implies that the $t_n$ exist. Moreover, if the sequence of $t_n$ were bounded, some $q(t_n)$ would then be divisible by an infinite number of primes, contradiction. So the statement also implies that the $t_n$ tend to infinity.
  Fix any isomorphism $\varphi: \prod_{\mathcal{F}} \overline{F}_{p_n} \cong \mathbb{C}$. Now that we know that the $p_n$ and the $t_n$ exist, by \Los{}' theorem, $\overline{t} := \prod_{\mathcal{F}} t_n$ is a root of $q$, hence algebraic. As in the above lemma, we can then compose $\varphi$ with an automorphism of $\mathbb{C}$ sending $\overline{t} \mapsto t$, so we are done.

  \begin{Proposition}{}{q_primes}
    There are infinitely many primes dividing the numbers $q(n)$.
  \end{Proposition}
  \begin{proof}
    Suppose for contradiction that only $k$ primes divide the numbers $q(n)$. Consider $Q_N$, the number of integers of the form $q(n)$, $n$ a nonnegative integer, so that $q(n) < N$. There exists some constant $C$ so that $q(n) < C \cdot n^{\deg q}$ for all positive $n$. Then
    \begin{equation*}
      Q_N \ge \frac{1}{C} \cdot N^{1/\deg q}.
    \end{equation*}
    However, by our assumption, $Q_N$ is bounded above by the number of integers factorizable into $k$ fixed primes. Each prime is at least $2$, so
    \begin{equation*}
      Q_N \le \log_2(N)^k.
    \end{equation*}
    The first expression is polynomial in $N$ while the second is logarithmic, so if we pick $N$ large enough, we have a contradiction.
  \end{proof}
\end{proof}

Now we can construct $\Rep(S_t)$ for all $t \in \mathbb{C}$.
\begin{Theorem}{}{ultraproduct_st_equiv}
  \begin{enumerate}
  \item Let $t$ be a transcendental complex number or nonnegative integer. Fix an isomorphism $\varphi: \prod_{\mathcal{F}} \overline{\mathbb{Q}} \cong \mathbb{C}$ so that $\prod_{\mathcal{F}} n = t$. Let $V_n$ be the permutation representation of $\Rep(S_n)$ and $\widehat{\mathcal{C}} = \prod_{\mathcal{F}} \Rep(S_n, \overline{\mathbb{Q}})$. Then the full subcategory of $\widehat{\mathcal{C}}$ generated by $V_t := \prod_{\mathcal{F}} V_n$ and closed under tensor products, direct sums, and direct summands is equivalent to $\Rep(S_t)$ in a way consistent with $\varphi$.
  \item Let $t$ be an algebraic number, but not a nonnegative integer, with minimal polynomial $q(x) \in \mathbb{Z}[x]$. Fix a sequence of distinct primes $p_n$ and a sequence $t_n$ tending to infinity so that $q(t_n) \equiv 0 \pmod{p_n}$. Fix an isomorphism $\varphi = \prod_{\mathcal{F}} \overline{F}_{p_n} \cong \mathbb{C}$ so that $\prod_{\mathcal{F}} t_n = t$. Let $V_{t_n}$ be the permutation representation of $\Rep(S_n, \overline{F}_{p_n})$ and $\widehat{\mathcal{C}} = \prod_{\mathcal{F}} \Rep(S_{t_n}, \mathbb{F}_{p_n})$. Then the full subcategory of $\widehat{\mathcal{C}}$ generated by $V_t := \prod_{\mathcal{F}} V_{t_n}$ and closed under tensor products, direct sums, and direct summands is equivalent to $\Rep(S_t)$ in a way consistent with $\varphi$.
  \end{enumerate}
\end{Theorem}
\begin{proof}
  The arguments for each part are very similar, so we only give the proof of (a); for (b), replace each $n$ with $t_n$.
  Since $V_n$ is a symmetric Frobenius algebra, $V_t$ will be a symmetric Frobenius algebra of dimension $t$. Therefore, by the universal property of $\Rep(S_t)$ given above, we have a symmetric monoidal functor $F: \Rep(S_t) \to \widehat{\mathcal{C}}$ with $F([1]) = V_t$, and the image of $F$ is the full subcategory of $\widehat{\mathcal{C}}$ generated from $V_t$ by taking tensor products, direct sums, and direct summands. Therefore, $F$ is essentially surjective and we also want to prove that it induces an isomorphism of Hom-spaces, i.e. it is fully faithful.

  Because both categories are the Karoubian envelope of the additive completion of a category generated by tensor powers of a single object, it suffices to show that
  \begin{equation*}
    \prod_{\mathcal{F}} \Hom_{S_n} (V_n^{\otimes r}, V_n^{\otimes s}) = \Hom_{\Rep(S_t)} ([r], [s]).
  \end{equation*}
  Recall that in our construction of $\Rep(S_t)$, we showed that for $N$ large enough, $\Hom_{S_n} (V_n^{\otimes r}, V_n^{\otimes s}) = FP_{r, s}$, i.e. it was a $F$-vector space with basis corresponding to the partitions of $r + s$. Therefore, this equality is true as vector spaces, since the left side is an ultraproduct with almost all terms $FP_{r, s}$ and the right side is $FP_{r, s}$.

  It remains for us to show that the composition law is the same. Again by our construction of $\Rep(S_t)$, we had for partitions $Q$ and $P$ that $e_Q \circ e_P = n^{N(P, Q)} e_{P * Q}$. Thus, for the ultraproduct, composition is given by
  \begin{equation*}
    \left(\prod_{\mathcal{F}} e_{Q}\right) \circ \left(\prod_{\mathcal{F}} e_P\right) = \prod_{\mathcal{F}} n^{N(P, Q)} e_{P * Q} = t^{N(P, Q)} e_{P * Q}.
  \end{equation*}
  This is the same composition law as $\Rep(S_t)$, so we are done.
\end{proof}

Likewise, we can construct $\Rep(GL_t)$, $\Rep(O_t)$, and $\Rep(Sp_t)$ via the ultraproduct; we only show the equivalence of constructions for $\Rep(GL_t)$ because the proof for the latter two categories is almost exactly the same.
\begin{Theorem}{}{}
  \begin{enumerate}
  \item Let $t$ be a transcendental complex number or nonnegative integer. Fix an isomorphism $\varphi: \prod_{\mathcal{F}} \overline{\mathbb{Q}} \cong \mathbb{C}$ so that $\prod_{\mathcal{F}} n = t$. Let $V_n$ be the permutation representation of $\Rep(GL_n)$ and $\widehat{\mathcal{C}} = \prod_{\mathcal{F}} \Rep(GL_n, \overline{\mathbb{Q}})$. Then the full subcategory of $\widehat{\mathcal{C}}$ generated by $V_t := \prod_{\mathcal{F}} V_n$ and closed under duals, tensor products, direct sums, and direct summands is equivalent to $\Rep(GL_t)$ in a way consistent with $\varphi$.
  \item Let $t$ be an algebraic number, but not a nonnegative integer, with minimal polynomial $q(x) \in \mathbb{Z}[x]$. Fix a sequence of distinct primes $p_n$ and a sequence $t_n$ tending to infinity so that $q(t_n) \equiv 0 \pmod{p_n}$. Fix an isomorphism $\varphi = \prod_{\mathcal{F}} \overline{F}_{p_n} \cong \mathbb{C}$ so that $\prod_{\mathcal{F}} t_n = t$. Let $V_{t_n}$ be the permutation representation of $\Rep(GL_n, \overline{F}_{p_n})$ and $\widehat{\mathcal{C}} = \prod_{\mathcal{F}} \Rep(GL_{t_n}, \mathbb{F}_{p_n})$. Then the full subcategory of $\widehat{\mathcal{C}}$ generated by $V_t := \prod_{\mathcal{F}} V_{t_n}$ and closed under duals, tensor products, direct sums, and direct summands is equivalent to $\Rep(GL_t)$ in a way consistent with $\varphi$.
  \end{enumerate}
\end{Theorem}
\begin{proof}
  Again, we only give the proof of (a); for (b), replace $n$ with $t_n$. Since $\dim V_t = t$, by the universal property of $\Rep(GL_t)$, we have an additive symmetric monoidal functor $F: \Rep(GL_t) \to \widehat{\mathcal{C}}$ with $F([1, 0]) = V_t$. The image of $F$ is the full subcategory of $\widehat{\mathcal{C}}$ generated from $V_t$ by taking duals, tensor products, direct sums, and direct summands. So $F$ is essentially surjective.

  Now we want to show $F$ is fully faithful. Again, it is enough to show that the Hom-spaces agree for tensor powers of $V$, that is,
  \begin{equation*}
    \prod_{\mathcal{F}} \Hom_{GL_n}(V_n^{\otimes r} \otimes V_n^{* \otimes s}, V_n^{\otimes p} \otimes V_n^{* \otimes q}) = \Hom_{\Rep(GL_t)}([r, s], [p, q]).
  \end{equation*}
  We showed in \cref{sec:gl_construction} that $\Hom_{GL_n}(V_n^{\otimes r} \otimes V_n^{* \otimes s}, V_n^{\otimes p} \otimes V_n^{* \otimes q}) \ne 0$ only if $r + q = p + s$, and in that case, for $n$ large enough, it is $\mathbb{C}S_{r + q}$. So indeed, as vector spaces, the Hom-spaces are the same.

  It remains to show that the composition laws are the same. Let $\sigma, \sigma' \in \mathbb{C}S_{r + q}$. In the $GL_n$ case, for $n$ large enough, we had $e_\sigma \circ e_{\sigma'} = n^l e_{\sigma \circ \sigma'}$, where $l$ is the number of loops removed after concatenating diagrams. Likewise, $GL_t$ has $e_\sigma \circ e_{\sigma'} = t^l e_{\sigma \circ \sigma'}$. So, the ultraproduct has composition law
  \begin{equation*}
    \left(\prod_{\mathcal{F}} e_\sigma\right) \circ \left(\prod_{\mathcal{F}} e_{\sigma'} \right) = \prod_{\mathcal{F}} n^l e_{\sigma \circ \sigma'} = t^l e_{\sigma \circ \sigma'}.
  \end{equation*}
  The composition laws agree as well, so we are done.
\end{proof}

\subsection{Deligne Categories in Positive Characteristic}
This ultraproduct construction will be most relevant in \cref{chap:research}, but it is also important for us because it allows us to construct the Deligne categories over fields of positive characteristic. We only briefly discuss this since the work focuses on the categories in the characteristic $0$ case.

The construction in \cref{chap:interp_constructions} doesn't work because our endomorphism algebras are no longer semisimple, so we will no longer get a symmetric tensor category. On the other hand, it is quite straightforward to write down the generalization using ultraproducts.
\begin{Definition}{}{}
  Fix a prime $p$ and let
  \begin{equation*}
    \widehat{\mathcal{C}} \coloneq \prod_{\mathcal{F}} \Rep(G_n; \overline{\mathbb{F}_p})
  \end{equation*}
  for $G$ one of $S, GL, O, Sp$. $\widehat{\mathcal{C}}$ will be a symmetric tensor category over $\mathbb{K} := \prod_{\mathcal{F}} \mathbb{F}_p$, which will be some algebraically closed field of characteristic $p$. Let $\mathcal{C}_{\mathcal{F}}$ be the full symmetric subcategory $\widehat{\mathcal{C}}$ generated by $V := \prod_{\mathcal{F}} V_n$, where $V_n$ is the permutation representation, and closed under direct sums, tensor products, and direct summands (and duals if $G = GL$). $\mathcal{C}_{\mathcal{F}}$ generalizes the Deligne categories to positive characteristic.
\end{Definition}
Unlike in the characteristic $0$ case, the dimension of $V$, which is $t := \prod_{\mathcal{F}} n$, is dependent on the choice of ultrafilter. Namely, $t \in \mathbb{F}_p$ is the value of $n \pmod{p}$ for almost all $n$. In fact, the dimensions of all objects $X \in \Ob(\mathcal{C}_{\mathcal{F}})$ will lie in $\mathbb{F}_p$. In addition, even if we fix $t$, $\mathcal{C}_{\mathcal{F}}$ is still determined by the choice of ultrafilter. It turns out $\mathcal{C}_{\mathcal{F}}$ depends only on the numbers $\wedge^{p^k} V$ for all $k$, a result shown in \cite{harman_stability}. The properties of this category are still an open question, but we direct the reader to \cite{etingof_padic} and \cite{etingof_stc} for further reading.


\chapter{Relating Deligne categories to the classical categories}
In this chapter, we describe some properties of the categories we constructed in the previous chapter to show that they reflect many properties of the classical representation categories. The results about $\Rep(S_t)$ are drawn from \cite{comes_ostrik}, and the results about $\Rep(GL_t)$ from \cite{comes_wilson}.

\section{Indecomposable objects}
Here we show that the indecomposable objects of the Deligne categories are classified by partitions, just as in the classical case. In both the symmetric and matrix group cases, we will describe the primitive idempotents of the endomorphism algebras (partition algebras, (walled) Brauer algebras) and show that these are in bijections with all (bi)partitions. In the classical case, for integer $n$, the irreducible representations correspond to partitions that are in some way bounded by $n$; but this will not be the case for the Deligne categories, even when $t$ is an integer, something we will discuss in the next section. Thus, this encourages thinking of Deligne categories as limits of the classical representation categories, as in the ultraproduct construction.

\subsection{Simple objects of $\Rep(S_t)$}

Recall from \cref{def:karoubi} and from the general theory of additive completions that we have the following facts about indecomposables in $\Rep(S_t)$ over any field $F$:
\begin{Proposition}{}{}
  \begin{enumerate}
  \item Every indecomposable object in $\Rep(S_t)$ is isomorphic to an object $([n], e)$ where $e$ is a primitive idempotent in $\Hom([n], [n]) =: FP_n(t)$.
  \item An object $([n], e)$ is indecomposable iff $e$ is primitive.
  \item If $e, e'$ are two primitive idempotents, then $([n], e) \cong ([n], e')$ iff $e, e'$ are conjugate idempotents.
  \end{enumerate}
  In other words, the primitive idempotents up to conjugation correspond to indecomposable objects of $\Rep(S_t)$.
\end{Proposition}

In order to classify these idempotents up to conjugation, we use a general result about idempotents in algebras:
\begin{Lemma}{}{}
  Let $A$ be a finite-dimensional $F$-algebra and $\mathbf{p} \in A$ an idempotent. Let $(\mathbf{p})$ denote the two-sided ideal in $A$ generated by $\mathbf{p}$. Then we have a correspondence between primitive idempotents in $A$ and the disjoint union
  \begin{equation*}
    \{\text{primitive idempotents of } A/(\mathbf{p})\} \sqcup \{\text{primitive idempotents in } \mathbf{p} A \mathbf{p}\}
  \end{equation*}
  (all up to conjugation). The correspondence works as follows:

  Let $e \in A$ be a primitive idempotent. If $e \in (\mathbf{p})$, it corresponds to an element in $\mathbf{p} A \mathbf{p}$. If $e \notin (\mathbf{p})$, then it corresponds to its equivalence class in $A/(\mathbf{p})$.
\end{Lemma}

Now we work only over a field $F$ with characteristic $0$ and classify idempotents of $FP_n(t)$.
We will choose $\mathbf{p}$ to be the following element of $FP_n(t)$ when $n > 1$, which is easily seen to be an idempotent:
\begin{center}
  \begin{tikzpicture}
    \filldraw[black] (0,0) circle (2pt)
    (0, 1) circle (2pt)
    (1.5, 0) circle (2pt)
    (1.5, 1) circle (2pt)
    (4.5, 0) circle (2pt)
    (4.5, 1) circle (2pt)
    (6, 0) circle (2pt)
    (6, 1) circle (2pt)
    (7.5, 0) circle (2pt)
    (7.5, 1) circle (2pt);
    \draw (0,0) node[below=1pt] {$1'$} -- (0,1) node[above=1pt] {$1$};
    \draw (1.5,0) node[below=1pt] {$2'$} -- (1.5,1) node[above=1pt] {$2$};
    \draw (4.5,0) node[below=1pt] {$(n-2)'$} -- (4.5,1) node[above=1pt] {$n-2$};      
    \draw (6,0) node[below=1pt] {$(n-1)'$} -- (6,1) node[above=1pt] {$n-1$} -- (7.5,1) node[above=1pt] {$n$} -- (7.5,0) node[below=1pt] {$n'$} -- cycle;
    \draw (3,0.5) node {$\cdots$};
  \end{tikzpicture}
\end{center}
\begin{Lemma}{}{idem_class}
  The chosen $\mathbf{p}$ has the following isomorphisms, which will allow us to inductively classify idempotents:
  \begin{enumerate}
  \item $\mathbf{p} FP_n(t) \mathbf{p} \cong FP_{n - 1}(t)$
  \item $FP_n(t)/(\mathbf{p}) \cong FS_n$.
  \end{enumerate}
\end{Lemma}
\begin{proof}
  ~\begin{enumerate}
  \item We can treat $FP_{n - 1}(t)$ as a subalgebra of $FP_n(t)$, namely as the subalgebra generated by the partitions of $(n, n)$ so that $\{n - 1, n\}$ are in the same part and $\{(n - 1)', n'\}$ are in the same part. But this is exactly $\mathbf{p} FP_n(t) \mathbf{p}$. Suppose $\pi \in FP_n(t)$ corresponds to a partition $P$ of $(n, n)$. Then $\mathbf{p} \pi \mathbf{p}$ corresponds to a partition refining $P$, so that the parts containing $n - 1$ and $n$ are merged, as are the parts containing $(n - 1)'$ and $n'$. Moreover, every $\pi \in FP_{n - 1}(t)$ (thought of as a subalgebra of $FP_n(t)$) can be obtained this way, since $\mathbf{p} \pi \mathbf{p} = \pi$.
  \item We can treat $FS_n$ as a subalgebra of $FP_n(t)$ by identifying every $\sigma \in S_n$ with the partition $\{\{1, \sigma(1)'\}, \dots, \{n, \sigma(n)'\}\}$, then taking the subalgebra generated by said partitions. Then $FS_n \cap (\mathbf{p}) = \{0\}$ (see the previous part, for example) and it suffices to show that elements of $P_{n, n}$ lie in exactly one of $(\mathbf{p})$ and $S_n$.

    First, we define some auxiliary partitions:
    \begin{equation*}
      \mu_{j, k} = \vcenter{\hbox{
          \begin{tikzpicture}
            \def \m{1.2};
            \foreach \i [evaluate=\i as \x using \i*\m] in {0, 2, 3, 4, 6, 7, 8, 10} {
              \filldraw[black] (\x, 0) circle (2pt) -- (\x, 1) circle (2pt);
            }
            \draw (0, 0) node[below=1pt] {$1'$};
            \draw (0, 1) node[above=1pt] {1};
            \draw ({\m * 3}, 0) node[below=1pt] {$j'$};
            \draw ({\m * 3}, 1) node[above=1pt] {$j$};
            \draw ({\m * 7}, 0) node[below=1pt] {$k'$};
            \draw ({\m * 7}, 1) node[above=1pt] {$k$};
            \draw ({\m * 10}, 0) node[below=1pt] {$n'$};
            \draw ({\m * 10}, 1) node[above=1pt] {$n$};
            \draw ({\m * 3}, 0) .. controls ({\m * 4}, 0.5) and ({\m * 6}, 0.5) .. ({\m * 7}, 0);
            \draw ({\m * 3}, 1) .. controls ({\m * 4}, 0.5) and ({\m * 6}, 0.5) .. ({\m * 7}, 1);            
          \end{tikzpicture}
        }
      }
    \end{equation*}
    and
    \begin{equation*}
      \nu_{i, j} = \vcenter{\hbox{
          \begin{tikzpicture}
            \def \m{1.2};
            \foreach \i [evaluate=\i as \x using \i*\m] in {0, 2, 4, 6, 7, 8, 10} {
              \filldraw[black] (\x, 0) circle (2pt) -- (\x, 1) circle (2pt);
            }
            \filldraw[black] ({\m * 3}, 0) circle (2pt);
            \filldraw[black] ({\m * 3}, 1) circle (2pt);
            \draw({\m * 1}, 0.5) node {$\cdots$};
            \draw (0, 0) node[below=1pt] {$1'$};
            \draw (0, 1) node[above=1pt] {1};
            \draw ({\m * 3}, 0) node[below=1pt] {$i'$};
            \draw ({\m * 3}, 1) node[above=1pt] {$i$};
            \draw ({\m * 7}, 0) node[below=1pt] {$j'$};
            \draw ({\m * 7}, 1) node[above=1pt] {$j$};
            \draw({\m * 9}, 0.5) node {$\cdots$};
            \draw ({\m * 10}, 0) node[below=1pt] {$n'$};
            \draw ({\m * 10}, 1) node[above=1pt] {$n$};
            \draw ({\m * 3}, 0) .. controls ({\m * 4}, 0.5) and ({\m * 6}, 0.5) .. ({\m * 7}, 0);
          \end{tikzpicture}
        }
      }.
    \end{equation*}
    
    Namely, $\mu_{j, k}$ is $\sigma \mathbf{p} \sigma$ where $\sigma = (j, n - 1)(k, n)$: right (resp. left) multiplication by $\mu_{j, k}$ merges the parts of $j$ and $k$ (resp. $j'$ and $k'$). The action of $\nu_{i, j}$ is a little more complicated: if $P_i$ is the part containing $i$ of a partition, then right-multiplying by $\nu_{i, j}$ forms a new partition with parts $\{i\}, P_i \setminus \{i\} \cup P_j$. Left-multiplying by $\nu_{i, j}$ forms a partition with parts $\{i'\} \cup P_{j'}, P_{i'} \setminus \{i'\}$ (possibly with a factor of $t$ if $P_{i'} = \{i'\}$).
    So suppose that $\pi \in P_{n, n} \setminus S_n$. Then $\pi$ falls into one of two categories:
    \begin{enumerate}
    \item $\pi$ has a part $\{i\}$ with $i \in \{1, \dots, n\}$. In this case, $\pi = \pi \mu_{i, j} \nu_{i, j}$.
    \item $\pi$ has $P_j = P_k$ for some $j, k \in \{1, \dots, n\}$: in other words, there exist $j, k$ lying in the same part of $\pi$. In that case, $\pi = \pi \mu_{j, k}$.
    \end{enumerate}
    In both cases, $\pi \in (\mathbf{p})$ because $\mu \in (\mathbf{p})$, so we are done.
  \end{enumerate}
\end{proof}

Therefore, the primitive idempotents of $FP_n(t)$ correspond to the primitive idempotents of $FP_{n - 1}(t)$ and $FS_n$. The primitive idempotents of $FS_n$ are already known, e.g. from the representation theory of $S_n$, to correspond to Young diagrams of weight $n$. So once we know about $FP_1(t)$, we can inductively classify the primitive idempotents of $FP_n(t)$. Moreover, a primitive idempotent $e$ corresponds to a Young diagram of weight $n$ iff $e \notin (\mathbf{p})$ and $[e] \in FP_n(t) / (\mathbf{p}) \cong FS_n$ is a primitive idempotent.
\begin{Theorem}{}{sn_idem}
  We split into two cases, when $t \ne 0$ and when $t = 0$.
  \begin{enumerate}
  \item If $t \ne 0$, primitive idempotents in $FP_n(t)$ correspond to Young diagrams $\lambda$ with $|\lambda| \le n$.
  \item If $t = 0$ and $n > 0$, primitive idempotents in $FP_n(t)$ correspond to Young diagrams $\lambda$ with $0 < |\lambda| \le n$.
  \end{enumerate}
\end{Theorem}
\begin{proof}
  ~\begin{enumerate}
  \item When $n = 0$, $FP_0(t) = F$ and there is only one idempotent, so the statement holds. When $n = 1$, $1$ is still an idempotent, but it decomposes into two primitive idempotents. Let $\pi \in P_{1, 1}$ be the partition $\{[1], [1']\}$; it is clear that $f := \frac{1}{t}\pi$ is an idempotent, and therefore so is $1 - f$. It is easy to see (recalling that $FP_1(t)$ has basis $1$ and $\pi$) that these are the only two primitive idempotents up to conjugation, so the statement is also true for $n = 1$. Now by \cref{lem:idem_class} and the fact stated above that the primitive idempotents of $FS_n$ correspond to Young diagrams $\lambda$ with $|\lambda| = n$, we are done by induction.
  \item When $n = 1$, $FP_1(0) = F[\pi]/(\pi^2)$ since $\pi^2 = 0 \cdot \pi$. Therefore, $1$ is the only idempotent and the statement is true for $n = 1$. Now the induction is the same as above.
  \end{enumerate}
\end{proof}

Now we can classify indecomposable objects of $\Rep(S_t)$.
\begin{Theorem}{}{sn_indecomp}
  Let $\lambda$ be a Young diagram with weight $m$ and $e_\lambda \in FP_m(t)$ be a primitive idempotent (unique up to conjugation) corresponding to $\lambda$, as constructed in \cref{theo:sn_idem}. If $m = t = 0$, we set $e_{\emptyset} = \id_0 \in FP_0(0)$. So $e_\lambda$ corresponds to an indecomposable object $L(\lambda) := ([m], e_\lambda)$, defined up to isomorphism. For $n \ge 0$, this map $\lambda \mapsto L(\lambda)$ is a bijection between Young diagrams $\lambda$ so that $0 \le |\lambda| \le n$ and (nonzero) indecomposable objects in $\Rep(S_t)$ of the form $([m], e)$ for $m \le n$, which has the following properties:
  \begin{enumerate}
  \item If $0 < |\lambda| \le n$, then there exists an idempotent $e \in FP_n(t)$ with $([n], e) \cong L(\lambda)$.
  \item If $t \ne 0$, there exists an idempotent $e \in FP_n(t)$ with $([n], e) \cong L(\emptyset)$ corresponding to the trivial empty partition of weight $0$.
  \item If $t = 0$, $([0], \id_0)$ is the unique object of the form $([m], e)$ isomorphic to $L(\emptyset)$.
  \end{enumerate}
\end{Theorem}
\begin{proof}
  Again, we induct on $n$. First, consider the $t \ne 0$ case. If $n = 0$, we know $\End([0]) = F$, so we are done.

  Now consider $n = 1$. We saw in \cref{theo:sn_idem} that $f, 1 - f$ were the only two primitive idempotents in $FP_1(t)$, so we know that $([1], f)$ and $([1], 1 - f)$ are non-isomorphic indecomposable objects. We claim that these form a complete set of pairwise non-isomorphic indecomposable objects by constructing an isomorphism $([0], \id_0) \cong ([1], f)$. Let $\mu \in P_{0, 1}$ and $\mu' \in P_{1, 0}$ be the only partitions of $1, 0$ and $0, 1$ respectively. The map $f\mu\id_0: ([0], \id_0) \to ([1], f)$ has inverse $\frac{1}{t}\id_0 \mu' f$, so we have the desired isomorphism. This shows the theorem for $n = 1$.

  Now let's prove it for $n > 1$. For every Young diagram $\lambda$ with $|\lambda| = n$, $\{L(\lambda) | |\lambda| = n\}$ is a set of pairwise non-isomorphic indecomposables. The remaining idempotents of $FP_n(t)$ must correspond to Young diagrams $\lambda$ with $|\lambda| < n$. By induction, for every Young diagram $\lambda$ with $|\lambda| < n$, there is a primitive idempotent $f_\lambda \in FP_{n - 1}(t)$ so that $\{([n - 1], f_\lambda) | 0 \le |\lambda| < n\}$ is the set of all the pairwise non-isomorphic indecomposable objects of the form $([m], e)$ when $m < n$. We want to find idempotents $\widetilde{f}_\lambda \in FP_n(t)$ so that $([n], \widetilde{f}_\lambda) \cong ([n - 1], f_\lambda)$. We will set $\widetilde{f}_\lambda = \varphi_n f_\lambda \varphi'_n$ where
  \begin{equation*}
    \varphi_n = \vcenter{\hbox{
        \begin{tikzpicture}
          \def \m{1.5};
          \foreach \i [evaluate=\i as \x using \i*\m] in {0, 1, 2, 3} {
            \filldraw[black] (\x, 0) circle (2pt) -- (\x, 1) circle (2pt);
          }
          \filldraw[black]({\m * 4}, 0) circle (2pt);
          \draw (0, 0) node[below=1pt] {$1'$};
          \draw (0, 1) node[above=1pt] {$1$};
          \draw ({\m * 1}, 0) node[below=1pt] {$2'$};
          \draw ({\m * 1}, 1) node[above=1pt] {$2$};
          \draw ({\m * 1.5}, 0.5) node {$\cdots$};
          \draw ({\m * 2}, 0) node[below=1pt] {$(n - 2)'$};
          \draw ({\m * 2}, 1) node[above=1pt] {$n - 2$};
          \draw ({\m * 3}, 0) node[below=1pt] {$(n - 1)'$};
          \draw ({\m * 3}, 1) node[above=1pt] {$n - 1$};

          \draw ({\m * 4}, 0) node[below=1pt] {$n'$} -- ({\m * 3}, 1);
        \end{tikzpicture}
      }
    },
    \varphi'_n = \vcenter{\hbox{
        \begin{tikzpicture}
          \def \m{1.5};
          \foreach \i [evaluate=\i as \x using \i*\m] in {0, 1, 2, 3} {
            \filldraw[black] (\x, 0) circle (2pt) -- (\x, 1) circle (2pt);
          }
          \filldraw[black]({\m * 4}, 1) circle (2pt);
          \draw (0, 0) node[below=1pt] {$1'$};
          \draw (0, 1) node[above=1pt] {$1$};
          \draw ({\m * 1}, 0) node[below=1pt] {$2'$};
          \draw ({\m * 1}, 1) node[above=1pt] {$2$};
          \draw ({\m * 1.5}, 0.5) node {$\cdots$};
          \draw ({\m * 2}, 0) node[below=1pt] {$(n - 2)'$};
          \draw ({\m * 2}, 1) node[above=1pt] {$n - 2$};
          \draw ({\m * 3}, 0) node[below=1pt] {$(n - 1)'$};
          \draw ({\m * 3}, 1) node[above=1pt] {$n - 1$};

          \draw ({\m * 4}, 1) node[above=1pt] {$n$} -- ({\m * 3}, 0);
        \end{tikzpicture}
      }
    }.
  \end{equation*}
  Then because $\varphi'_n\varphi_n = \id_{n - 1}$, it follows that $\widetilde{f}_\lambda$ is an idempotent and that
  \begin{equation*}
  f_\lambda \varphi'_n \widetilde{f}_\lambda: ([n], \widetilde{f}_\lambda) \to ([n - 1], f_\lambda)
\end{equation*}
  is an isomorphism with inverse $\widetilde{f}_\lambda \varphi_n f_\lambda$. Therefore, $([n], \widetilde{f}_\lambda)$ is indecomposable and so $\widetilde{f}_\lambda$ is primitive in $FP_n(t)$. Moreover, by construction $\widetilde{f}_\lambda = \mathbf{p} \widetilde{f}_\lambda \in (\mathbf{p})$. So we have found a complete set of pairwise non-conjugate primitive idempotents in $FP_n(t)$ that correspond to Young diagrams of weight at most $n$, namely
  \begin{equation*}
    \{\widetilde{f}_\lambda | 0 \le |\lambda| < n\} \cup \{e_\lambda | |\lambda| = n\}.
  \end{equation*}
  Thus we also have a complete set of pairwise non-isomorphic indecomposable objects of the form $([m], e)$ with $m \le n$, namely
  \begin{equation*}
    \{([n], \widetilde{f}_\lambda) | 0 \le |\lambda| < n\} \cup \{([n], e_\lambda) = L(\lambda) | |\lambda| = n\}.
  \end{equation*}

  The $t = 0$ case is very similar: for finding irreducibles corresponding to $0 < |\lambda| \le n$, the proof is the same as above. For the uniqueness of the indecomposable corresponding to $\lambda = \emptyset$, $([0], \id_0)$, notice that every composition $([0], \id_0) \to ([m], e) \to ([0], \id_0)$ in $\Rep(S_0)$ must be $0$ unless $m = 0$, $e = \id_0$: multiplying an element of $P_{0, m}$ and an element of $P_{m, k}$ always has at least one factor of $t = 0$ if $m > 0$.
\end{proof}
\begin{Corollary}{}{}
  The map $\lambda \mapsto L(\lambda)$ induces a bijection between Young diagrams of arbitrary size and nonzero indecomposable objects of $\Rep(S_t)$ up to isomorphism.
\end{Corollary}
\begin{Remark}{}{t_young_diagram}
  If $t = n \in \mathbb{Z}_+$ and $\lambda$ is a partition with $|\lambda| + \lambda_1 \le n$, the indecomposable object $L(\lambda)$ can be thought of as the $S_n$ representation corresponding to the Young diagram with top row of size $n - |\lambda|$. Thus, we should think of $L(\lambda)$ as corresponding to the Young diagram with a top row of ``length'' $t - |\lambda|$ and $\lambda$ below.
\end{Remark}

\subsection{Indecomposable objects of $\Rep(GL_t)$}
Our analysis of the indecomposable objects of $\Rep(GL_t)$ is analogous to the previous part on $\Rep(S_t)$, except that we will now consider bipartitions rather than partitions and $B_{r, s}(t)$ instead of $FP_n(t)$. We first establish some notation.
\begin{Definition}{}{}
  A \textbf{bipartition} $\lambda$ is a pair of partitions $(\lambda^\bullet, \lambda^\circ)$. We define $|\lambda| = (|\lambda^\bullet|, |\lambda^\circ|)$.
\end{Definition}
\begin{Proposition}{}{}
  The bipartitions $\lambda$ with $|\lambda| = (r, s)$ index primitive idempotents (up to conjugation) of $\mathbb{C}(S_r \times S_s)$.
\end{Proposition}
The next result is the analogue of \cref{theo:sn_idem} in the previous section.
\begin{Theorem}{}{}
  \begin{enumerate}
  \item If $t \ne 0$ or $r \ne s$, primitive idempotents in $B_{r, s}$ correspond to bipartitions $\lambda$ so that $|\lambda| = (r - i, s - i)$ for $0 \le i \le \min(r, s)$.
  \item If $t = 0$ and $r > 0$, primitive idempotents in $B_{r, s}$ correspond to bipartitions $\lambda$ so that $|\lambda| = (r - i, r - i)$ for $0 \le i < r$.
  \end{enumerate}
\end{Theorem}
\begin{proof}
  The proof is also analogous, but involves different choices of $\mathbf{p}$, $\varphi$, and $\varphi'$, and we will need to split into $t = 0$ and $t \ne 0$ cases.
  
  To prove the analogue of \cref{lem:idem_class}: in the $t \ne 0$ case, we pick
  \begin{equation*}
    \mathbf{p} = \frac{1}{t}\vcenter{\hbox{
        \begin{tikzpicture}
          \def \m{1.2};
          \foreach \i [evaluate=\i as \x using \i*\m] in {0, 2, 3, 4, 6, 7} {
            \ifnum \i<4
            \ifnum \i<3
            \filldraw[black] (\x, 0) circle (2pt) -- (\x, 1) circle (2pt);
            \else
            \filldraw[black] (\x, 0) circle (2pt);
            \filldraw[black] (\x, 1) circle (2pt);
            \fi
            \else
            \ifnum \i<7
            \draw (\x, 0) circle (2pt) -- (\x, 1) circle (2pt);
            \else
            \draw (\x, 0) circle (2pt);
            \draw (\x, 1) circle (2pt);
            \fi
            \fi
          }
          \draw (0, 0) node[below=1pt] {$1'$};
          \draw (0, 1) node[above=1pt] {1};
          
          \draw ({\m * 1}, 0.5) node {$\cdots$};
          
          \draw ({\m * 2}, 0) node[below=1pt] {$(r - 1)'$};
          \draw ({\m * 2}, 1) node[above=1pt] {$r - 1$};
          
          \draw ({\m * 3}, 0) node[below=1pt] {$r'$};
          \draw ({\m * 3}, 1) node[above=1pt] {$r$};

          \draw ({\m * 4}, 0) node[below=1pt] {$1'$};
          \draw ({\m * 4}, 1) node[above=1pt] {$1$};

          \draw ({\m * 5}, 0.5) node {$\cdots$};
          
          \draw ({\m * 6}, 0) node[below=1pt] {$(s - 1)'$};
          \draw ({\m * 6}, 1) node[above=1pt] {$s - 1$};
          
          \draw ({\m * 7}, 0) node[below=1pt] {$s'$};
          \draw ({\m * 7}, 1) node[above=1pt] {$s$};

          \foreach \x in {0, 1} {
            \draw ({\m * 3}, \x) .. controls ({\m * 4}, 0.5) and ({\m * 6}, 0.5) .. ({\m * 7}, \x);
          }
        \end{tikzpicture}
      }}
  \end{equation*}
  If $t = 0$, if one of $r$ or $s$ is at least $2$, WLOG $s$, we define
  \begin{equation*}
    \mathbf{p} = \vcenter{\hbox{
        \begin{tikzpicture}
          \def \m{1.2};
          \foreach \i [evaluate=\i as \x using \i*\m] in {0, 2, 3, 4, 6, 7, 8} {
            \ifnum \i<4
            \ifnum \i<3
            \filldraw[black] (\x, 0) circle (2pt) -- (\x, 1) circle (2pt);
            \else
            \filldraw[black] (\x, 0) circle (2pt);
            \filldraw[black] (\x, 1) circle (2pt);
            \fi
            \else
            \ifnum \i<7
            \draw (\x, 0) circle (2pt) -- (\x, 1) circle (2pt);
            \else
            \draw (\x, 0) circle (2pt);
            \draw (\x, 1) circle (2pt);
            \fi
            \fi
          }
          \draw (0, 0) node[below=1pt] {$1'$};
          \draw (0, 1) node[above=1pt] {1};
          
          \draw ({\m * 1}, 0.5) node {$\cdots$};
          
          \draw ({\m * 2}, 0) node[below=1pt] {$(r - 1)'$};
          \draw ({\m * 2}, 1) node[above=1pt] {$r - 1$};
          
          \draw ({\m * 3}, 0) node[below=1pt] {$r'$};
          \draw ({\m * 3}, 1) node[above=1pt] {$r$};

          \draw ({\m * 4}, 0) node[below=1pt] {$1'$};
          \draw ({\m * 4}, 1) node[above=1pt] {$1$};

          \draw ({\m * 5}, 0.5) node {$\cdots$};
          
          \draw ({\m * 6}, 0) node[below=1pt] {$(s - 2)'$};
          \draw ({\m * 6}, 1) node[above=1pt] {$s - 2$};
          
          \draw ({\m * 7}, 0) node[below=1pt] {$(s - 1)'$};
          \draw ({\m * 7}, 1) node[above=1pt] {$s - 1$};

          \draw ({\m * 8}, 0) node[below=1pt] {$s'$};
          \draw ({\m * 8}, 1) node[above=1pt] {$s$};
          
          \draw ({\m * 3}, 1) .. controls ({\m * 4}, 0.5) and ({\m * 6}, 0.5) .. ({\m * 7}, 1);
          \draw ({\m * 3}, 0) .. controls ({\m * 4}, 0.5) and ({\m * 7}, 0.5) .. ({\m * 8}, 0);
          \draw ({\m * 8}, 1) -- ({\m * 7}, 0);
        \end{tikzpicture}
      }}
  \end{equation*}
  $(\mathbf{p})$ is spanned by diagrams with less than $r + s$ edges between dots of the same color, and $B_{r, s}(t) / (\mathbf{p}) = \mathbb{C}(S_r \times S_s)$.

  Now the proof is similar to that of \cref{theo:sn_idem}.
\end{proof}

\begin{Theorem}{}{}
  For every bipartition $\lambda$ with weight $(r, s)$, define $L(\lambda) := ([r, s], e_\lambda)$ where $e_\lambda \in B_{r, s}$ is a primitive idempotent, unique up to conjugation, corresponding to $\lambda$. The assignment $\lambda \mapsto L(\lambda)$ induces a bijection between bipartitions of arbitrary size and nonzero indecomposable objects of $\Rep(GL_t)$, up to isomorphism.
\end{Theorem}
\begin{proof}
  Also analogous to \cref{theo:sn_indecomp}, but we need different definitions of $\varphi$ and $\varphi'$ and need to separate the $t = 0$ and $t \ne 0$ cases.

  If $t \ne 0$, we will choose
  \begin{equation*}
    \varphi_{r, s} = \vcenter{\hbox{
        \begin{tikzpicture}
      \foreach \x in {0, 2} {
        \filldraw[black] (\x, 0) circle (2pt) -- (\x, 1) circle (2pt);
      }
      \foreach \x in {5, 7} {
        \draw (\x, 0) -- (\x, 1);
        \filldraw[white] (\x, 0) circle (2pt);
        \filldraw[white] (\x, 1) circle (2pt);
        \draw (\x, 0) circle (2pt);
        \draw (\x, 1) circle (2pt);
      }
      \foreach \x in {1, 6} {
        \draw (\x, 0.5) node {$\cdots$};
      }
      \draw (4, 0) arc (0:180:0.5);
      \filldraw[black] (3, 0) circle (2pt);
      \filldraw[white] (4, 0) circle (2pt);
      \draw (4, 0) circle (2pt);
      \draw[decorate,decoration={brace,amplitude=10pt},xshift=0pt,yshift=4pt] (0, 1) -- (2, 1) node[midway,yshift=20pt] {$r - 1$};
      \draw[decorate,decoration={brace,amplitude=10pt},xshift=0pt,yshift=4pt] (5, 1) -- (7, 1) node[midway,yshift=20pt] {$s - 1$};
    \end{tikzpicture}
  }
}
  \end{equation*}
  and
  \begin{equation*}
    \varphi'_{r, s} = \vcenter{\hbox{
    \begin{tikzpicture}
      \foreach \x in {0, 2} {
        \filldraw[black] (\x, 0) circle (2pt) -- (\x, 1) circle (2pt);
      }
      \foreach \x in {5, 7} {
        \draw (\x, 0) -- (\x, 1);
        \filldraw[white] (\x, 0) circle (2pt);
        \filldraw[white] (\x, 1) circle (2pt);
        \draw (\x, 0) circle (2pt);
        \draw (\x, 1) circle (2pt);
      }
      \foreach \x in {1, 6} {
        \draw (\x, 0.5) node {$\cdots$};
      }
      \draw (3, 1) arc (180:360:0.5);
      \filldraw[black] (3, 1) circle (2pt);
      \filldraw[white] (4, 1) circle (2pt);
      \draw (4, 1) circle (2pt);
      \draw[decorate,decoration={brace,mirror,amplitude=10pt},xshift=0pt,yshift=-4pt] (0, 0) -- (2, 0) node[midway,yshift=-20pt] {$r - 1$};
      \draw[decorate,decoration={brace,mirror,amplitude=10pt},xshift=0pt,yshift=-4pt] (5, 0) -- (7, 0) node[midway,yshift=-20pt] {$s - 1$};
    \end{tikzpicture}
  }
  }.
  \end{equation*}
  Here we see that $\varphi'_{r, s}\varphi_{r, s} = t \id_{r - 1, s - 1}$, so an idempotent $f_\lambda$ of $[r - 1, s - 1]$ corresponds to the idempotent $\frac{1}{t} \varphi_{r, s} f_\lambda \varphi'_{r, s}$. If $t = 0$, we use the same $\varphi'$, but we take
  \begin{equation*}
    \varphi_{r, s} =
    \begin{cases}
      \vcenter{\hbox{
          \begin{tikzpicture}
            \foreach \x in {0, 2, 3} {
              \filldraw[black] (\x, 0) circle (2pt);
              \filldraw[black] (\x, 1) circle (2pt);
            }
            \filldraw[black] (4, 0) circle (2pt);
            \draw (1, 0.5) node {$\cdots$};
            \draw (0, 1) -- (0, 0);
            \draw (2, 1) -- (2, 0);
            \draw (4, 0) -- (3, 1);
            \draw (5, 0) arc (0:180:1 and 0.5);
            \filldraw[white] (5, 0) circle (2pt);
            \draw (5, 0) circle (2pt);                  \draw[decorate,decoration={brace,amplitude=10pt},xshift=0pt,yshift=4pt] (0, 1) -- (2, 1) node[midway,yshift=20pt] {$r - 2$};
          \end{tikzpicture}
        }} & s = 0, r > 0 \\
      \vcenter{\hbox{
          \begin{tikzpicture}
            \foreach \x in {0, 2} {
              \filldraw[black] (\x, 0) circle (2pt);
              \filldraw[black] (\x, 1) circle (2pt);
            }
            \filldraw[black] (4, 0) circle (2pt);
            \foreach \x in {1, 7} {
              \draw (1, 0.5) node {$\cdots$};
            }
            \filldraw[black] (3, 0) circle (2pt);
            \draw (4, 0) -- (4, 1);
            \draw (5, 0) arc (0:180:1 and 0.5);
            \filldraw[white] (5, 0) circle (2pt);
            \draw (5, 0) circle (2pt);
            \foreach \x in {4, 6, 8} {
              \draw (\x, 0) -- (\x, 1);
              \filldraw[white] (\x, 0) circle (2pt);
              \draw (\x, 0) circle (2pt);
                            \filldraw[white] (\x, 1) circle (2pt);
              \draw (\x, 1) circle (2pt);
            }
                  \draw[decorate,decoration={brace,amplitude=10pt},xshift=0pt,yshift=4pt] (0, 1) -- (2, 1) node[midway,yshift=20pt] {$r - 1$};
      \draw[decorate,decoration={brace,amplitude=10pt},xshift=0pt,yshift=4pt] (6, 1) -- (8, 1) node[midway,yshift=20pt] {$s - 2$};
          \end{tikzpicture}
        }} & s > 0.
    \end{cases}
  \end{equation*}
  Now here $\varphi'_{r, s}\varphi_{r, s} = \id_{r - 1, s - 1}$, so an idempotent $f_\lambda$ of $[r - 1, s - 1]$ corresponds to the idempotent $ \varphi_{r, s} f_\lambda \varphi'_{r, s}$.
  The rest follows like the proof of \cref{theo:sn_indecomp}.
\end{proof}

We can likewise prove similar theorems for $\Rep(O_t)$ and $\Rep(Sp_t)$ by studying the Brauer algebras, whose structure is discussed in \cite{brauer_visscher}; we have
\begin{Theorem}{}{}
  The simple objects of $\Rep(O_t)$ (and therefore $\Rep(Sp_t)$) are labeled by all partitions $\lambda$ as well.
\end{Theorem}

\begin{Remark}{}{}
  As one might expect, we can also construct the indecomposable objects of all these categories by taking ultraproducts of irreducible representations in the classical categories. For example, consider $\Rep(S_t)$. Fix a non-principal ultrafilter $\mathcal{F}$ on $\mathbb{N}$ and either an isomorphism $\prod_{\mathcal{F}} \overline{\mathbb{Q}} \cong \mathbb{C}$ or $\prod_{\mathcal{F}} \overline{\mathbb{F}}_{p_n} \cong \mathbb{C}$ depending on whether $t$ is algebraic or not, so that $t$ corresponds to $\prod_{\mathcal{F}} t_n$ where $t_n \to \infty$. Let $\lambda$ be a partition, which we know corresponds to an indecomposable object $W_\lambda \in \Ob(\Rep(S_t))$, and define $\lambda|_{t_n}$ to be the partition $(t_n - |\lambda|, \lambda_1, \dots, \lambda_{l(\lambda)})$; this is well-defined for almost all $t_n$. Then
  \begin{equation*}
    W_\lambda = \prod_{\mathcal{F}} V_{\lambda|_{t_n}}
  \end{equation*}
  where $V_{\lambda|_{t_n}}$ is the irreducible representation in $S_{t_n}$ corresponding to the partition $\lambda|_{t_n}$. This follows from \cref{theo:ultraproduct_st_equiv}, which shows that the basis elements of $\mathbb{C}P_n(t)$ are ultraproducts of basis elements of $\mathbb{C}P_n(t_n)$. Therefore, the idempotents of $\mathbb{C}P_n(t)$ are the ultraproducts of the same idempotents of $\mathbb{C}P_n(t_n)$ for almost all $n$.
\end{Remark}

\section{Semisimplicity}

We have classified the indecomposable objects above; when $t$ is not an integer, the Deligne categories are semisimple, so the above construction actually lists all the simple objects.

\begin{Theorem}{}{}
  If $t \notin \mathbb{Z}_+$, then $\Rep(S_t)$ is a semisimple, hence abelian, category.
\end{Theorem}
\begin{proof}
  This follows from the fact that $\mathbb{C}P_n(t)$ is semisimple for $t \notin \mathbb{Z}_+$, proved in \cite[Section 2]{halverson_ram}. Suppose $t \notin \mathbb{Z}_+$. Let $\lambda, \mu$ be two partitions so that $\Hom(L(\lambda), L(\mu)) \ne 0$; by \cref{theo:sn_indecomp} we can pick an $n$ so that $L(\lambda) \cong ([n], e_1)$ and $L(\mu) \cong ([n], e_2)$. By definition, then, $\Hom(L(\lambda), L(\mu)) = e_2 \mathbb{C}P_n(t) e_1 \ne 0$. Now the semisimplicity of $\mathbb{C}P_n(t)$ implies that $e_2, e_1$ are conjugate, so $L(\lambda) \cong L(\mu)$ and $\lambda = \mu$. Moreover, for $\lambda$ a partition of $n$, $\End(L(\lambda)) = e_\lambda \mathbb{C}P_n(t) e_\lambda$ is a division algebra (again by semisimplicity of $\mathbb{C}P_n(t)$), so $\Rep(S_t)$ is indeed semisimple.
\end{proof}
\begin{Remark}{}{}
  In the next section, we will discuss what happens in the $t \in \mathbb{Z}_+$ case.
\end{Remark}

We have similar results on $\Rep(GL_t)$, $\Rep(O_t)$, and $\Rep(Sp_t)$:
\begin{Theorem}{}{}
  \begin{itemize}
  \item If $t \notin \mathbb{Z}$, then $\Rep(GL_t)$ is an abelian semisimple category.
  \item If $t \notin \mathbb{Z}$, then $\Rep(O_t)$ and $\Rep(Sp_t)$ are abelian semisimple categories.
  \end{itemize}
\end{Theorem}
\begin{proof}
  $B_{r, s}(t)$ and $B_r(t)$ are semisimple when $t \notin \mathbb{Z}$, proved in \cite[Section 6]{cddm} and \cite{rui}. Then we can make the same argument as above to show that $\Rep(GL_t)$, $\Rep(O_t)$, and $\Rep(Sp_t)$ are semisimple in this case.
\end{proof}

\section{Recovering the classical categories} \label{sec:interpolation}
In this section, we relate the Deligne categories when $t \in \mathbb{Z}$ to the classical categories. To distinguish these, the Deligne categories will be denoted by $\Rep$, while the classical categories will be denoted by $\overline{\Rep}$. When $t$ is an integer (or positive integer, in the symmetric group case), the Deligne categories are not abelian semisimple. However, for a Deligne category $\mathcal{C}_t$ and its classical counterpart, recall that we can use the universal properties described in \cref{sec:univ_prop_st} and \cref{sec:glosp_interp} to write down the interpolation functor defined in \cref{def:interp_functor}. This tensor functor in fact kills the ``negligible morphisms'' - those with trace zero - in the Deligne category.

How does the interpolation functor act on other objects besides $[m] \in \Ob(\Rep(S_n))$? It suffices to define $\mathcal{F}$ on objects of the form $([m], e)$ where $e$ is an idempotent. Recall from \cref{propn:t_compn} the definition of $T_P$ as the composition of iterated multiplication and its dual. Since we are working with representations of $S_n$, we just need to use the correspondence in \cref{xmpl:correspondence}, where we demonstrated how elements of partition algebras correspond to maps between tensor powers of $V$. Denote this correspondence via $e \in \Hom([l], [m]) \mapsto f_e : V^{\otimes l} \to V^{\otimes m}$. Then
\begin{equation*}
  \mathcal{F}([m], e) = f_e V^{\otimes m}.
\end{equation*}
A morphism $\alpha: ([l], d) \to ([m], e)$ is of the form $e \beta d$ for $\beta \in \Hom([l], [m])$. So
\begin{equation*}
  \mathcal{F}(\alpha) = f_e f_\beta f_d.
\end{equation*}

\begin{Proposition}{}{}
  $\mathcal{F}$ is surjective on objects and morphisms.
\end{Proposition}
\begin{proof}
  We know the image of $\mathcal{F}$ contains $V$. Because $\mathcal{F}$ is an additive symmetric monoidal functor between Karoubian categories, the image must be closed under taking tensor products, duals, direct sums, and direct summands. So the image of $\mathcal{F}$ is in fact all of $\overline{\Rep}(S_n)$ and $\mathcal{F}$ is surjective on objects.

  To see that $\mathcal{F}$ is surjective on morphisms, recall that the basis we wrote down for the partition algebra $FP_{l, m}$ always spans $\Hom(V^{\otimes l}, V^{\otimes m})$, so the map $e \mapsto f_e$ is surjective.
\end{proof}

However, $\mathcal{F}$ is not an equivalence of categories. One way to see this is to go back to our construction of the correspondence between $FP_{l, m}$ and $\Hom(V^{\otimes l}, V^{\otimes m})$. If a partition $P$ of $l + m$ has more than $n$ parts, which can occur whenever $l + m > n$, it does not correspond to an orbit of $S_n$. Then at least $f_{\delta_P} = 0 : V^{\otimes l} \to V^{\otimes m}$. To precisely describe by how much $\mathcal{F}$ is not an equivalence of categories, we want to describe what morphisms $\mathcal{F}$ sends to $0$: these will be the ``negligible morphisms.''

\begin{Definition}{}{}
  A morphism $f: X \to Y$ in a symmetric tensor category $\mathcal{C}$ is \textbf{negligible} if $\Tr(fg) = 0$ for all $g: Y \to X$. Let $\mathcal{N}(X, Y)$ denote all the negligible morphisms in $\Hom(X, Y)$ and $\mathcal{N}$ denote the collection of $\mathcal{N}(X, Y)$ for all $X, Y$.
\end{Definition}

\begin{Proposition}{}{propn}
  $\mathcal{N}$ is a tensor ideal. That is, the following properties hold:
  \begin{enumerate}
  \item if $f \in \mathcal{N}(X, Y)$ and $g \in \Hom(Z, X)$, then $f \circ g \in \mathcal{N}(Z, Y)$.
  \item if $f \in \mathcal{N}(X, Y)$ and $g \in \Hom(Y, Z)$, then $g \circ f \in \mathcal{N}(X, Z)$.
    In other words, $\mathcal{N}$ ``absorbs'' composition.
  \end{enumerate}
  \begin{enumerate}[resume]
  \item if $f \in \mathcal{N}(X, Y)$ and $g \in \Hom(W, Z)$, then $f \otimes g \in \mathcal{N}(X \otimes W, Y \otimes Z)$. So $\mathcal{N}$ absorbs tensor products.
  \end{enumerate}
\end{Proposition}
\begin{proof}
  ~\begin{enumerate}
  \item Take $h \in \Hom(Y, Z)$. Then $\Tr((f \circ g) \circ h) = \Tr(f \circ (g \circ h)) = 0$.
  \item Take $h \in \Hom(Z, X)$. We repeatedly use the fact that $\Tr(A \circ B) = \Tr(B \circ A)$:
    \begin{align*}
      \Tr((g \circ f) \circ h) &= \Tr(h \circ (g \circ f)) \\
                               &= \Tr((h \circ g) \circ f) \\
                               &= \Tr(f \circ (h \circ g)) \\
                               &= 0.
    \end{align*}
  \item Take $h \in \Hom(Y \otimes Z, X \otimes W)$. Recall that we have an isomorphism $\Xi: \Hom(X, Y) \cong \Hom(\mathbf{1}, Y \otimes X^*)$, described in \cref{eq:dual_hom}. Then $\Tr((f \otimes g) \circ h) = \Tr(f \circ H)$, where $H \in \Hom(Y, X)$ corresponds to the map in $\Hom(\mathbf{1}, X \otimes Y^*)$ given by
    \begin{equation*}
      \mathbf{1} \xrightarrow{\Xi(h) \otimes \Xi(g)} X \otimes W \otimes Y^* \otimes Z^* \otimes Z \otimes W^* \xrightarrow{\ev_{Z} \otimes \ev_{W}} Y^* \otimes X.
    \end{equation*}
    We omit the diagram chasing for proving the correspondence. Thus, $f \otimes g$ is also negligible.
  \end{enumerate}
\end{proof}

\begin{Definition}{}{}
  As one might expect from the terminology, we can define the quotient $\mathcal{C}/\mathcal{J}$ with $\mathcal{C}$ a tensor category and $\mathcal{J}$ a tensor ideal
  \begin{itemize}
  \item $\Ob(\mathcal{C}/\mathcal{J}) = \Ob(\mathcal{C})$
  \item $\Hom_{\mathcal{C}/\mathcal{J}}(X, Y) = \Hom_{\mathcal{C}}(X, Y) / \mathcal{J}(X, Y)$.
  \end{itemize}
\end{Definition}
This quotient will also be a tensor category; moreover, if $\mathcal{C}$ was a Karoubian category, then so will the quotient - for a finite-dimensional algebra $A$, an idempotent of a quotient $A'$ of $A$ will lift to an idempotent of $A$.

\begin{Theorem}{}{}
  The functor $\Rep(S_n) / \mathcal{N} \to \overline{\Rep}(S_n)$ induced by the interpolation functor $\mathcal{F}$ is an equivalence of categories.
\end{Theorem}
\begin{proof}
  We claim that semisimple categories have no nonzero negligible morphisms, proved in \cref{propn:semisimple_negligible} below. Now $\mathcal{F}$ is full, and it preserves traces because it is a tensor functor. Thus a morphism $f$ is negligible iff $\mathcal{F}(f)$ is negligible, and the morphisms $\mathcal{F}$ maps to $0$ are exactly the negligible morphisms.
\end{proof}

\begin{Proposition}{}{semisimple_negligible}
  If a symmetric tensor category $\mathcal{C}$ is semisimple, it has no nonzero negligible morphisms.
\end{Proposition}
\begin{proof}
  If $f \ne 0 \in \mathcal{N}(X, Y)$, then the bilinear form $(f, g) \mapsto \Tr(fg)$ is degenerate. So it suffices to show that this bilinear form is actually nondegenerate for every nonzero $f$ in $\mathbb{C}$. This pairing is equivalent to the pairing
  \begin{align*}
    \Hom(\mathbf{1}, X^* \otimes Y) \times \Hom(\mathbf{1}, X \otimes Y^*) &\to \mathbf{1} \\
    (\Xi(f), \Xi(g)) &\mapsto \ev_{X \otimes Y^*} \circ (\Xi(f) \otimes \Xi(g))
  \end{align*}
  which follows from the definition of trace and of $\Xi$. So, it suffices to show that for all $Z \in \Ob(\mathcal{C})$, the pairing
  \begin{align*}
    \Hom(\mathbf{1}, Z) \times \Hom(\mathbf{1}, Z^*) &\to \mathbf{1} \\
    (f, g) &\mapsto \ev_Z \circ (f \otimes g)
  \end{align*}
  is nondegenerate. But by semisimplicity of $\mathcal{C}$, we can decompose $Z = \mathbf{1}^{\oplus k} \oplus \overline{Z}$ where $\overline{Z}$ contains no copies of $\mathbf{1}$. So $\Hom(\mathbf{1}, Z) = \Hom(\mathbf{1}, Z^*) = F^k$ and $\Hom(\mathbf{1}, \overline{Z}) = \Hom(\overline{Z}, \mathbf{1}) = 0$. Thus the above pairing is actually just the dot product on $F^k$, given on basis vectors by $(e_i, e_j) \mapsto \delta_{ij}$, which is clearly nondegenerate.
\end{proof}
\begin{Remark}{}{}
  After quotienting, some objects will become isomorphic to $0$, namely $X$ so that the only morphisms into and out of $X$ are zero morphisms. In this case, $\Tr \id_X = 0$ where $\id_X \in \End_{\overline{\Rep}(S_n)}(X)$. That means that $\id_X \in \End_{\Rep(S_n)}(X)$ must have been negligible, and in particular $\Tr(\id_X \circ \id_X) = 0$ and $X$ has dimension $0$ in $\Rep(S_n)$. The simple objects that get annihilated are the $L(\lambda)$ so that $|\lambda| + \lambda_1 > n$: see \cref{rmk:t_young_diagram}.
\end{Remark}
\begin{Remark}{}{}
  This construction, namely quotienting a symmetric tensor category $\mathcal{C}$ by its negligible morphisms, is known as the semisimplification $\overline{\mathcal{C}}$ of $\mathcal{C}$. The semisimplification is semisimple, abelian, and its simple objects are the indecomposable objects of $\mathcal{C}$ of nonzero dimension. See \cite{semisimplification} for more on this construction.
\end{Remark}

\subsection{Interpolation functors for $\Rep(GL_t)$, $\Rep(O_t)$, and $\Rep(Sp_{t})$}

Unlike the case of $\Rep(S_t)$, the Deligne categories from matrix groups also have interpolation functors for $t = n \in \mathbb{Z}_{< 0}$. We describe what happens in the $\Rep(GL_t)$ case, but supergroups are outside the scope of this paper, so we direct the reader to \cite{comes_wilson} to see proofs.

Consider the category of representations of $GL_{|n|}$ with a $\mathbb{Z}/2\mathbb{Z}$-grading given by the action of $-I \in GL_{|n|}$: that is, every representation has an ``odd'' part (where $-I$ acts by $-\id$) and an ``even'' part (where $-I$ acts by $\id$), which we denote by $\underline{\Rep}(GL_{|n|}, \mathbb{Z}/2\mathbb{Z})$. For example, the permutation representation is completely odd.

\begin{Theorem}{}{}
Let $\mathcal{N}$ be the negligible morphisms of $\Rep(GL_{t = n})$, as defined above. We have $\Rep(GL_{t = n})/\mathcal{N} \cong \underline{\Rep}(GL_{|n|}, \mathbb{Z}/2\mathbb{Z})$ via the functor $\Rep(GL_{t = n}) \to \underline{\Rep}(GL_{|n|}, \mathbb{Z}/2\mathbb{Z})$ induced by the universal property of $\Rep(GL_{t = n})$.
\end{Theorem}

We can likewise recover $\underline{\Rep}(O_n)$ and $\underline{\Rep}(Sp_{2n})$ by quotienting out negligible morphisms.


\chapter{Further reading and research}
\label{chap:research}
In this section, we survey some of the further literature on the subject and describe a joint work with Alexandra Utiralova (\cite{utiralova2022harishchandra}).

Deligne categories were first introduced in the papers of Deligne and Milne, namely \textit{Tannakian categories} (\cite{deligne_milne}) and \textit{La categorie des repr{\'e}sentations du groupe sym{\'e}trique $S_t$, lorsque $t$ n'est pas un entier naturel} (\cite{deligne_st}). While we discussed the parts of \cite{comes_ostrik} and \cite{comes_wilson} on indecomposable objects of these categories, the papers also classify blocks of the Deligne categories and generalize from complex number parameters to indeterminate ones. \cite{comes_ideals} classifies the tensor ideals of $\Rep(GL_t)$. Comes and Ostrik described the structure of the abelian envelope of $\Rep(S_t)$ in \cite{rep_ab_sd}, and Entova, Hinich, and Serganova did the same for $\Rep(GL_t)$ in \cite{gl_ab_env}. Knop generalized Deligne's constructions in \cite{knop_semisimple} and \cite{knop_tensor_envelopes}, constructing families of semisimple tensor categories that interpolated other representation categories.

Etingof, in his papers \cite{etingof_cr1} and \cite{etingof_cr2}, explained how to transfer problems from classical representation theory to complex rank. The ultraproduct construction is very useful in this setting. For example, \cite{aizenbud_sw} describes Schur-Weyl duality for Deligne categories and \cite{hk_simple_algebras} classifies simple algebras for $\Rep(S_t)$. Etingof also proposed ways to study infinite-dimensional representation theory in complex rank via looking at representations of ind-objects in Deligne categories, and there are now  many papers on the subject, e.g. \cite{kalinov}, \cite{entova_cherednik}, \cite{sciarappa_sca}, \cite{utiralova_characters}, and \cite{pakharev_weylkac}.

One example of this ``problem transfer'' is the definition and investigation of Harish-Chandra bimodules in Deligne categories. We first give some definitions and motivation for their study following the exposition in \cite{humphreys}, then present some results in \cite{utiralova_characters} and \cite{utiralova2022harishchandra}. The material here requires more background than is assumed in the rest of the thesis.

\section{Definitions}

\begin{Definition}{}{}
  Let $\mathfrak{g}$ be a semisimple Lie algebra and let $\mathfrak{g}_d \subset \mathfrak{g} \oplus \mathfrak{g}^{op}$ be the diagonal copy of $\mathfrak{g}$, embedded via the map $x \mapsto (x, -x)$. A \textbf{Harish-Chandra bimodule} $M$ is a finitely generated $U(\mathfrak{g})$-bimodule with a locally finite action of $\mathfrak{g}_d$. The local finiteness condition is equivalent to saying that the adjoint action of $\mathfrak{g}$ integrates to a $G$-action.
\end{Definition}

Now we give some motivation for the study of these bimodules. Suppose we have a complex Lie algebra $\mathfrak{g}$ which is the complexification of the Lie algebra of some semisimple Lie group $G$ with maximal compact subgroup $K$. Let $\mathfrak{l}$ be the complexified Lie algebra of $K$. Consider a unitary representation of $G$ on a Hilbert space; the $K$-finite vectors form a $\mathfrak{g}$-module, and the action of $\mathfrak{g}_d$ is locally finite. This module has the properties of a Harish-Chandra module (finitely generated, $\mathfrak{l}$-finite).

The bimodules arise when we begin with a complex group $G$ with Lie algebra $\mathfrak{g}$, but treat it as real; then the complexification of the Lie algebra is $\mathfrak{g} \times \mathfrak{g}$ with universal enveloping algebra $U(\mathfrak{g}) \otimes U(\mathfrak{g})$. Let $U := U(\mathfrak{g})$ for convenience of notation. The center of $U \otimes U$ is $Z(\mathfrak{g}) \otimes Z(\mathfrak{g})$. $U$ has an antiautomorphism (induced by the transpose map) that fixes $Z(\mathfrak{g})$ and turns left $U$-modules into right $U$-modules, so every $\mathfrak{g} \times \mathfrak{g}$-module has the structure of a $U \otimes U$-bimodule. Now $\mathfrak{l}$ is a copy of $\mathfrak{g}$ via the $x \mapsto (x, -x)$ embedding given above.

In the classical case, all Harish-Chandra bimodules are of \textbf{finite $K$-type}, that is, every simple $\mathfrak{g}$-module appears finitely many times. However, this condition is not necessarily true in Deligne categories.

We set some notation we use throughout.
\begin{Definition}{}{}
  Let $G_t$ be one of $GL_t$, $O_t$, or $Sp_t$, and $G_n$ be $GL_n$, $O_{2n + 1}$ or $Sp_{2n}$ respectively, with Lie algebra $\mathfrak{g}_n$. Let $\mathcal{C}_t$ be $\Rep(G_t)$.
\end{Definition}
\begin{Definition}{}{}
  Let $\rho$ be the half-sum of the positive roots of $G_n$, which has root system $X_n$ for $X$ type $A$, $B$, $C$, or $D$, with Weyl group $W$. Let $E_n$ be an $n$-dimensional vector space with orthonormal basis $e_i$; we can realize the root systems as
  \begin{align*}
    A_{n - 1} &= \{e_i - e_j\}_{i, j} \\
    B_n &= \{e_i \pm e_j\}_{i, j} \cup \{e_i\}_i \\
    C_n &= \{e_i \pm e_j\}_{i, j} \cup \{2e_i\}_i \\
    D_n &= \{e_i \pm e_j\}_{i, j}.
  \end{align*}
\end{Definition}
Now we define Harish-Chandra bimodules in Deligne categories. Fix a non-principal ultrafilter $\mathcal{F}$ on $\mathbb{N}$, let $\mathbb{k} := \overline{\mathbb{Q}}$, and fix an isomorphism $\prod_{\mathcal{F}} \mathbb{k} \cong \mathbb{C}$. First, we need to define the universal enveloping algebra and describe some properties of its center.
\begin{Definition}{}{}
  The Lie algebra of $G_t$ is $\mathfrak{g}_t := \prod_{\mathcal{F}} \mathfrak{g}_n$; recall that it will inherit the bracket operation, which will still satisfy the Jacobi identity. The Lie algebra $\mathfrak{g}_t$ also inherits an action on objects of $\mathcal{C}_t$, given by the ultraproduct of $\mathfrak{g}_n$ actions. As objects of their corresponding categories, we have
  \begin{align*}
    \mathfrak{gl}_t &= V \otimes V^* \\
    \mathfrak{o}_t &= \wedge^2 V \\
    \mathfrak{sp}_t &= S^2 V
  \end{align*}
  where $V$ is the generating object of $\mathcal{C}_t$.
\end{Definition}
\begin{Definition}{}{}
  The universal enveloping algebra $U(\mathfrak{g}_t)$ is the filtered ultraproduct of $U(\mathfrak{g}_n)$ with respect to the PBW filtration $\prod^r_{\mathcal{F}} U(\mathfrak{g}_n)$; that is, we define
  \begin{equation*}
    F^k U(\mathfrak{g}_t) = \prod_{\mathcal{F}} F^k U(\mathfrak{g}_n)
  \end{equation*}
  and
  \begin{equation*}
    U(\mathfrak{g}_t) = \colim_k F^k U(\mathfrak{g}_t).
  \end{equation*}
  Note that $U(\mathfrak{g}_t)$ is not an object in $\mathcal{C}_t$, but rather in its ind-completion.

  The center $Z(U(\mathfrak{g}_t))$ is likewise the filtered ultraproduct of $Z(U(\mathfrak{g}_n))$ with respect to the PBW filtration.
\end{Definition}
For each $n$, we have the Harish-Chandra isomorphism $Z(U(\mathfrak{g}_n)) \cong \mathbb{k}[E_n]^W$, where the filtration on the right side is given by polynomial degree. We will use this isomorphism to characterize $Z(U(\mathfrak{g}_t))$.

If $G_n = GL_n$, it has Weyl group $S_n$, so $Z(U(\mathfrak{g}_n))$ is isomorphic to the ring of symmetric polynomials in $n$ variables. Then, the power sums $p_k := \sum x_i^k$ correspond to central elements $C_k \in Z(U(\mathfrak{g}_n))$ which generate $Z(U(\mathfrak{g}_n))$; in other words, $Z(U(\mathfrak{g}_n)) = \mathbb{k}[C_1, \dots, C_n]$. Moreover, the $C_k$ act on the irreducible $\mathfrak{g}_n$-module $L(\nu)$ with highest weight $\nu$ by
\begin{equation*}
  C_k|_{L(\nu)} = \sum (\nu_i + \rho_i)^k - \rho_i^k.
\end{equation*}
Therefore,
\begin{equation*}
  Z(U(\mathfrak{gl}_t)) = \prod^r_{\mathcal{F}} \mathbb{k}[C_1, \dots, C_n] \cong \mathbb{C}[C_1, C_2, \dots]
\end{equation*}
and $\deg C_k = k$.

If $G_n = O_{2n + 1}$ or $Sp_{2n}$, it has Weyl group $(\mathbb{Z} / 2\mathbb{Z})^n \rtimes S_n$, so $Z(U(\mathfrak{g}_n))$ is isomorphic to the symmetric polynomials in $n$ variables that are invariant under negating variables. In other words, it is generated by the even power sums $p_{2k}$, corresponding to a generating set $C_{2k}$ for $Z(U(\mathfrak{g}_n))$, and $Z(U(\mathfrak{g}_n)) = \mathbb{k}[C_2, \dots, C_{2n}]$. Again, the $C_{2k}$ act on $L(\nu)$ by $\sum (\nu_i + \rho_i)^{2k} - \rho_i^{2k}$, and
\begin{equation*}
  Z(U(\mathfrak{g}_t)) = \mathbb{C}[C_2, C_4, \dots].
\end{equation*}
We can thus also talk about central characters in the Deligne category setting.
\begin{Definition}{}{}
A central character of $U(\mathfrak{g}_t)$ is a homomorphism $\chi: U(\mathfrak{g}_t) \to \mathbb{C}$; by definition, it is an ultraproduct of characters $\chi^{(n)} : Z(U(\mathfrak{g}_n)) \to \mathbb{C}$ with $\chi(C_k) = \prod_{\mathcal{F}} \chi^{(n)}(C_k)$. If $\mathfrak{g} = \mathfrak{o}, \mathfrak{sp}$, we set $\chi(C_k) = 0$ for odd $k$.
\end{Definition}
\begin{Definition}{}{}
  We will say that a central character $\chi$ of $U(\mathfrak{g}_n)$ \textbf{corresponds} to a weight $\lambda$ if $Z(U(\mathfrak{g}_n))$ acts on $M_\lambda$, the Verma module with highest weight $\lambda - \rho$, by $\chi$.
\end{Definition}
\begin{Note}{}{}
  This is not a one-to-one correspondence - the center can act on multiple modules with the same central character - but this notation will be useful later on.
\end{Note}
\begin{Definition}{}{exp_cc}
  For ease of notation, we will also associate with each central character $\chi$ an \textbf{exponential central character} $\chi(z)$ defined as
  \begin{equation*}
    \chi(z) := \sum_i \frac{\chi(C_k)}{k!}z^k.
  \end{equation*}
In the $\mathfrak{g} = \mathfrak{o}, \mathfrak{sp}$ cases, it's more accurate to name this a ``hyperbolic cosine character.''
\end{Definition}
\begin{Definition}{}{}
  For a fixed central character $\chi$, we can define a module
  \begin{equation*}
  U_\chi := U(\mathfrak{g}_t) / \langle z - \chi(z) | z \in Z(U(\mathfrak{g}_t)) \rangle,
\end{equation*}
  i.e. we quotient the universal enveloping algebra by the ideal generated by the kernel of $\chi$.
\end{Definition}
Now we can define Harish-Chandra bimodules in $\mathfrak{g}_t$. It is analogous to the classical definition we gave above.
\begin{Definition}{}{}
  A Harish-Chandra bimodule of $\mathfrak{g}_t$ is a finitely generated $\mathfrak{g}_t \oplus \mathfrak{g}_t^{op}$-module $M \in \Ind \mathcal{C}_t$ so that $\mathfrak{g}_d$ acts naturally on $M$ and both copies of $Z(U(\mathfrak{g}_t))$ act locally finitely. We will use $\mathfrak{g}_t^{(l)}$ and $\mathfrak{g}_t^{(r)}$ to denote the right and left actions of $\mathfrak{g}_t$, respectively. By finitely generated, we mean that $M$ is a quotient of $(U(\mathfrak{g}_t) \otimes U(\mathfrak{g}_t))^{op} \otimes X$ for $X \in \Ob(\mathcal{C}_t)$. By a locally finite action of $Z(U(\mathfrak{g}_t))$, we mean that it has a finite filtration on whose successive quotients $Z(U(\mathfrak{g}_t))$ acts as a scalar.
\end{Definition}

For Harish-Chandra bimodules, central characters correspond to pairs of central characters $(\chi, \psi)$ of $Z$, as they  are homomorphisms $Z(U(\mathfrak{g}_t)) \otimes Z(U(\mathfrak{g}_t)) \to \mathbb{C}$.
\begin{Definition}{}{}
  Let $HC(\mathfrak{g}_t)$ denote the full subcategory of Harish-Chandra bimodules of $\mathfrak{g}_t$, and let $HC_{\chi, \psi}(\mathfrak{g}_t)$ denote the full subcategory of $HC(\mathfrak{g}_t)$ so that the left copy of $Z(U(\mathfrak{g}_t))$ acts by $\chi$ and the right by $\psi$.
\end{Definition}
\begin{Example}{}{}
Each $U_\chi$ is a Harish-Chandra bimodule; by definition, both copies of $Z(U(\mathfrak{g}_t))$ act via $\chi$, so $U_\chi \in \Ob(HC_{\chi, \chi}(\mathfrak{g}_t))$. However, it doesn't have finite $K$-type, as shown in \cite[Remark~3.21]{etingof_cr2}.
\end{Example}
\begin{Example}{}{}
  We can also come up with some objects in $HC_{\chi, \psi}(\mathfrak{g}_t)$. Let $X \in \Ob(\mathcal{C})$; then $U_\psi \otimes X$ is a $\mathfrak{g}_t$-bimodule, where $\mathfrak{g}_t$ acts on the left via left multiplication on $U_\psi$ and on the right via minus multiplication on the right on $U_\psi$ (this makes the action natural in $X$). Then, $(U_\psi \otimes X)_\chi$, i.e. the quotient $(U_\psi \otimes X) / \langle z - \chi(z) | z \in Z(U(\mathfrak{g}_t)) \rangle$, is an object in $HC_{\chi, \psi}(\mathfrak{g}_t)$. For ease of notation, we will call this module $N(\chi, \psi, X)$.
\end{Example}

\section{Central characters of Harish-Chandra bimodules}
We have the following result describing when, given a central character $(\chi, \psi)$, there exist Harish-Chandra bimodules with this central character.
\begin{Theorem}{}{central_characters}
  If $\mathfrak{g}_t = \mathfrak{gl}_t$, then $HC_{\chi, \psi}(\mathfrak{g}_t)$ is nonzero iff there exist complex numbers $b_1, \dots, b_r, c_1, \dots, c_s$ so that
  \begin{equation*}
    \chi(z) - \psi(z) = (e^z - 1)\left(\sum_{i = 1}^r e^{b_i z} - \sum_{i = 1}^s e^{c_i z}\right).
  \end{equation*}
  If $\mathfrak{g}_t = \mathfrak{o}_t$ or $\mathfrak{sp}_t$, then $HC_{\chi, \psi}(\mathfrak{g}_t)$ is nonzero iff there exist complex numbers $b_1, \dots, b_r$ so that
  \begin{equation*}
    \chi(z) - \psi(z) = \sum_{i = 1}^r (\cosh((b_i + 1)z) - \cosh(b_i z)) = 2\sinh\left(\frac{z}{2}\right) \sum_{i = 1}^r \sinh\left(\frac{(2b_i + 1)z}{2}\right).
  \end{equation*}
\end{Theorem}
To tackle this problem, we make liberal use of ultraproducts to work first in the classical case, then lift the results we find to $\mathfrak{g}_t$. First we show that $HC_{\chi, \psi}$ is nonzero iff $N(\chi, \psi, [r, s]) \ne 0$ if $\mathfrak{g} = \mathfrak{gl}$ and $N(\chi, \psi, [r]) \ne 0$ otherwise. This will allow us to analyze $N(\chi, \psi, V)$ where $V$ is the generating object (analog of the permutation representation).

\begin{Lemma}{}{}
  $HC_{\chi, \psi} \ne 0$ iff $N(\chi, \psi, X) \ne 0$ for some $X \in \Ob(\mathcal{C})$.
\end{Lemma}
\begin{proof}
  Suppose that $M \in HC_{\chi, \psi} \ne 0$. The bimodule $M$ is finitely generated, i.e. it is a quotient of the $U(\mathfrak{g}_t)$-bimodule $(U(\mathfrak{g_t}) \otimes U(\mathfrak{g}_t)^{op}) \otimes X$ for $X \in \Ob(\mathcal{C})$. We also have
  \begin{equation*}
    N(\chi, \psi, X) \coloneq U_\chi \otimes_{U(\mathfrak{g}_t^{(l)})} \otimes (U_\psi \otimes X) \cong (U_\chi \otimes U_\psi^{op}) \otimes_{U(\mathfrak{l}_t)} X.
  \end{equation*}
  So we have a natural map $N(\chi, \psi, X) \to M$ which is surjective because $X$ generates $M$, in which case $N(\chi, \psi, X) \ne 0$.
\end{proof}
\begin{Corollary}{}{}
  $HC_{\chi, \psi} \ne 0$ iff $N(\chi, \psi, [r, s]) \ne 0$ if $\mathfrak{g} = \mathfrak{gl}$, and $N(\chi, \psi, [r]) \ne 0$ otherwise.
\end{Corollary}
How does the center act on $N(\chi, \psi, [r, s]) = (U_\psi \otimes [r, s])_\chi$ (or $N(\chi, \psi, [r]$)? By the above analysis, a central element $C$ acts as $(C \otimes 1)_{U_\chi \otimes [r, s]}$ on the left and as $\Delta(C)_{U_\chi, [r, s]}$ on the right, and likewise for $[r]$. So it will be helpful to somehow relate $\Delta(C)$ and $C \otimes 1$, which we will first do for just $N(\chi, \psi, V) = (U_\psi \otimes V)_\chi$, where $V$ is the permutation representation. Using the ultraproduct construction, it suffices to look at the classical case. Let $\chi$ now be a central character of $Z(U(\mathfrak{g}_n))$. There, $U_\chi \hookrightarrow \End(M_\lambda)$ where $\lambda = (\lambda_1, \dots, \lambda_n)$ and $M_\lambda$ is the Verma module with highest weight $\lambda - \rho$; the map is injective by Duflo's theorem. Then we also have an injective map
\begin{equation*}
  U_\chi \otimes V \to \Hom(M_\lambda, M_\lambda \otimes V)
\end{equation*}
where $\mathfrak{g}^{(r)}$ acts on $M_\lambda$ and $\mathfrak{g}^{(l)}$ acts on $M_\lambda \otimes V$; thus, it suffices for us to consider the center's action on $M_\lambda \otimes V$. As long as $\lambda$ has trivial stabilizer in $W$, $M_\lambda \otimes V$ has a nice decomposition into Verma modules:
\begin{equation*}
  M_\lambda \otimes V =
  \begin{cases}
    \bigoplus_{i = 1}^n M_{\lambda + e_i} & \mathfrak{g}_n = \mathfrak{gl}_n \\
    \bigoplus_{i = 1}^n (M_{\lambda + e_i} \oplus M_{\lambda - e_i}) \oplus M_\lambda & \mathfrak{g}_n = \mathfrak{so}_{2k + 1} \\
    \bigoplus_{i = 1}^n (M_{\lambda + e_i} \oplus M_{\lambda - e_i}) & \mathfrak{g}_n = \mathfrak{sp}_{2k}.
  \end{cases}
\end{equation*}
Here $\lambda + e_i$ is the weight $(\lambda_1, \dots, \lambda_i + 1, \dots, \lambda_n)$. Then we can ``extract'' $\lambda_i$: define an operator $\Omega := \frac{1}{2}(\Delta(C_2) - C_2 \otimes 1 - 1)$. Recall that $C_k$ acts on $M_\lambda$ as $\left(\sum \lambda_i^k\right) - \rho_i^k$, though in the $\mathfrak{o}, \mathfrak{sp}$ cases we are only considering when $k$ is even.
\begin{eqnarray*}
  (2\Omega + 1)|_{M_{\lambda + e_i}} &= 2\lambda_i + 1 \\
  (2\Omega + 1)|_{M_\lambda} &= -\frac{1}{2}
\end{eqnarray*}
and
\begin{equation*}
  (\Delta(C_k) - C_k \otimes 1)_{U_\psi \otimes V} = (\Omega_{U_\psi \otimes V} + 1)^k - (\Omega_{U_\psi \otimes V})^k
\end{equation*}
for all $k$ if $\mathfrak{g} = \mathfrak{gl}$, and only for even $k$ if $\mathfrak{g} = \mathfrak{o}, \mathfrak{sp}$.
\begin{Definition}{}{}
  Define $P_k$ to be the polynomial $P_k(x) = (x + 1)^k - x^k$ and $f_x(u) = \sum \frac{P_k(x)}{k!} u^k$. In the $\mathfrak{g} = \mathfrak{o}, \mathfrak{sp}$ case, this is only defined for even $k$. We can write down a closed form for $f$:
  \begin{equation*}
    f_x(u) =
\begin{cases}
  e^{(x + 1)u} - e^{xu} & \mathfrak{g} = \mathfrak{gl} \\
  \cosh((x + 1)u) - \cosh(xu) & \mathfrak{g} = \mathfrak{o}, \mathfrak{sp}.
  \end{cases}
  \end{equation*}
\end{Definition}
Then we can write the above more compactly as
\begin{equation*}
  (\Delta(C_k) - C_k \otimes 1)_{U_\psi \otimes V} = P_k(\Omega_{U_\psi \otimes V}).
\end{equation*}
Now we have still assumed that $\psi$ corresponds to a generic weight (trivial stabilizers in $W$), but one can show that actually this is true for all central characters $\psi$. So we have
\begin{Lemma}{}{}
  For any central character $\psi$ of $U(\mathfrak{g}_n)$,
\begin{equation*}
  (\Delta(C_k) - C_k \otimes 1)_{U_\psi \otimes V} = P_k(\Omega_{U_\psi \otimes V}).
\end{equation*}  
\end{Lemma}
Now because central characters $\psi$ of $U(\mathfrak{g}_t)$ are ultraproducts of central characters $\psi^{(n)}$ of $U(\mathfrak{g}_n)$, we can take an ultraproduct of this equation over the $\psi^{(n)}$ to obtain
\begin{Corollary}{}{}
  For any central character $\psi$ of $U(\mathfrak{g}_k)$,
\begin{equation*}
  (\Delta(C_k) - C_k \otimes 1)_{U_\psi \otimes V} = P_k(\Omega_{U_\psi \otimes V}).
\end{equation*}
\end{Corollary}
Now we can move ``back'' to arbitrary finitely generated $(\mathfrak{g}_t, \mathfrak{g}_t)$-bimodules in $\Ind(\mathcal{C})$ by extending our result to $U(\mathfrak{g}_t) \otimes U(\mathfrak{g}_t)^{op} \otimes X$ for $X \in \Ob(\mathcal{C})$.
\begin{Lemma}{}{}
For any object $X \in \Ob(\mathcal{C})$, let $M = U(\mathfrak{g}_t) \otimes U(\mathfrak{g}_t)^{op} \otimes X$. Then
\begin{equation*}
  (\Delta(C_k) - C_k \otimes 1)_{M \otimes V} = P_k(\Omega_{M \otimes V}).
\end{equation*}
\end{Lemma}
\begin{proof}
  First, extend our result to modules of the form $M' = U_\chi \otimes U(\mathfrak{g}_t)^{op} \otimes X$ (where $U(\mathfrak{g}_t^{(l)})$ acts only on $U_\chi$). Passing again to the classical case, $X = \prod_{\mathcal{F}} X_i$, and for each $X_i$ we have
  \begin{equation*}
    U_\chi \otimes U(\mathfrak{g}_n^{op}) \otimes X_n \otimes V \hookrightarrow \Hom(M_\lambda, M_\lambda \otimes V) \otimes U(\mathfrak{g}_n) \otimes X_n.
  \end{equation*}
  Now again the left action $\mathfrak{g}_n^{(l)}$ only acts on $M_\lambda \oplus V$, so repeating the proof from above, we have that
\begin{equation*}
  (\Delta(C_k) - C_k \otimes 1)_{M' \otimes V} = P_k(\Omega_{M' \otimes V}).
\end{equation*}
  To further get this for $M$, notice that $U(\mathfrak{g}_n)$ embeds into $\bigoplus_{\chi^{(n)}} U_{\chi^{(n)}}$, where the sum is taken over all central characters of $U(\mathfrak{g}_n)$, so the result follows.
\end{proof}
\begin{Corollary}{}{}
  For any finitely generated $(\mathfrak{g}_t, \mathfrak{g}_t)$-bimodule $M$,
  \begin{equation*}
  (\Delta(C_k) - C_k \otimes 1)_{M \otimes V} = P_k(\Omega_{M \otimes V}).
\end{equation*}
\end{Corollary}

\begin{Note}{}{}
In the $\mathfrak{g} = \mathfrak{gl}$ case, because $V$ is not self-dual, we also have to look at what happens when we tensor with $V^*$; it is not hard to see that
  \begin{equation*}
  (\Delta(C_k) - C_k \otimes 1)_{M \otimes V^*} = \overline{P}_k(\Omega_{M \otimes V^*})
\end{equation*}
where $\overline{P}_k(c) = (c - 1)^k - c^k$.
\end{Note}

Now we look at how $\Delta(C_k) - C_k \otimes 1$ acts on $U_\psi \otimes [r, s]$ or $U_\psi \otimes [r]$. Let
\begin{equation*}
  \tau_j^m = (\Delta^{j - 1} \otimes \id) \otimes 1 \otimes \cdots \otimes 1: U(\mathfrak{g}_t) \otimes U(\mathfrak{g}_t) \to U(\mathfrak{g}_t)^{m + 1}
\end{equation*}
where there are $m - j$ $1$s in the tensor product. It will be evident from context what $m$ is, so we omit it. For an element $A \in Z(U(\mathfrak{g}_t) \otimes U(\mathfrak{g}_t))$, define $A_j = \tau_j$ as an operator on $X_0 \otimes \cdots \otimes X_m$ for $U(\mathfrak{g_t})$-modules $X_i$. Treating $[r, s] = V^{\otimes r} \otimes V^{* \otimes s}$ and $[r] = V^{\otimes r}$, we can use a telescoping argument to show that for $A \coloneq \Delta(C_k) - C_k \otimes 1$,
\begin{align}
  \label{eq:telescope}
  A_{U_\psi \otimes [r, s]} &= A_{r + s} + A_{r + s - 1} + \cdots + A_1 \nonumber \\
                            &= \overline{P}_k(\Omega_{r + s}) + \cdots + \overline{P}_k(\Omega_{r + 1}) + P_k(\Omega_r) + \cdots + P_k(\Omega_1); \\
  A_{U_\psi \otimes [r]} &= A_r + A_{r - 1} + \cdots + A_1 \nonumber \\
  &= P_k(\Omega_r) + \cdots + P_k(\Omega_1).
\end{align}

\begin{Theorem}{}{}
  If $\mathfrak{g} = \mathfrak{gl}$ and $(U_\psi \otimes [r, s])_\chi \ne 0$, then there exist $b_1, \dots, b_r, c_1, \dots, c_s \in \mathbb{C}$ so that
  \begin{equation*}
    \chi(u) - \psi(u) = f_{b_1}(u) + \cdots + f_{b_r}(u) - f_{c_1 - 1}(u) + \cdots + f_{c_s - 1}(u).
  \end{equation*}
  If $\mathfrak{g} = \mathfrak{o}, \mathfrak{sp}$ and $(U_\psi \otimes [r])_\chi \ne 0$, then there exist $b_1, \dots, b_r \in \mathbb{C}$ so that
  \begin{equation*}
    \chi(u) - \psi(u) = f_{b_1}(u) + \cdots + f_{b_r}(u).
  \end{equation*}
  In other words, if $HC_{\chi, \psi} \ne 0$, there exist complex numbers satisfying the above conditions.
\end{Theorem}
\begin{proof}
  For convenience of notation, let $N = N(\chi, \psi, [r, s])$ or $N(\chi, \psi, [r])$ depending on $\mathfrak{g}$, and likewise $m = r + s$ or $r$. We know $\Delta(C_k) - C_k \otimes 1$ acts on $N$ as $\chi_k - \psi_k$. By assumption, because $N \ne 0$, $\End(N) \ne 0$ also, and we have a homomorphism $\mathbb{C}[\Omega_1, \dots, \Omega_m] \to \End(N)$. The image of this map is a nonzero finitely generated commutative algebra $A$; abusing notation and letting $\Omega_i$ also be the image of $\Omega_i$ under the map, \cref{eq:telescope} implies
  \begin{equation*}
    P_k(\Omega_1) + \cdots + P_k(\Omega_r) + \overline{P}_k(\Omega_{r + 1}) + \cdots + \overline{P}_k(\Omega_{r + s}) = \chi_k - \psi_k
  \end{equation*}
  when $\mathfrak{g} = \mathfrak{gl}$ and
  \begin{equation*}
    P_k(\Omega_1) + \cdots + P_k(\Omega_r) = \chi_k - \psi_k.
  \end{equation*}
  By the Nullstellensatz, there exists a maximal ideal in $A$. If $\mathfrak{g} = \mathfrak{gl}$, that ideal corresponds to complex numbers $b_1, \dots, b_r$ and $c_1, \dots, c_s$ where
  \begin{equation*}
    P_k(b_1) + \cdots + P_k(b_r) + \overline{P}_k(c_1) + \cdots + \overline{P}_k(c_s) = \chi_k - \psi_k.
  \end{equation*}
  If $\mathfrak{g} = \mathfrak{o}, \mathfrak{sp}$, that ideal corresponds to $b_1, \dots, b_r \in \mathbb{C}$ with
  \begin{equation*}
    P_k(b_1) + \cdots + P_k(b_r) = \chi_k - \psi_k.
  \end{equation*}
  By the definition of exponential central characters in \cref{def:exp_cc}, this implies the result.
\end{proof}

Now we want to show the converse, which we will do by constructing an object $X$ with nonzero $N(\chi, \psi, X)$ if $\chi$ and $\psi$ satisfy the above condition.

First, we will find a sufficient condition for $N(\chi, \psi, X) \ne 0$ in the classical case, i.e. representations and characters of $\mathfrak{g}_n$. We will relate this condition to our expression for $\chi(u) - \psi(u)$ above. Then, we will show that characters in the $\chi, \psi \in \mathfrak{g}_t$ case can be written as ultraproducts of central characters $\prod_{\mathcal{F}} \chi^{(n)}, \prod_{\mathcal{F}} \psi^{(n)}$ of $U(\mathfrak{g}_n)$ so that $(\chi^{(n)}, \psi^{(n)})$ satisfy the condition we find. Finally, we will construct an $X \in \Ob(\mathcal{C}_t)$ so that $N(\chi, \psi, X) \ne 0$.

\begin{Lemma}{}{nonzero_cond}
  Suppose that $\chi, \psi$ are central characters of $U(\mathfrak{g}_n)$ corresponding to weights $\lambda, \mu$ respectively with $\lambda - \mu \in \Lambda$, the root lattice. Let $X$ be a finite-dimensional $\mathfrak{g}_n$-module with $X[\lambda - \mu] \ne 0$, i.e. the $\lambda - \mu$ weight space of $X$ is nonzero. Then $N(\chi, \psi, X) \ne 0$.
\end{Lemma}
\begin{proof}
  It suffices to show that $(M_\mu \otimes X)_\chi \ne 0$, since we have a surjective map $U_\psi \otimes X \twoheadrightarrow M_\mu \otimes X$ that is still surjective after tensoring with $U_\chi$. We have $\ch M_\mu \otimes X = \sum_{\nu} (\dim X[\nu]) \ch M_{\mu + \nu}$. Let $(M_\mu \otimes X)^\varphi$ denote the projection of $M_\mu \otimes X$ to the block $\mathcal{O}_\varphi$, i.e. modules on which $U(\mathfrak{g}_n)$ acts by the generalized central character $\varphi$. Since $\ch N_1 = \ch N_2 \implies \ch N_1^\varphi = \ch N_2^\varphi$ for $N_1, N_2 \in \mathcal{O}$, and $M_\lambda$ has central character $\chi$,
  \begin{equation*}
    \ch (M_\mu \otimes X)^\chi = \dim X[\lambda - \mu] \cdot \ch M_\lambda + \text{ characters of other representations}.
  \end{equation*}
  Now $\dim X[\lambda - \mu] \ne 0$ by assumption, so $\ch(M_\mu \otimes X)^\chi \ne 0$ also. Then $(M_\mu \otimes X)^\chi$ surjects onto a simple module $L \in \mathcal{O}_\chi$, and tensoring with $U_\chi$, we have a surjective map $(M_\mu \otimes X)_\chi \twoheadrightarrow L$. This implies $(M_\mu \otimes X)_\chi \ne 0$, so we are done.
\end{proof}

\begin{Lemma}{}{lambda_mu_diff}
  Let
  \begin{equation*}
    \alpha =
    \begin{cases}
      e_1 + \cdots + e_r - e_{n - s + 1} - \cdots - e_n & \mathfrak{g} = \mathfrak{gl} \\
      e_1 + \cdots + e_r & \mathfrak{g} = \mathfrak{o}, \mathfrak{sp}.
    \end{cases}
  \end{equation*}
  Let $n > r + s$ if $\mathfrak{g} = \mathfrak{gl}$ and $n > r$ otherwise; pick $\lambda, \mu$ corresponding to central characters $\chi, \psi$ of $U(\mathfrak{g}_n)$ so that $\lambda - \mu = \alpha$.
  Then
  \begin{equation*}
    \chi(u) - \psi(u) =
    \begin{cases}
      f_{\mu_1}(u) + \cdots + f_{\mu_r}(u) - f_{\mu_{n - s + 1} - 1}(u) - \cdots - f_{\mu_n - 1}(u) & \mathfrak{g} = \mathfrak{gl} \\
      f_{\mu_1}(u) + \cdots + f_{\mu_r}(u) & \mathfrak{g} = \mathfrak{o}, \mathfrak{sp}.
    \end{cases}
  \end{equation*}
\end{Lemma}
\begin{proof}
  We have
  \begin{equation}
    \label{eq:ck_shift}
    C_k|_{M_{\mu + e_l}} - C_k|_{M_\mu} = P_k(\mu_l)
  \end{equation}
  (only considering even $k$ if $\mathfrak{g} = \mathfrak{o}, \mathfrak{sp}$), and if $\mathfrak{g} = \mathfrak{gl}$, also that
  \begin{equation*}
    C_k|_{M_{\mu - e_l}} - C_k|_{M_\mu} = \overline{P}_k(\mu_l).
  \end{equation*}
  Now we travel from $\mu$ to $\lambda$, taking steps $e_i$ for the various $i$. Each step will add a $P_k$ or an $\overline{P}_k$. Inductively define $\lambda^{[0]} := \lambda$ and $\lambda^{[i]} := \lambda^{[i - 1]} - e_i$ for $1 \le i \le r$. If $\mathfrak{g} = \mathfrak{o}, \mathfrak{sp}$, then $\lambda^{[r]} = \mu$. Let $\chi^{[i]}$ be the central character corresponding to $\lambda^{[i]}$, with $\chi^{[r]} = \psi$; thus, \cref{eq:ck_shift} implies $\chi_k^{[i]} - \chi_k^{[i + 1]} = P_k(\mu_{i + 1})$. If we sum up these equations for $0 \le i < r$, we obtain that
  \begin{equation*}
    \chi_k - \psi_k = P_k(\mu_1) + \cdots + P_k(\mu_r).
  \end{equation*}
  If $\mathfrak{g} = \mathfrak{gl}$, we also have to move from $\lambda^{[r]}$ to $\mu$, so we also inductively define $\mu^{[0]} := \mu$ and $\mu^{[i]} := \mu^{[i - 1]} - e_{n - i + 1}$ for $1 \le i \le s$. Let $\psi^{[i]}$ be the central character corresponding to $\mu^{[i]}$. So $\mu^{[s]} = \lambda^{[r]}$, $\psi^{[s]} = \chi^{[r]}$. We have $\psi_k^{[i + 1]} - \psi_k^{[i]} = \overline{P}_k(\mu_{n - i})$, so summing up these differences for the $\chi$s and $\psi$s gives
  \begin{equation*}
    \chi_k - \psi_k = P_k(\mu_1) + \cdots + P_k(\mu_r) + \overline{P}_k(\mu_{n - s + 1}) + \cdots + \overline{P}_k(\mu_n).
  \end{equation*}
\end{proof}

Now we show that we can pick ultraproduct presentations for $\chi, \psi$ so that $\chi^{(n)}, \psi^{(n)}$ correspond to $\lambda^{(n)}, \mu^{(n)}$ with $\lambda^{(n)} - \mu^{(n)}$ as in \cref{lem:lambda_mu_diff}. First, we prove a more general statement:
\begin{Lemma}{}{b_c_ultraproduct}
  First suppose $\mathfrak{g} = \mathfrak{o}, \mathfrak{sp}$. Let $b_i = \prod_{\mathcal{F}} b_i^{(n)}$ for $0 \le i \le r$. Then any central character $\psi$ can be written as $\prod_{\mathcal{F}} \psi^{(n)}$ with
  \begin{equation*}
    \psi^{(n)} \begin{cases}
      = 0 & \text{if } n \le r \\
      \text{corresponds to } \mu^{(n)}, \mu^{(n)}_i = b_i^{(n)}, 1 \le i \le r & \text{if } n > r.
    \end{cases}
  \end{equation*}

  so that $\psi^{(n)} = 0$ when $n \le r$ and corresponds to a weight $\mu^{(n)}$ with $\mu^{(n)}_i = b_i^{(n)}$ for $1 \le i \le r$ when $n > r$.
  If $\mathfrak{g} = \mathfrak{gl}$, we also want to consider $c_j = \prod_{\mathcal{F}} c_j^{(n)}$ for $0 \le j \le s$, and instead we have
  \begin{equation*}
    \psi^{(n)} \begin{cases}
      = 0 & \text{if } n \le r + s\\
      \text{corresponds to } \mu^{(n)}, \mu^{(n)}_i = b_i^{(n)}, 1 \le i \le r , \mu^{(n)}_{n - j + 1} = c_j^{(n)}, 1 \le j \le s & \text{if } n > r + s.
    \end{cases}
  \end{equation*}
\end{Lemma}
\begin{proof}
  Begin with $\psi_k = \prod_{\mathcal{F}} \varphi^{(n)}$, where there are no conditions on the $\varphi^{{n}}$. For each $k$, we show that we only need to change $\varphi^{(n)}$ when $n < k + r$, or $k + r + s$ if $\mathfrak{g} = \mathfrak{gl}$, so that the resulting central characters satisfy the condition we want; then, $\psi_k$ is unchanged because we fixed all but finitely many terms of the ultraproduct. So fix a $k$.

  First suppose $\mathfrak{g} = \mathfrak{o}, \mathfrak{sp}$. If $n \le r$, we set $\psi_k^{(n)} = 0$; if $r < n \le k + r$, pick a $\psi^{(n)}$ with a corresponding $\mu^{(n)}$ so that $\mu^{(n)}_i = b_i^{(n)}$ for $1 \le i \le r$.
  Fix an $n > r$, let $m := n - r$, and for $1 \le k \le n - r$, we claim that we can keep $\psi_k^{(n)} = \varphi_k^{(n)}$ even if we set the first $r$ coordinates of $\mu^{(n)}$ to the $b_i^{(n)}$. We want to set the other $m$ coordinates of $\mu^{(n)}$, which we denote by $x_1, \dots, x_m$, so that

  \begin{equation*}
    (b_1^{(n)})^k + \cdots + (b_r^{(n)})^k + x_1^k + \cdots + x_m^k = \varphi_k^{(n)}
  \end{equation*}
  when $1 \le k \le m$. So we have requirements for the first $m$ power sums of the $x_i$, and the fundamental theorem of symmetric functions tells us that we have a solution $x_i = a_i^{(n)}$. Set $\psi_k^{(n)}$ to be the central character corresponding to $(b_1^{(n)}, \dots, b_r^{(n)}, a_1^{(n)}, \dots, a_m^{(n)})$. By construction $\psi_k^{(n)} = \varphi_k^{(n)}$. Thus, if $n > k + r$, we can just let $\psi_k^{(n)} = \varphi_k^{(n)}$.
  So, $\prod_{\mathcal{F}} \psi_k^{(n)} = \psi_k$ as well.

  If $\mathfrak{g} = \mathfrak{gl}$, if $n \le r + s$, set $\psi_k^{(n)} = 0$, and if $r + s< n \le k + r + s$, pick $\psi^{(n)}$ with corresponding $\mu^{(n)}$ so $\mu^{(n)}_i = b_i^{(n)}$ and $\mu^{(n)}_{n - j + 1} = c_j^{(n)}$ for $1 \le i \le r$, $1 \le j \le s$. We can still construct the $a_i^{(n)}$ for $x_1, \dots, x_{m}$, $m := n - r - s$, by asking that the $x_i$ satisfy
  \begin{equation*}
    (b_1^{(n)})^k + \cdots + (b_r^{(n)})^k + x_1^k + \cdots + x_m^k + (c_1^{(n)})^k + \cdots + (c_s^{(n)})^k = \varphi_k^{(n)}.
  \end{equation*}
  So in this case, we have $\psi_k^{(n)} = \varphi_k^{(n)}$ when $n \ge k + r + s$, the conditions we want are satisfied for all $n$, and $\prod_{\mathcal{F}} \psi_k^{(n)} = \psi_k$.
\end{proof}
\begin{Corollary}{}{set_ultraproduct_alpha}
  Let $m := r + s$ if $\mathfrak{g} = \mathfrak{gl}$ and $m := r$ if $\mathfrak{g} = \mathfrak{o}, \mathfrak{sp}$.
  Suppose that
  \begin{equation*}
    \chi(u) - \psi(u) =
    \begin{cases}
      f_{b_1}(u) + \cdots + f_{b_r}(u) - f_{c_1 - 1}(u) - \cdots - f_{c_s - 1}(u) & \mathfrak{g} = \mathfrak{gl} \\
      f_{b_1}(u) + \cdots + f_{b_r}(u) & \mathfrak{g} = \mathfrak{o}, \mathfrak{sp}.
  \end{cases}
\end{equation*}
Then, there exist ultraproduct presentations $\chi = \prod_{\mathcal{F}} \chi^{(n)}$, $\psi = \prod_{\mathcal{F}} \psi^{(n)}$ so that for any $n > m$, there exist $\lambda^{(n)}, \mu^{(n)}$ weights with
\begin{equation*}
  \lambda^{(n)} - \mu^{(n)} =
  \begin{cases}
    e_1 + \cdots + e_r - e_{n - s + 1} - \cdots - e_n & \mathfrak{g} = \mathfrak{gl} \\
    e_1 + \cdots + e_r & \mathfrak{g} = \mathfrak{o}, \mathfrak{sp}
  \end{cases}
\end{equation*}
so that $\chi^{(n)}, \psi^{(n)}$ correspond with $\lambda^{(n)}, \mu^{(n)}$.
\end{Corollary}
\begin{proof}
  We pursue a similar strategy as in the proof of \cref{lem:b_c_ultraproduct}: take advantage of how the ultraproduct only sees what happens on almost all terms. First, pick $\psi^{(n)}$ and corresponding weights $\mu^{(n)}$ satisfying the conditions in \cref{lem:b_c_ultraproduct}. Now for $n > m$, set $\lambda^{(n)}$ so that $\lambda^{(n)} - \mu^{(n)}$ is the difference we want, and let $\chi^{(n)}$ be the central character corresponding to the $\lambda^{(n)}$. By \cref{lem:lambda_mu_diff},
  \begin{equation*}
    \chi^{(n)}(u) - \psi^{(n)}(u) =
    \begin{cases}
      f_{b_1^{(n)}}(u) + \cdots + f_{b_r^{(n)}}(u) - f_{c_1^{(n)} - 1}(u) - \cdots - f_{c_s^{(n)}- 1}(u) & \mathfrak{g} = \mathfrak{gl} \\
      f_{b_1^{(n)}}(u) + \cdots + f_{b_r^{(n)}}(u) & \mathfrak{g} = \mathfrak{o}, \mathfrak{sp}.
    \end{cases}
  \end{equation*}
  Set $\chi^{(n)} = 0$ for $n \le r + s$. Then
  \begin{equation*}
    \prod_{\mathcal{F}} \chi^{(n)}(u) - \psi(u) = \chi(u)- \psi(u),
  \end{equation*}
  so $\chi(u) = \prod_{\mathcal{F}} \chi^{(n)}(u)$ and we have found our desired ultraproduct presentations.
\end{proof}
\begin{Theorem}{}{}
  Let $\chi, \psi$ be central characters of $\mathfrak{g}_t$. Suppose that there exist $b_1, \dots, b_r, c_1, \dots, c_s \in \mathbb{C}$ so that
  \begin{equation*}
    \chi(u) - \psi(u) = 
    \begin{cases}
      f_{b_1}(u) + \cdots + f_{b_r}(u) - f_{c_1 - 1}(u) - \cdots - f_{c_s - 1}(u) & \mathfrak{g} = \mathfrak{gl} \\
      f_{b_1}(u) + \cdots + f_{b_r}(u) & \mathfrak{g} = \mathfrak{o}, \mathfrak{sp}.
    \end{cases}
  \end{equation*}
  Then, $HC_{\chi, \psi} \ne 0$.
\end{Theorem}
\begin{proof}
Let
  \begin{align*}
    m &=
        \begin{cases}
          r + s & \mathfrak{g} = \mathfrak{gl} \\
          r & \mathfrak{g} = \mathfrak{o}, \mathfrak{sp}
        \end{cases} \\
    X &=
        \begin{cases}
          S^r V \otimes S^s V^* & \mathfrak{g} = \mathfrak{gl} \\
          S^r V & \mathfrak{g} = \mathfrak{o}, \mathfrak{sp}
        \end{cases} \\
    X^{(n)} &=
              \begin{cases}
          S^r V^{(n)} \otimes S^s V^{*(n)} & \mathfrak{g} = \mathfrak{gl} \\
          S^r V^{(n)} & \mathfrak{g} = \mathfrak{o}, \mathfrak{sp}.
        \end{cases}
  \end{align*}
  We show that $N(\chi, \psi, X) \ne 0$. Since $X = \prod_{\mathcal{F}} X^{(n)}$, and for $n > m$ $\dim X^{(n)}[\alpha] = 1$, we can pick $\chi^{(n)}, \psi^{(n)}$ satisfying the conditions of \cref{cor:set_ultraproduct_alpha} so that \cref{lem:nonzero_cond} implies that $N(\chi^{(n)}, \psi^{(n)}, X^{(n)}) \ne 0$. Then, $N(\chi, \psi, X) \ne 0$ as well and we are done.
\end{proof}

Finally, by using facts about $\exp$, $\cosh$, and $\sinh$, we can arrive at the statement of \cref{theo:central_characters}.

\section{Harish-Chandra bimodules of finite $K$-type}

Analyzing objects in $HC_{\chi, \psi}$ also allows us to construct families of Harish-Chandra bimodules with finite $K$-type. We cannot just take an arbitrary ultraproduct of the classical Harish-Chandra bimodules, although they are of finite $K$-type, for the same reason as in \cref{xmpl:ultra_vec}, but if we ensure that the modules in our ultraproduct have universally bounded multiplicities of irreducibles, then we can indeed obtain a Harish-Chandra bimodule of $\mathfrak{g}_t$ that has finite $K$-type. The proof that ultraproduct of these modules actually lies in the ind-completion is rather technical, and we refer the reader to \cite{utiralova2022harishchandra}.

It turns out that there are already closed-form combinatorial expressions for the multiplicities of irreducibles in tensor products of irreducible modules for $GL$, $O$, and $Sp$, due to Brylinski (\cite{brylinski}) and King (\cite{king}). Let $\lambda$ and $\mu$ be two partitions, corresponding to dominant integral weights, and let $V_\lambda, V_\mu$ be their corresponding irreducible representations for $\mathfrak{g}_n$ (if $\mathfrak{g} = \mathfrak{gl}$, we will abuse notation and also let $\lambda$ denote the bipartition $(\lambda, \emptyset)$). Then we have
\begin{Theorem}{}{brylinski}
  Let $\boldsymbol{\nu} = (\nu, \overline{\nu})$ be a bipartition. For $n$ large enough,
  \begin{equation*}
    \dim \Hom_{\mathfrak{gl}_n} (V_{\boldsymbol{\nu}}, V_\lambda \otimes V_\mu^*) = (s_{\lambda / \nu}, s_{\mu / \overline{\nu}}) = \sum_{\eta \subset \lambda, \mu} c^\lambda_{\nu, \eta} c^\mu_{\overline{\nu}, \eta}.
  \end{equation*}
  Here $s_{\alpha / \beta}$ is the skew Schur polynomial associated with the skew diagram $\alpha / \beta$, and $c^\alpha_{\beta, \gamma}$ is the Littlewood-Richardson coefficient.
\end{Theorem}
Brylinski proved this when $|\lambda| = |\mu|$ and Utiralova extended this to all $\lambda, \mu$ in her forthcoming paper.

\begin{Theorem}{King}{king}
  Let $\lambda, \mu, \nu$ be arbitrary partitions. For $n$ large enough,
  \begin{equation*}
    \dim \Hom_{\mathfrak{g}_n}(V_\nu, V_\lambda \otimes V_\mu) = \sum_{\zeta, \sigma, \tau} c^\lambda_{\zeta, \sigma} c^\mu_{\zeta, \tau} c^\nu_{\sigma, \tau}
  \end{equation*}
  for $\mathfrak{g} = \mathfrak{o}, \mathfrak{sp}$.
\end{Theorem}

Using these, we will associate to some pairs of partitions $\lambda, \mu$ a space denoted $\Hom(\mu, \lambda)$ that is a Harish-Chandra bimodule of finite $K$-type.

\begin{Definition}{}{}
  We can think of every partition $\lambda$ as determined by a triple of partitions $(\alpha, \beta, \gamma)$. Cut $\lambda$ into four parts by cutting under the $k^{th}$ row and to the right of the $l^{th}$ column. Let $\alpha$ be the top right part, $\beta'$ be the bottom left part, and $\gamma$ be the bottom right part. Then we have the following constraints on $\alpha, \beta, \gamma$, which also allow us to recover $\lambda$:
  \begin{enumerate}
  \item $k := l(\alpha)$ and $l := l(\beta) \le d(\lambda)$, the length of the main diagonal of $\lambda$;
  \item $l(\gamma) = \lambda'_{l + 1} - k$;
  \item $\alpha_i = \lambda_i - l$ when $i \le k$ and $\beta_j = \lambda'_j - k$ when $j \le l$;
  \item $\gamma_i = \lambda_{k + i} - l$ when $i \le l(\gamma)$;
  \item $\gamma_1 \le \alpha_k, \gamma'_1 \le \beta_l$.
  \end{enumerate}
  We will write $\lambda = [\alpha, \beta, \gamma]$.
\end{Definition}
\begin{Note}{}{}
  This triple representation is not unique, nor does every triple of partitions result in a well-defined $\lambda$.
\end{Note}
This notation will be convenient when we think about multiplicities of irreducibles. Now we define $\Hom(\mu, \lambda)$. First, we must describe a filtration on $\mathfrak{g}_t$.
\begin{Definition}{}{s_filtration}
  Let
  \begin{equation*}
    S =
    \begin{cases}
      \mathcal{P} \times \mathcal{P} & \mathfrak{g} = \mathfrak{gl} \\
      \mathcal{P} & \mathfrak{g} = \mathfrak{o}, \mathfrak{sp}.
    \end{cases}
  \end{equation*}
  Let $F$ be a filtration by finite subsets on $S$ (e.g. partitions of length at most $n$ if $S = \mathcal{P}$), which will induce a filtration we will also denote by $F$ on $\Rep(\mathfrak{g}_n)$ and $\Rep(\mathfrak{g}_t)$:
  \begin{equation*}
    F^k X = \bigoplus_{\lambda \in F^k S} V_\lambda \otimes \Hom(V_\lambda, X).
  \end{equation*}
\end{Definition}
\begin{Definition}{}{}
  Fix integers $k, l$; partitions $\gamma$, $\delta$; and integer sequences $a \in \mathbb{Z}^k, b \in \mathbb{Z}^l$. Let $|a| = \sum a_i$ and $|b| = \sum b_i$. Let $\lambda^{(n)} = [\alpha^{(n)}, \beta^{(n)}, \gamma]$ be a sequence of partitions with $l(\lambda^{(n)}) \ll n$, $l(\alpha^{(n)}) = k$, $l(\beta^{(n)}) = l$, and assume that
  \begin{equation*}
    \lim_{n \to \infty}(\theta_i^{(n)} - \theta_{i + 1}^{(n)}) = \infty, i \le l(\theta), \theta \in \{\alpha, \beta\}.
  \end{equation*}
  Then, for large enough $n$, the partition
  \begin{equation*}
    \mu^{(n)} := [\alpha^{(n)} + a, \beta^{(n)} + b, \delta]
  \end{equation*}
  is well-defined.

  Now we can define $\Hom(\lambda, \mu)$: let
  \begin{align*}
    \lambda &\coloneq (\alpha, \beta, \gamma) &\mu &\coloneq (\alpha + a, \beta + b, \delta) \\
    \alpha &\coloneq \prod_{\mathcal{F}} \alpha^{(n)} &\beta &\coloneq \prod_{\mathcal{F}} \beta^{(n)}
  \end{align*}
  and
  \begin{equation*}
    \Hom_{\mathfrak{g}_t}(\mu, \lambda) \coloneq \prod_{\mathcal{F}}^r \Hom_{\mathfrak{g}_n}(V_{\mu^{(n)}}, V_{\lambda^{(n)}}) = \prod_{\mathcal{F}}^r V_{\lambda^{(n)}} \otimes V_{\mu^{(n)}}^*.
  \end{equation*}
\end{Definition}
\begin{Lemma}{}{stabilized_nu}
For any (bi)partition $\boldsymbol{\nu}$, the multiplicity of $V_{\boldsymbol{\nu}}$ in $\Hom_{\mathcal{C}_n}(V_{\mu^{(n)}}, V_{\lambda^{(n)}})$ is constant for large enough $n$.
\end{Lemma}
\begin{proof}
  We have different formulas for the multiplicity of $V_{\boldsymbol{\nu}}$ depending on which Lie algebra we are considering, though it will turn out that the casework is similar. First fix a $\boldsymbol{\nu}$.

  \textit{General linear case.} Suppose $\mathfrak{g} = \mathfrak{gl}$; then $\boldsymbol{\nu} = (\nu, \overline{\nu})$. By \cref{theo:brylinski}, we have that the multiplicity of $V_{\boldsymbol{\nu}}$ is
  \begin{equation*}
    \sum_{\eta \subset \lambda^{(n)}, \mu^{(n)}} c_{\nu, \eta}^{\lambda^{(n)}} c_{\overline{\nu}, \eta}^{\mu^{(n)}}.
  \end{equation*}
  We claim that for $n$ large enough, this sum depends only on $k, l, a, b, \gamma, \delta$, the constants we chose in our construction of $\lambda^{(n)}$ and $\mu^{(n)}$. In other words, it is constant for $n$ large enough. The only $\eta$ that contribute to the sum are those with $|\lambda^{(n)}| - |\eta| = |\nu|$ and $|\mu^{(n)}| - |\eta| = |\overline{\nu}|$, and we want to say something about $\lambda^{(n)} / \eta$ and $\mu^{(n)} / \eta$. For $\eta = [\sigma, \tau, \epsilon]$; we have
  \begin{align*}
    \sigma &\subset \alpha^{(n)}, \alpha^{(n)} + a \\
    \tau &\subset \beta^{(n)}, \beta^{(n)} + b \\
    \epsilon &\subset \gamma, \delta.
  \end{align*}
  We can choose $n$ large enough so that
  \begin{equation*}
    \theta_i^{(n)} - \theta_{i + 1}^{(n)} > \max(|\nu|, |\overline{\nu}|), i \le l(\theta), \theta \in \{\alpha, \alpha + a, \beta, \beta + b\}.
  \end{equation*}
  Then, we have that
  \begin{align*}
    \lambda^{(n)} / \eta &= \alpha^{(n)} / \sigma \sqcup (\beta^{(n)})' / \tau' \sqcup \gamma / \epsilon \\
    \mu^{(n)} / \eta &= \alpha^{(n)} + a / \sigma \sqcup (\beta^{(n)} + b)' / \tau' \sqcup \delta / \epsilon
  \end{align*}
  where we are taking the disjoint union of three skew diagrams. $\alpha^{(n)} / \sigma$ is the disjoint union of rows of lengths $c_1, \dots, c_k \in \mathbb{Z}_{\ge 0}$; for ease of notation, let $c \coloneq (c_1, \dots, c_k)$. Likewise, $\beta'^{(n)} / \tau'$ is the disjoint union of columns of lengths $d_1, \dots, d_l \in \mathbb{Z}_{\ge 0}$ and we let $d \coloneq (d_1, \dots, d_l)$. $\eta$ is uniquely defined by $c, d, \epsilon$. Likewise, $\alpha^{(n)} + a / \sigma$ is a disjoint union of rows of length $c + a$, and $(\beta^{(n)} + b)' / \tau'$ is a disjoint union of rows of length $d + b$. Therefore, we have constraints on $c$ and $d$ depending only on $a, b$ respectively:
  \begin{align*}
    c &\in (-a + \mathbb{Z}^k_{\ge 0}) \cap \mathbb{Z}^k_{\ge 0} \\
    d &\in (-b + \mathbb{Z}^l_{\ge 0}) \cap \mathbb{Z}^l_{\ge 0}.
  \end{align*}
  This motivates us to construct skew diagrams with the same disjoint parts, in the same order, as $\lambda^{(n)} / \eta$ and $\mu^{(n)} / \eta$, but that depend only on $c$, $d$, $\epsilon$, $\gamma$, and $\delta$. Define
  \begin{align*}
    \widetilde{\alpha}_i &\coloneq \gamma_1 + \sum_{j = i}^k c_j \\
    \widetilde{\beta}_i &\coloneq \gamma'_1 + \sum_{j = i}^l d_j \\
    \widetilde{\lambda}(c, d, \gamma) &\coloneq [\widetilde{\alpha}, \widetilde{\beta}, \gamma] \\
    \widetilde{\eta}(c, d, \epsilon) &\coloneq [\widetilde{\alpha} - c, \widetilde{\beta} - d, \epsilon].
  \end{align*}
  By construction, $\widetilde{\lambda} / \widetilde{\eta}$ have the same disjoint parts in the same order as $\lambda^{(n)} / \eta$, so they have the same number of Littlewood-Richardson tableau; thus $c_{\nu, \eta}^{\lambda^{(n)}}$ is a function of $\nu, c, d, \gamma, \epsilon$, which we will denote as $c_\nu(c, d, \gamma, \epsilon)$. So
  \begin{equation*}
    \dim \Hom_{\mathfrak{g}_n}(V_{\boldsymbol{\nu}}, V_{\lambda^{(n)}}\otimes V^*_{\mu^{(n)}}) = \sum_{c, d, \epsilon} c_\nu(c, d, \gamma, \epsilon) c_{\overline{\nu}}(c + a, d + b, \delta, \epsilon)
  \end{equation*}
  where $\epsilon \subset \gamma, \delta$ and
  \begin{align*}
    c &\in (-a + \mathbb{Z}^k_{\ge 0}) \cap \mathbb{Z}^k_{\ge 0} \\
    d &\in (-b + \mathbb{Z}^l_{\ge 0}) \cap \mathbb{Z}^l_{\ge 0},
  \end{align*}
  \begin{align*}
    &|c| + |d| + |\gamma| - |\epsilon| = |\nu| \\
    &|c| + |d| + |a| + |b| + |\delta| - |\epsilon| = |\overline{\nu}|.
  \end{align*}
  Thus, as promised, the sum depends only on $k$, $l$, $a$, $b$, $\gamma$, and $\delta$ for $n$ large enough.

  \textit{Orthogonal and symplectic cases.} Because $\boldsymbol{\nu}$ is a partition, we just denote it by $\nu$. As stated in \cref{theo:king}, the multiplicity of $V_\nu$ in $V_{\lambda^{(n)}} \otimes V_{\mu^{(n)}}^*$ is
    \begin{equation*}
      \sum_{\eta, \omega, \xi} c^{\lambda^{(n)}}_{\eta, \omega} c^{\mu^{(n)}}_{\eta, \xi} c^\nu_{\omega, \xi}.
    \end{equation*}
    Again, we claim that this sum depends only on $k$, $l$, $a$, $b$, $\gamma$, and $\delta$. Though the sum is over all partitions $\eta, \omega, \xi$, we see that for a term to be nonzero, $\eta \subset \lambda^{(n)}, \mu^{(n)}$. There are only finitely many $(\omega, \xi)$ pairs so that $c^\nu_{\omega, \xi} \ne 0$ since $\nu$ is fixed and $|\omega| + |\xi| = |\nu|$. Thus, we will focus on $c^{\lambda^{(n)}}_{\eta, \omega}$ and $c^{\mu^{(n)}}_{\eta, \xi}$. As above, we see that the possible $\eta = [\sigma, \tau, \epsilon]$ are constrained by $c$, $d$, where
  \begin{align*}
    c &\in (-a + \mathbb{Z}^k_{\ge 0}) \cap \mathbb{Z}^k_{\ge 0} \\
    d &\in (-b + \mathbb{Z}^l_{\ge 0}) \cap \mathbb{Z}^l_{\ge 0},
  \end{align*}
  and we construct the same $\widetilde{\lambda}$ and $\widetilde{\eta}$. Using the same notation as above,
\begin{equation*}
    \dim \Hom_{\mathfrak{g}_n}(V_{\nu}, V_{\lambda^{(n)}}\otimes V^*_{\mu^{(n)}}) \le \sum_{c, d, \epsilon, |\omega| + |\xi| = |\nu|} c_\omega(c, d, \gamma, \epsilon) c_{\xi}(c + a, d + b, \delta, \epsilon) c^\nu_{\omega, \xi}.
  \end{equation*}
  Here $\epsilon \subset \gamma, \delta$,
  \begin{align*}
    c &\in (-a + \mathbb{Z}^k_{\ge 0}) \cap \mathbb{Z}^k_{\ge 0} \\
    d &\in (-b + \mathbb{Z}^l_{\ge 0}) \cap \mathbb{Z}^l_{\ge 0},
  \end{align*}
  and
  \begin{align*}
    &|c| + |d| + |\gamma| - |\epsilon| = |\omega|, \\
    &|c| + |d| + |a| + |b| + |\delta| - |\epsilon| = |\xi| \\
    \implies &2(|c| + |d|) + |a| + |b| + |\gamma| + |\delta| - 2|\epsilon| = |\nu|.
  \end{align*}
  So, our sum still depends only on $k$, $l$, $a$, $b$, $\gamma$, and $\delta$ for $n$ large enough.
\end{proof}
\begin{Theorem}{}{}
  $\Hom(\mu, \lambda)$ is a Harish-Chandra bimodule of finite $K$-type.
\end{Theorem}
\begin{proof}
  By \cref{lem:stabilized_nu}, the multiplicity of $V_\nu \in \Hom(\mu^{(n)}, \lambda^{(n)})$ doesn't depend on $n$ for any $\nu \in S$, where $S$ is as in \cref{def:s_filtration}. This implies that $\Hom(\mu, \lambda)$ is a well-defined Harish-Chandra bimodule in $\mathcal{C}_t$. To show it has finite $K$-type, take $k$ so that $\nu \in F^k S$; then
  \begin{align*}
    [\Hom(\mu, \lambda) : V_\nu] &= [F^k \Hom(\mu, \lambda) : V_\nu] \\
                                 &= [F^k \Hom(\mu^{(n)}, \lambda^{(n)}) : V_\nu] \\
                                 &= [\Hom(\mu^{(n)}, \lambda^{(n)}): V_\nu] < \infty.
  \end{align*}
\end{proof}
What do partitions $\mu, \lambda$ look like? We have
\begin{align*}
  \lambda &= (\alpha, \beta, \gamma) \in \mathbb{C}^k \times \mathbb{C}^l \times \mathcal{P} \\
  \mu &= (\alpha, \beta, \gamma) \in \mathbb{C}^k \times \mathbb{C}^l \times \mathcal{P}, a \in \mathbb{Z}^k, b \in \mathbb{Z}^l. \\
\end{align*}
Because the parts of $\alpha^{(n)}, \beta^{(n)}$ tend to infinity, as do their consecutive differences, as we showed in \cref{lem:transcendental_ultraproduct}, we need the parts of $\alpha$ and $\beta$ to be transcendental with transcendental difference. We also need $l(\lambda^{(n)}) \ll n$, so $t$ should be algebraically independent from the parts of $\alpha, \beta$. However, it's not clear if there are more necessary conditions. For now, we have the following condition on $\alpha, \beta$:
\begin{Theorem}{}{}
  Let $\lambda = (\alpha, \beta, \gamma) \in \mathbb{C}^k \times \mathbb{C}^l \times \mathcal{P}$ be a triple, so that $\{\theta_i\}_{i \le l(\theta)} \cup \{1, t\}$ are algebraically independent numbers for $\theta = \alpha, \beta$. Then, we can find $\lambda^{(n)} = [\alpha^{(n)}, \beta^{(n)}, \gamma]$ for almost all $n$ with $l(\lambda^{(n)}) \ll n$ so that
  \begin{align*}
    \alpha_i &= \prod_{\mathcal{F}} \alpha_i^{(n)} \\
    \beta_i &= \prod_{\mathcal{F}} \beta_i^{(n)}
  \end{align*}
  and
  \begin{equation*}
    \lim_{n \to \infty}(\theta_i^{(n)} - \theta_{i + 1}^{(n)}) = \infty, i \le l(\theta), \theta \in \{\alpha, \beta\}.    
  \end{equation*}
  For any $a \in \mathbb{Z}^k, b \in \mathbb{Z}^l, \delta \in \mathcal{P}$, let $\mu \coloneq (\alpha + a, \beta + b, \delta)$. Then,
  \begin{equation*}
    \Hom(\mu, \lambda) = \prod_{\mathcal{F}}^r \Hom(V_{\mu^{(n)}}, V_{\lambda^{(n)}})
  \end{equation*}
  is a Harish-Chandra bimodule in $\mathcal{C}_t$ of finite $K$-type.
\end{Theorem}

This construction gives a nontrivial family of Harish-Chandra bimodules of finite $K$-type by generalizing classical examples of such bimodules and provides some insight into how one can construct manageable ``infinite-dimensional'' objects of Deligne categories. These may not be all the Harish-Chandra bimodules of finite $K$-type. For example, we don't know of formulas for multiplicities of irreducibles of non-polynomial irreducible representations of $GL_n$; such a formula might imply that for arbitrary bipartitions $\lambda$ and $\mu$, $\Hom_{\Rep(GL_t)}(\mu, \lambda)$ is a Harish-Chandra bimodule of finite K-type. Constructing more such bimodules or classifying all such bimodules are directions for future research.


\printbibliography[heading=bibintoc]

\chapter*{Acknowledgements}

First and foremost, I thank my amazing advisors, Michael Hopkins and Pavel Etingof, for suggesting this thesis topic, inspiring me to understand mathematics from a variety of perspectives, and guiding me throughout both the thesis process and my mathematical growth in general. Despite their busy schedules, they met weekly with me to answer my questions, recommend references, and suggest new avenues of thought. I also thank Professor Hopkins for his brilliant and witty teaching in both Math 231A and 231BR. I thank Alexandra Utiralova, my graduate student mentor, for many fun and thoughtful mathematical discussions and journeys.

I thank all the people who proofread my thesis - Professors Hopkins and Etingof, Sasha, Joyce, Lux, Gaurav, Raluca, Lucy, Henry, Alison, Gwyn, Davis, Daniel, and William.

I thank Professor Noam Elkies, who inspired my decision to become a math concentrator through his excellent teaching in Math 55, and who first introduced me to representation theory; Professor Joe Harris, who graciously agreed to do a representation theory reading project with me at PRISE 2018, leading to many conversations about topology, algebraic geometry, and representation theory, and who has advised me since that summer; Morgan Opie, who introduced me to category theory through her 2018 summer tutorial and who has always brightened my day in the math department; and Professor Wilfried Schmid, who taught me a great deal about Lie groups and algebras.

I thank all my friends, who supported me throughout the coronavirus pandemic and without whom my Harvard and Harvard mathematics experiences would be woefully incomplete - in particular, Gwyn, Davis, and Natalia for all the math, games, anime, and random thoughts we shared; Joyce, Alison, Annie, Henry, Franklyn, and Liana for all the emotional support and deep conversations; Crc's for making PRISE an unforgettably good summer; the residents of the math lounge, who made the lounge a vibrant second home; and my roommate Lux, whose many talents and good cheer never cease to amaze.

I thank the students, staff, tutors, and faculty of Adams House and the math department for creating a home for me at Harvard; and I thank URAF for supporting me and welcoming me to PRISE and the summer research village.

Finally, I thank my family for their love and support, for their stalwart belief in me, and for everything.

\textcolor{white}{To L.C., my heart.}

\end{document}